\documentclass[11pt]{article}
\usepackage{graphics , graphicx}
\usepackage{epsfig}
\usepackage{subfigure}
\usepackage{amsmath, amssymb, amsfonts}
\usepackage{a4, a4wide}
\usepackage{lscape}
\usepackage{color}
\usepackage{verbatim}
\usepackage{array}
\def\ie{{\rm i.e.,\/}\ }
\def\etc{{\rm etc.\/}\ }

\newcommand{\nco}{\newcommand}
\nco{\one}{\ensuremath{\,\,\mathrm{l}\!\!\!1}} 
\nco{\ZZ}{\mathbb{Z}}
\nco{\CC}{\mathbb{C}}
\nco{\red}{\color{red}}
\nco{\redend}{\normalcolor}
\nco{\0}{\underline{0}}

\begin{document}

\noindent
\begin{flushright}
{\small CERN-PH-TH-2010-019} \\
\end{flushright}
 
 \vspace{0.2cm} 
 
\begin{center}
\begin{Large}
    {Exceptional quantum subgroups \\  \vspace{0.20cm}  for the rank two Lie algebras $B_2$ and $G_2$}
\end{Large}
\end{center}


\begin{center}
\renewcommand{\thefootnote}{\arabic{footnote}}
\bf{R. Coquereaux} \footnotemark[1] ${ }^{ ,}$\footnotemark[4] ,
\bf{R. Rais} \footnotemark[2] ${ }^{ ,}$\footnotemark[3] ${ }^{ ,}$\footnotemark[4] ,
\bf{E.H. Tahri} \footnotemark[3] ${ }^{ ,}$\footnotemark[4]
\\ \renewcommand{\thefootnote}{\arabic{footnote}} 
\end{center}

 \vspace{0.2cm} 
 
\begin{center}
{January 2010}
\end{center}

 \vspace{0.3cm}

\abstract{Exceptional modular invariants for the Lie algebras $B_2$ (at levels $2,3,7,12$) and $G_2$ (at levels $3,4$) can be obtained from conformal embeddings. We determine the associated algebras of quantum symmetries and discover or recover, as a by-product, the graphs describing exceptional quantum subgroups of type $B_2$ or $G_2$ which encode their module structure over the associated fusion category. Global dimensions are given.}

\vspace{0.3cm}


\vspace{0.2cm}

\noindent{\bf{Classification}}: 81R50; 16W30; 18D10. 

\addtocounter{footnote}{0}
\footnotetext[1]{ {\scriptsize{\it CERN, Geneva, Switzerland}. On leave from {\it CPT, Luminy, Marseille, France}.}}
\footnotetext[2]{ {\scriptsize{\it Centre de Physique ThŽorique (CPT)  Luminy, Marseille, France}.}}
\footnotetext[3]{ {\scriptsize{\it D\'epartement de Physique, Facult\'e des Sciences, Universit\'e Mohamed I,  Oujda, Morocco}.}}
\footnotetext[4]{\scriptsize{coque@cpt.univ-mrs.fr}, \scriptsize{rochdi.rais@gmail.com}, \scriptsize{hassanfa@yahoo.com}}

\vspace{0.4cm}

\unitlength = 1mm

\section{Introduction}

What is written below summarizes what is needed to obtain and study quantum subgroups (self-fusion)  associated with existence of conformal embeddings $G \subset K$.  In this section we do not assume that $G$ coincides with $B_2$ or $G_2$. Warning : several properties reminded below do not necessarily hold in other types of situations like non-simple conformal embeddings followed by contraction, or other types of mechanisms also leading to quantum modules or modular invariants.

To every complex Lie group $G$ and to every positive integer $k$ (the level), one associates a fusion category  that can be described either in terms of integrable representations of affine Lie algebras, or in terms of a particular class of representations of quantum groups at roots of unity \cite{KazhdanLusztig}. Fusion by fundamental integrable representations can be described by graphs called "fusion graphs". These fusion categories somehow generalize, at the quantum level, the theory of representations of Lie groups and Lie algebras, and the process of "fusion" replaces the usual notion of tensor product of representations. 

To every such fusion category is associated a Wess-Zumino-Witten model describing a conformal quantum field  theory of diagonal type, defined on a torus. By ``diagonal type'', one means that the vacuum partition function of the model, written in terms of characters, is a sesquilinear form fully determined by a modular invariant matrix $M$, which, in those cases, is just the identity matrix. 

In classical group theory, groups may have non trivial subgroups and  the space of characters of a subgroup is a module over the ring of characters of the group. In the present situation, we have a fusion ring, specified by the pair $(G,k)$ and one may look for quantum analogs of subgroups (or modules) by looking at ``module-categories'' over which the given fusion category acts \cite{FuchsSchweigert, KirilovOstrik, Ostrik}. To every such quantum module is also associated a partition function, which is modular invariant but not diagonal, and that physically describes a boundary conformal field theory. Two such quantum modules may be associated to the same partition function,  the classification issues are therefore different. The module structure is entirely determined by the action of the generators of the fusion ring on the chosen module, and this can be encoded by matrices with non negative integer entries, or, equivalently, by one or several graphs.

A last ingredient to be introduced is the concept of algebra of quantum symmetries or double triangle algebra DTA \cite{Ocneanu:paths} : 
given a module-category ${\mathcal E}_k$ over a fusion category  ${\mathcal A}_k(G)$, one may define its endomorphism category  ${\mathcal O(E)} = End_{{\mathcal A}} {\mathcal E}$. In practice one prefers to think in terms of (Grothendieck)  algebras and modules, and use the same notations to denote them:  given a quantum module ${\mathcal E}_k$ over the fusion algebra ${\mathcal A}_k(G)$, one may define an algebra ${\mathcal O(E)}$ ``of quantum symmetries'' which  acts on ${\mathcal E}$ in a way compatible with the ${\mathcal A}$ action. 
The fundamental objects of ${\mathcal A}_k(G)$, \ie to those fundamental representations of $G$ that exist at the given level, also act on ${\mathcal O(E)}$, and in particular on its simple objects,  in two possible ways (the later is a bimodule over the former), so that one can associate to each of them two matrices with non negative entries describing left and right multiplication, and therefore two graphs, called left and right chiral graphs whose union is a non connected graph called the Ocneanu graph. The chiral graphs associated with one or another fundamental object of ${\mathcal A}_k(G)$ can themselves be disconnected. More generally, the bimodule structure of  ${\mathcal O(E)}$ is described by ``toric matrices'' that can be physically interpreted as twisted partition functions (presence of defects, see \cite{PetkovaZuber:Oc}). 

The classification of $A_1=su(2)$ modular invariant partition functions (which coincides with the classification of $A_1$ quantum modules) was performed more that twenty years ago  \cite{CIZ} and was shown to be in one to one correspondence with the $ADE$ Coxeter-Dynkin diagrams (\ie with the simply laced complex Lie algebras).  The classification of $A_2=su(3)$ modular invariant partition functions was obtained by \cite{Gannon},  and the graphs describing the action of the fundamental representations of $SU(3)$ on the associated modules were discovered by \cite{DiFrancescoZuber}, see also \cite{YellowBook, GilDahmaneHassan}  (the classification was slightly amended in \cite{Ocneanu:Bariloche}). In the case of $A_3=su(4)$, the classification of quantum modules (and corresponding graphs) was obtained at the end of the nineties and presented in \cite{Ocneanu:Bariloche}.
For every single choice of $G$ we have the diagonal theories ${\mathcal A}_k(G)$ themselves (with ${\mathcal A}_k(SU(2)) = A_{k+1}$) but we also meet several non diagonal infinite series, usually denoted ${\mathcal D}$, like
 the $D_{even}$ or $D_{odd}$ members of the $SU(2)$ family,  and also a finite number of exceptional cases, like the $E_6$, $E_7$ and $E_8$ examples of the same $G=SU(2)$ family.
  
Conformal embeddings of a  Lie group $G$ into simple Lie groups $K$\footnote{They are listed in references \cite{BaisBouwknegt, KacWakimoto, SchellekensWarner}}  give rise to modular invariants that are sums of squares (called type\footnote{When the modular invariant is not block diagonal, it is of type $II$} $I$). This happens for instance for the $D_4$, $E_6$ and $E_8$ invariants of the $SU(2)$ family, or for the ${\mathcal D}_3$, ${\mathcal E}_5$, ${\mathcal E}_9$ and ${\mathcal E}_{21}$ invariants of the $SU(3)$ family.
These cases share several nice features: first of all, from the study of the connected parts of the chiral graphs,  one can obtain one or several quantum modules for ${\mathcal A}_k$; moreover the module $\mathcal{E}_k$ associated with the connected component of the identity, in chiral graphs, can be endowed with an associative multiplication with non-negative structure constants, compatible with the action of the fusion algebra (such modules are said to have self-fusion\footnote{We use this terminology only when ${\mathcal E}_k$ is a subalgebra of its algebra of quantum symmetries} and define ``quantum subgroups'' ). In other words, the choice of a conformal embedding  singles out a quantum subgroup (called exceptional), and sometimes, like for $SU(3)$ at level $9$, an exceptional quantum module.  

 When $G$ is semi-simple, conformal embeddings followed by contraction\footnote{or ``contracted embeddings''} with respect to a simple factor lead to modular invariants that may or may not be of type $I$. For instance, when $G=SU(2)$, the cases $E_7$ and $D$ (even or odd) can be obtained from non-simple  conformal embeddings followed by contraction, but $E_7$ and $D_{odd}$ do not enjoy self-fusion whereas the $D_{even}$ do; notice that  the lowest member ($k=4$) of the $D$ series with self-fusion can be obtained both from direct (\ie non followed by contraction) conformal embedding, namely $SU(2)_4 \subset SU(3)_1$,  or from a contracted embedding $SU(2)_4 \times SU(4)_2 \subset SU(8)_1$.
  
   Quantum symmetries (DTA) for quantum subgroups of type $G=SU(2)$ were first obtained in \cite{Ocneanu:Unpublished} and \cite{Ocneanu:paths}, see also \cite{Gil:thesis,PetkovaZuber:Oc}. For $G=SU(3)$ they were also first obtained in \cite{Ocneanu:Unpublished}, see also \cite{OujdaRDGH, RobertGilSL3Categories, EstebanGil}.  Quantum symmetries for exceptional quantum subgroups of $G=SU(4)$ are described in \cite{RobertGil:E4SU4, RobertGil:E4E6E8ofSU4}.
Exceptional quantum subgroups of the rank two Lie algebras $G=B_2 \sim C_2$ and $G=G_2$ occur at levels $k=2,3,7,12$ for the former, and $k=3,4$ for the later. 
 The associated modular invariants themselves were obtained long ago by \cite{ChristeRavanini}, \cite{Verstegen:B2G2}, and some fusion graphs for $G_2$ (in the case of $A_2 \sim SU(3)$ see also  \cite{DiFrancescoZuber}) were obtained first in \cite{DiFrancesco:thesis}, by imposing various types of constraints, in particular eigenvalues constraints, on their adjacency matrices. 
 Our purpose, in this paper, is to determine  their algebras of quantum symmetries (this was not done before), or equivalently their Ocneanu graphs, and to deduce the fusion graphs from the later.  For instance, in the $G_2$ case,  fusion graphs  already obtained in \cite{DiFrancesco:thesis} (figures $27$ and $35$)  correspond to  subgraphs with blue lines of the Ocneanu graphs given in figures \ref{fig:fusiongraphsE3G2} and \ref{fig:fusiongraphsE4G2} of the the present paper.   These algebras are determined by solving the so-called modular splitting equation that expresses associativity of their bi-module structures over the fusion algebras ${\mathcal A}_k(G)$. Solving the full set of modular splitting equations -- what we do here -- gives full information, not only, as a by-product, the fusion graph of the quantum subgroup\footnote{\label{misleadingnotation}Although standard practice, we shall usually refrain from using the notation $G_k$ to denote ${\mathcal A}_k(G)$.}, but the collection of toric matrices (twisted partition functions for boundary conformal field theory with defects), the quantum symmetries, and, in some cases,  exceptional quantum modules associated with the same modular invariant (they appear at levels $3,7$ for $B_2$ and at levels $3,4$ for $G_2$). As far as we know, the fusion graphs for quantum subgroups and quantum modules of $B_2$ did not appear before.

There exist quantum modules and quantum subgroups of type $B_2$ and $G_2$ associated with conformal embeddings of non-simple type followed by contraction, however we do not call them ``exceptional'' (by definition), and therefore we do not study them here. The case ${\mathcal E}_2(B_2)$ is  special, since, once constructed as an exceptional quantum subgroup of $B_2$ at level $2$, it can be recognized as an orbifold ${\mathcal D}_2(Spin(5))$ \ie as ``a structure of type ${\mathcal D}$'',  a kind of higher analog of $D_4={\mathcal D}_4(SU(2))$, ${\mathcal D}_3(SU(3))$, or ${\mathcal D}_2(SU(4))$.
There is another partition function for $G_2$ at level $4$  found in \cite{Verstegen:B2G2},  which is of type $II$,  and should give rise to a module\footnote{Its fusion graph is denoted $\mathcal{E}_8^\star$ in \cite{DiFrancesco:thesis}: the subscript $k$ (level) is shifted by the dual Coxeter number.} denoted  ${\mathcal E}^\star_4(G_2)$ 
isomorphic to a module ${\mathcal E}^\star_3(F_4)$; its algebra of quantum symmetries is actually isomorphic\footnote{Equality of fusion matrices $N_{01}(G_2) \simeq N_{1000}(F_4)$ and  $N_{10}(G_2) \simeq N_{0001}(F_4)$ at respective levels $4$ and $3$ can be seen from tables of fusion graphs and quantum groupo\"\i d structures given in \cite{Coquereaux:FusionGraphs}.} with ${\mathcal A}_4(G_2) \simeq {\mathcal A}_3(F_4) $ but we shall not study it here because it is not obtained from a conformal embedding.

The size of calculations involved in the study of any single example is huge, and many interesting results will have to be put aside. For instance it would be unreasonable to present the  $81$ distinct  twisted partition functions (matrices of size $91 \times 91$) associated with the $144$ quantum symmetries of ${\mathcal E}_{12}(B_2)$. However fusion graphs and Ocneanu graphs are reasonably small,  so that they can be presented and commented. More details will be given in a thesis \cite{RochdiThesis}.
 
\section{Terminology, notations and techniques}

\subsection{General methods}

Irreducible representations of the Lie group $G$ are denoted $m,n\ldots$ The level of an irreducible representation $n$ of highest weight $\lambda_n$ is the integer $\langle \lambda_n,  \theta  \rangle$ where  $\theta$ is the highest root of $Lie G$, and $\langle .. , .. \rangle $ is the fundamental quadratic form. These representations appear\footnote{They are often called integrable, in relation with the theory of affine Lie algebras.}  as simple objects of  ${\mathcal A}_k={\mathcal A}_k(G) $, if and only if their level is smaller or equal to $k$. 
 This condition selects a finite set, with cardinal $r_A$, of irreducible representations $m,n\ldots$ of $G$ that are parametrized by $r$-tuples in the basis of fundamental weights ($r$ is the rank of $G$);  the objects $m$ are ordered by increasing values of the level, and, at a given level,  the $r$-tuples are sorted in a standard way, starting from the end (see examples later). One says that a representation of $G$ ``exists at level $k$'' if its level is smaller than $k$ or equal to $k$. Fundamental representations of $G$ correspond to the fundamental weights. In the case $G=SU(N)$ all fundamental representations already appear at level $1$. This is usually not so for other choices of $G$.
We call ``basic representations'' those fundamental representations of $G$ that have smallest classical dimension. 
Quantum dimension\footnote{or ``categorical dimension''. We shall often drop the adjective and just write ``dimension''.} $qdim(n)$ of a simple object $n$ can be calculated from the quantum version of the Weyl formula, together with the choice of a root of unity $q=exp(i \pi/\kappa)$, with an altitude $\kappa=k+g$, where $g$ is the dual Coxeter number of $G$. The order of ${\mathcal A}_k$ (or global dimension) is $\vert {\mathcal A}_k \vert = \sum_n (qdim(n))^2$.

Products of simple objects of ${\mathcal A}_k$ read $m \, n = \sum N_{mnp} \, p$, they are described by fusion coefficients and fusion matrices $N_m$ with $(N_m)_{np} = N_{mnp}$ of dimension $r_A \times r_A$.
Fusion coefficients depend on the level $k$, they can be obtained from the Verlinde formula (sec. \ref{Verlinde}). 
To use this expression, one needs to know the $r_A \times r_A$ matrix representing $S$,  one of the two generators ($S$ and $T$) of the modular group $SL(2,\ZZ)$ in a vector space of dimension $r_A$; it is calculated from the Kac-Peterson formula \cite{Kac-Peter-81, Kac-Peterson} (see sec. \ref{KacPeterson}),  or from character polynomials, and the calculation involves in any case a summation over the Weyl group of $G$, which is in general tedious. Using Verlinde formula to determine all fusion coefficients turns out to be time consuming; it is more efficient to use it to calculate only the fundamental fusion matrices $N_f$ (those associated with fundamental representations $f$ existing at level $k$) and to obtain the other fusion matrices $N_n$ from the fusion polynomials that express them polynomially, with rational coefficients, in terms of the fundamental ones (introduce a vector $u= (u_1,\ldots,u_{r_A})$ with  $r_A$ components, write for each fundamental $f$ the family  of equations  $u_f \, \times \, (u) = N_f \, . \, (u)$ and solve for $u_n$ in terms of the $u_f$'s). At this point, one can recover the quantum dimension of simple objects $n$ by calculating the normalized Perron Frobenius vector associated with the basic fundamental fusion matrix (its Perron Frobenius norm is denoted $\beta$), or from the entries of the matrix $S$, and also recover the quantity $\vert {\mathcal A}_k \vert$.

From the chosen conformal embedding $(G,k) \subset (K,1)$, one determines, by reduction, a non trivial modular invariant $M$: a matrix of dimension $r_A\times r_A$, with non negative integer entries and $M_{11}=1$, commuting with $S$ and $T$. From that piece of data, one reads: $r_E=Tr(M)$, the number of simple objects of ${\mathcal E}_k$, $r_O=Tr(M  M^{tr})$ the number of simple objects of ${\mathcal O(E)}$, and $r_W = \#( \{ (m,n) / M_{m,n} \neq 0) \} $, the rank of the family of toric matrices $W_{x{\0}} $ (see below), an integer  smaller or equal to $r_O$. One can also calculate quantities like $\vert {\mathcal A}_k/{\mathcal E}_k \vert$ or $\vert {\mathcal E}_k \vert$ that are quantum analogues of cardinalities for discrete subgroups or cosets.

Call $x,y,\ldots$ the generators associated with simple objects of ${\mathcal O(E)}$. Its bimodule structure  over ${\mathcal A}_k$ is encoded by toric matrices $W_{xy}$ defined by $m \, x \, n = \sum_y \, (W_{xy})_{mn} \, y$, or by matrices $V_{mn} $ defined by $(V_{mn})_{xy}= (W_{xy})_{mn}$.   In general, existence of a conjugation (a ``bar operation'' $m \mapsto \overline m$) sends $W_{x{\0}} $ to $W_{{\0} x}$, but we shall not have to worry about that here because this operation is trivial for $B_2$ and $G_2$, so we set $W_x=W_{x{\0}} =W_{{\0} x}$.
Imposing bimodule associativity $(m_1 m_2) x (n_1 n_2) =  m_1 (m_2  x n_1) n_2$ and taking\footnote{Our indices often refer to some chosen ordering (starting with $1$) on the collection of simple objects of  ${\mathcal A}$, ${\mathcal E}$ or ${\mathcal O}$, but we also denote their unit by $\underline{0}$, so that, for instance, $N_{\{0,0\}}$,  $N_1$ and $N_{\underline 0}$ denote the same identity matrix $\one$. This will be clear from the context and should not confuse the reader.}  $x={\0}$ leads to the modular splitting equation\footnote{This equation was first displayed in\cite{Ocneanu:Bariloche}, the derivation presented above was found in \cite{Gil:thesis}.} 
$$N_{m}\,M\,N_{n}^{tr} = \sum_{x} (W_{{\0} x})_{mn} \, W_{x{\0}} $$ 
where we used the fact that $W_{{\0}{\0}}=M$.  For every choice of the pair $(m,n)$ (there are $r_A^2$ choices) the left hand side ${\mathcal K}_{mn}$, itself a matrix $r_A \times r_A$, is known. The (often arduous) problem is to determine the $W_{x{\0}} $ matrices.  One uses the property that, for each choice of $(m,n)$,  the coefficients $(W_{{\0} x})_{mn}$ on the r.h.s. of the above equations are strongly constrained by the fact that $ ({\mathcal K} _{m n})_{\overline{m} \overline{n}}  =  \sum_{p,q} (N_{m})_{\overline{m} p}\, ({M})_{pq}  \,(N_{n}^{tr})_{{q} \overline{n}}$ is a calculable integer equal to a sum of squares measuring the norm $\sum_x  |(W_{{\0} x})_{mn}|^2$ of a vector $\sum_x \, (W_{{\0}x})_{mn} \, x \in {\mathcal O(E)}$,  its components being themselves non-negative integers.   The toric matrices $W_{x{\0}} $ are therefore determined from the matrices ${\mathcal K}_{mn}$ associated with vectors of increasing norms $1,2,3\ldots$  More details about this technique can be found in \cite{RobertGil:E4E6E8ofSU4} and  \cite{EstebanGil}.  One difficulty, apart from the size of matrices involved, comes from the fact that the rank $r_W$ of the family of  matrices $W_{x{\0}} $ is usually smaller than $r_O$. 

Once the $W_{x{\0}} $ matrices are determined, one notices that another particular case of the (bi) associativity constraints gives intertwining equations 
$$  
N_f \, W_{x{\0}}  \,  {N_{\0}}^{tr} = \sum_y \,( V_{f{\0}}  )_{xy}\, W_{y{\0}}   \qquad \hbox{and} \qquad   N_{\0} \, W_{x{\0}}  \,  {N_f}^{tr} = \sum_y \,( V_{{\0} f} )_{xy}\, W_{y{\0}}    
$$ 
This gives a set of linear equations for the $V_{f{\0}} $.
A priori the solution could be non unique, but  it turns out that it is so, up to isomorphism (at least in the studied cases), after imposing that matrix elements should be non-negative integers.
So, for every fundamental object $f$ of ${\mathcal A}_k$ one obtains  two matrices $V_{f{\0}} $ and $V_{{\0} f}$ called chiral left and chiral right generators. They are adjacency matrices for the various graphs that build the Ocneanu graph of quantum symmetries.

The fusion graph of ${\mathcal E}_k$ is obtained from any one of the two chiral parts of the later (for instance left) by requiring that its vertices, to be the simple objects $a,b,\ldots$ of  ${\mathcal E}_k$,  should span a space of dimension $r_E$ stable under the  action of  ${\mathcal A}_k$. The Cayley graph(s) for this action is (are) deduced from a study\footnote{${\mathcal E}_k$  contains the identity $\underline{0}$, but in some cases one can discover other quantum modules (not subgroups).} of the connected components of the chiral generators $V_{f{\0}} $.  As a by-product we obtain the adjacency matrices $F_f$ describing this action. 
More generally, the action of simple objects of  ${\mathcal A}_k$ on simple objects $a,b,\ldots$ of   ${\mathcal E}_k$ reads $n \, a = \sum F_{nab} \, b$ and is characterized by the so-called annular matrices $F_n$ with $(F_n)_{ab}= F_{nab}$ of dimension $r_E \times r_E$. Since they provide a representation of the fusion algebra, the easiest way to obtain all of them, once we have the fundamental ones, $F_f$, is to use the already determined fusion polynomials $u_n$ since they also express $F_n$ in terms of the $F_f$'s. We introduce the (rectangular) essential matrices $E_a$ with $(E_a)_{nb} =F_{nab}$ of dimension $r_A \times r_E$.  

The matrix $E_{\underline{0}}$, where $\underline{0}$ is the unit object of  ${\mathcal E}_k$, plays a special role and is often called ``the intertwiner''.  One can then obtain quantum dimensions for all simple objects of ${\mathcal E}_k$, either by calculating the normalized Perron Frobenius eigenvector of the basic annular matrix, or by another technique using the induction-restriction rules read from the intertwiner $E_{\underline{0}}$ (see later). As a check, one can in particular calculate  $\vert {\mathcal E}_k \vert = \sum_a (qdim(a))^2$ and recover the result already directly obtained in a previous step.

At this point, there are still several interesting things to do, for instance :  1) obtain explicitly a system of generators $O_x$ for ${\mathcal O(E)}$,  2) Obtain, for each simple object $x$ of ${\mathcal O(E)}$, a matrix $S_x$ that encodes structure constants of the module action of this algebra on the module ${\mathcal E}$, 3) obtain a system of generators $G_a$ for ${\mathcal E}$; remember that in the present situation, the module enjoys self-fusion. This multiplication should be compatible with the module structure: $(n a) b = n ( a b)$. It is not always commutative\footnote{Existence of classical symmetries for the graph ${\mathcal E}_k$ often leads to non-commutativity of the self-fusion and of the algebra of quantum symmetries (nevertheless, chiral left and right parts commute).}, nevertheless $1 \, a = a\, 1$, so that  $(n a) 1 = (n 1) a$ requires the constraint $E_a = E_{\underline{0}} . G_a$. Bi-module associativity  $m(xy)n = (mx)(yn)$ implies in particular the relation $V_{n{\0}}  = \sum_z \, (V_{n{\0}} )_{{\0} z} \, O_z$, so that if the first line of the matrix $V_{f{\0}} $ contains a simple non-zero entry equal to $1$ (it is so for quantum subgroups of type $SU(n), n=2,3,4$), one may define $O_f^L \doteq V_{f{\0}} $. In the same way the relation $n(ab)=(na)b$ implies in particular $F_n=\sum_c \, (F_n)_{\underline{0}c} \, G_c$, and if the first line of 
$F_{f}$ contains a simple non-zero entry equal to $1$, one obtains a  graph matrix $G_f$ equal to $F_f$, but there is no justification for doing so in general\footnote{We thank G. Schieber for this comment.} and we shall actually meet a counter-example in studying  $\mathcal{E}_2(B_2)$.

Simple objects $a,b, \ldots$ of the module-category  ${\mathcal E}_k$ can also be thought  as right modules over a Frobenius  algebra  ${\mathcal F}$, which is a particular object in the monoidal category ${\mathcal A}_k$. This object is not simple in general, and its decomposition in terms of simple objects can be read from the first column of the intertwiner $E_{\underline{0}}$ (the first essential matrix describes induction rules). Actually, for quantum subgroups examples coming from conformal embeddings, this decomposition can also be read from the modular block of identity in the partition function (the modular invariant), and this is simpler.

Finally,  following the lines of  \cite{Ocneanu:paths},  one can associate a quantum groupo\"\i d ${\mathcal B}$ to every module-category ${\mathcal E}_k$. It is  a finite dimensional weak Hopf algebra which is simple and co-semisimple. 
One can think of the algebra ${\mathcal B}$ as a direct sum of $r_A$ matrix simple components,  and of its dual, the algebra $\widehat {\mathcal B}$, as a sum of $r_O$ matrix simple components. 
The dimensions $d_n$ (and $d_x$) of these blocks, called horizontal or vertical dimensions, or dimensions of generalized spaces of essential paths, or spaces of admissible triples or generalized triangles, \etc,  can be obtained from the annular (or dual annular) matrices. We shall not give more details about the structure of this quantum groupo\"\i d in the present paper but, for each of the studied examples,  we shall list the dimensions of its blocks as well as  its total dimension $d_{\mathcal B} = \sum_n d_n^2 = \sum_x d_x^2 $.

 \subsubsection*{Sequence of  steps} 
In order to study a particular quantum subgroup of $G$ at level $k$, determined by a conformal embedding $G \subset K$,   the strategy is the following.

\begin{itemize}
\item Step 1: Determine the $SL(2,\ZZ)$ generators $S$ and $T$, the fusion matrices $N_m$ at level $k$, the quantum dimensions $qdim(n)$ of simple objects of ${\mathcal A}_k(G)$ and the global dimension (quantum order) $\vert {\mathcal A}_k(G) \vert$. Draw the fundamental fusion graph(s).  Determine fusion polynomials for the fusion algebra. 
\item Step 2: Determine the modular invariant $M$ from the embedding $G \subset K$. Find $r_E$, $r_O$, $r_W$ and calculate the global dimension $\vert {\mathcal E}_k\vert $. 
\item Step 3: Solve the modular splitting equation and determine all $r_O$ toric matrices $W_{x{\0}} $ (with a single twist).
This is technically the most involved part of the whole procedure.
\item Step 4: For each fundamental representation of $G$ existing at that level, solve the intertwining equations for the left and right chiral generators $V_{f{\0}} $ and $V_{{\0} f}$ of the algebra of quantum symmetries  ${\mathcal O(E)}$. Draw the graph(s) of quantum symmetries. As a by-product, obtain the graph defining the quantum subgroup ${\mathcal E}$ and its adjacency matrix. In some cases, one also obtains, at this step, an exceptional module.
\item Step 6: Study ${\mathcal E}$. Draw its graph(s). Calculate the quantum dimension of its simple objects. Check the already obtained value of its global dimension.
\item Step 7: Using the fusion polynomials (step 1), obtain the $r_A$ annular matrices $F_n$ and in particular the intertwiner $E_{\underline{0}}$ (tables of induction). Calculate the size of the horizontal blocks (horizontal dimensions) and the dimension of the bialgebra ${\mathcal B}$.
\item Steps 8: Obtain generators $O_x$, dual annular matrices $S_x$ and graph matrices $G_x$.
\end{itemize}

\subsection{Remarks}

 \subsubsection{The chiral modular splitting technique} 
One can replace the modular splitting equation (step 3) by a system of equations  that is technically simpler than the other but gives less information since it determines the adjacency matrices describing the graph(s) of ${\mathcal E}$, as well as the intertwiner $E_{\underline{0}}$ and the annular matrices, but does not provide a description of the algebra of quantum symmetries.
The chiral equation of modular splitting is a simplified form of the full system of equations,  reflecting the fact that the following equality, for $m,n \in {\mathcal A}_k, a \in {\mathcal E}_k$, holds: 
$
(m  (n a) ) = (m . n) a
$.
For cases where ${(F_p)}_{{\0}{\0}}= {M}_{p {\0}}$, a condition that holds for cases associated with conformal embeddings, this associativity constraint implies immediately 
$
\sum_p ({N_n})_{mp} \,  {M}_{p{\0}}  =  \sum_b \, {(F_n)}_{{\0} b} \, {(F_m)}_{b{\0}}  
$.
The left hand side (that involves only the first line of the modular invariant matrix) is known and the right hand side (the annular matrices) can be determined thanks to methods analog to those used to solve the full set of equations (Step 3), but this is technically simpler since the previous identity  describes only $r_A^2$ equations instead of $r_A^4$. 

 \subsubsection{The relative modular splitting technique} 
One way to obtain all toric matrices without having to solve the equation of modular splitting is to use the following technique (it was briefly explained in one section of \cite{OujdaRDGH}). In order to use it, we need not only the adjacency matrices describing the graph of ${\mathcal E}_k$ (this can be obtained in the simplified manner just explained above) but also all graph matrices $G_a$ describing the self-fusion. We introduce a ``relative modular invariant'' $m$ which is a matrix $r_G \times r_G$ such that $E_{\underline{0}} \, m \, {E_{\underline{0}}}^{tr} = M$, as well as a collection of ``relative toric matrices'' defined by the property $a \, x \, b = \sum_y \, (w_{xy})_{ab} \, y$ since, in our cases, ${\mathcal O(E)}$ is trivially an $\mathcal E$ bimodule. Using self-fusion and  bi-associativity $(a \, a') \, x \, (b \, b') = a \, (a' \, x \, b) \, b' $ one obtains a system of equations that read $G_{a}\, m \, G_{b}^{tr} = \sum_{x} (w_{{\0} x})_{ab} \, w_{x{\0}} $. The discussion leading to a determination of the right hand side (the ``relative'' toric matrices $w_{x{\0}} $) is similar to the general one (Step 3), but it is technically simpler because these $w$'s are only of size $r_E \times r_E$. Once they are obtained, one can recover the usual toric matrices from the equality $W_{x{\0}}  = E_{\underline{0}} \, w_{x{\0}} $.  

\subsubsection{About fundamental representations}
 Fundamental representations of finite dimensional Lie groups are ``algebraic generators'' in the sense that 
 all other simple objects of the algebra of representation are factors of tensor products of these distinguished objects (they appear on the right hand side when we tensor multiply fundamental representations and decompose the result into irreducible ones), however the set of these distinguished generators is not a system of generators for the algebra, one needs to add -- infinitely many -- others. However, at a fixed level $k$, the category ${\mathcal A}_k(G)$ contains a finite number of simple objects and they all can be written polynomially in terms of the (integrable) fundamental ones thank's to the fusion polynomials already mentioned; it may even happen (this depends upon the choice of $G$ and $k$) that one can obtain a system of generators smaller than the one provided by the fundamental representations existing at the chosen level. In the case of ${\mathcal O(E)}$ and ${\mathcal E}_k$, when the later enjoys self-fusion (which is the case in this paper), the situation is more subtle, in particular when their graphs possess classical symmetries,  because it may not be possible to express all simple objects in terms of those called  ``fundamental'' : one need to add a few others in order to build a system of generators for the studied algebra(s). 

\subsubsection{About the calculation of $\vert {\mathcal E}_k \vert$ from the modular invariant}

It is convenient to think of $\mathcal{A}_k(G)/\mathcal{E}_k$ as a homogenous space, both discrete and quantum.
Remember that, classically, the dimension of the space of vector valued functions defined on the quotient of a finite group $A$ by a subgroup $E$, and valued in a vector space $a$, can be calculated either trivially as $dim(a)\times \vert A/E \vert$ or non trivially, when $a$ is a representation space for $E$,  by decomposing this space of functions  into a sum of irreducible representations of $A$; this applies in particular to the case where $a$ is the trivial representation of $E$, so that we obtain a decomposition of the space (actually the algebra) of complex valued functions over $A/E$. 
Here, vertices of the graph $\mathcal{E}_k$ do not only label irreducible objects $a$ of $\mathcal{E}$ but also spaces $\Gamma_a$ of sections of quantum vector bundles  which can be decomposed, using induction, into irreducible objects of $\mathcal{A}_k$.
Define $\Gamma_a = \sum _{n  \uparrow \Gamma_a} \,  n$ where the sum runs over those simple objects $n$ of  ${\mathcal A}_k$ that appear in the column $a$ of the intertwiner matrix $E_{\underline{0}}$. Its dimension $\vert \Gamma_a \vert$ can be calculated as the sum of quantum dimensions $\sum _{n  \uparrow \Gamma_a} \,  qdim(n)$ but also as:  $qdim(a) \times \vert {\mathcal A}_k/{\mathcal E}_k  \vert$.
When $a={\0}$, we get  ${\mathcal F}=\Gamma_{\0} =\sum _{n  \uparrow \Gamma_{\0}}$
 and its dimension reads $\vert {\mathcal F} \vert = \sum _{n  \uparrow \Gamma_{\0}} \,  qdim(n) =  \vert {\mathcal A}_k/{\mathcal E}_k  \vert = \vert{\mathcal A}_k \vert/\vert{\mathcal E}_k \vert$, and ${\mathcal F}$ plays the role of an algebra of functions over a non commutative space. 
Using the fact that, in our cases, the decomposition of ${\mathcal F}$ into simple objects can also be read from the first modular block of the partition function, we obtain directly ${\mathcal E}_k$ from $ {\mathcal A}_k$ and from the first block of the modular invariant. At a later stage we calculate the quantum dimensions for all simple objects of 
$a$ of ${\mathcal E}_k$, so that the equality $ \vert{\mathcal E}_k  \vert= \sum_a (qdim(a))^2$ provides a nice check of the calculations.

Intersection of the chiral subalgebras of ${\mathcal O(E)}$ defines the (commutative) ambichiral subalgebra ${\mathcal J}$. Since\footnote{We repeat that we only deal here with cases obtained from conformal embeddings (type $I$ and self-fusion)} the quantum subgroup ${\mathcal E}_k$ is obtained from the (left for instance) identity component of ${\mathcal O(E)}$, ambichiral vertices ${\mathcal J}$ of the later define a particular subset $J$ of vertices of the former, that we call ``modular'',  which is also a subalgebra of ${\mathcal E}_k$. 
The simple objects of  $ {\mathcal A}_k$ labelling non-zero diagonal entries  of the modular invariant matrix $M$ define the family (possible multiplicity) ${\mathcal E}xp$ of generalized exponents of ${\mathcal E}_k$,  and $J$ introduces a partition of this family, since if $n$ is an exponent, then ${n  \uparrow \Gamma_a}$ for some modular vertex $a$ of  ${\mathcal E}_k$.

From unitarity properties of $S$ one proves the generalized $ADE$ trigonometrical identity $\sum_{m,n}   qdim(m) \,  Z_{m,n} \, qdim(n) = \sum_n  qdim(n)^2 \,   = \vert {\mathcal A}_k \vert $ which, incidentally, can be used to provide a check of the correctness of the modular invariant $Z$. If the later is of type $I$, it is a sum of modular blocks indexed by $J$, so $\vert {\mathcal A}_k \vert =  \sum_{a \in J} \vert \Gamma_a \vert^2$. Using the value of  $ \vert \Gamma_a \vert$ obtained by induction, one finds $\vert {\mathcal A}_k \vert = \sum_{a \in J} qdim(a)^2 \vert {\mathcal A}_k \vert^2 / \vert {\mathcal E}_k \vert^2  = \vert J \vert  \, \vert {\mathcal A}_k \vert^2 / \vert {\mathcal E}_k \vert^2 $ that simplifies to give   $\vert {\mathcal A}_k \vert = \vert {\mathcal E}_k \vert \, \vert {\mathcal E}_k \vert / \vert {\mathcal J} \vert$. The conclusion is that  $\vert {\mathcal J} \vert =  \vert J \vert$, the sum of squares of quantum dimensions of $J$, or ``ambichiral quantum dimension'', is equal to  $\vert {\mathcal E}_k \vert ^2/ \vert {\mathcal A}_k \vert$ and can therefore be obtained without having to calculate independently the quantum dimensions of the simple objects belonging to $J$. This is not only a double check but a triple check since this value  could also have been calculated, at the very beginning, from the conformal embedding $G\subset K$ as the quantum dimension of $K$ at level $1$ (a sum of squares for another root of unity). This family of non trivial identities that can be verified in each case  and provide a powerful check of the calculations. 
  On general grounds one also proves that  $\vert {\mathcal A}_k(G) \vert  = \vert {\mathcal O(E)} \vert$. 
  
  \subsubsection{The Kac-Peterson formulae for  $S$ and $T$ \label{KacPeterson}}
   
\begin{eqnarray*}
S_{mn}= \frac{i^{\Sigma_{+}} \sqrt{\Delta} }{ (g+k)^{r/2} }\left(\sum _{w} \, \epsilon_{w}  \, e^{-\frac{2 i \pi  \langle w(m+\varrho ),n+\varrho \rangle }{g+k}}\right) 
&,  &
T_{mn}=e^{2 i \pi \left[\frac{\langle m+\varrho ,m+\varrho \rangle }{2 (g+k)}-\frac{\langle \varrho ,\varrho \rangle }{2 g}\right]} \; \delta_{mn}
\end{eqnarray*}
where $g$ is the dual Coxeter number,  $k$ is the level, $w$ runs over the Weyl group of $G$, $\epsilon_{w}$ is its signature, $r$ is the rank of $G$, $\varrho$ is the Weyl vector, $\Sigma_{+}$ is the number of positive roots (also equal to the sum of exponents, or to $r \, \gamma/2$ where $\gamma$ is the Coxeter number),  and $\Delta$ is the determinant of the fundamental quadratic form. With these definitions one has $(S\,T)^3 = S^2 = \mathcal{C}$, the ``charge matrix'' satisfying $\mathcal{C}^2=\one$.
 It may be the right place to remember that $T$ is related as follows to other group theoretical quantities:
 The  eigenvalue of the quadratic Casimir for a representation $\lambda$ is $C(\lambda)=\langle \lambda ,\lambda +2 \varrho \rangle$, the classical Dynkin index is $dim(\lambda)\, C(\lambda)/2 d$, where $d$ is $dim(G)$, and, at level $k$, the conformal weight of $\lambda$ is\footnote{The coefficient $2$ stands for $\langle \lambda ,\lambda \rangle$, the usual convention for long roots.} $h(\lambda)=C(\lambda)/2(k+g)$. Some authors (like \cite{BakalovKirillov} ) prefer to use another $t$-matrix that differs from the one introduced previously by a modular phase: $t = T \, exp(2 i \pi c/24)$ where $c$ is the central charge $c=\frac{d \, k}{g+k}$.
 
    \label{Verlinde}
    For the convenience of the reader; we give the Verlinde formula expressing fusion matrices in terms of the modular generator $S$  :
    $ (N_m)_{np} = \sum_q \, \frac{S_{mq} S_{nq} S_{pq}^\star}{S_{1q}}$
    
\section{Quantum subgroups of $B_2$}

\subsection{General properties of $B_2$ (or  $C_2$)}

$B_2 = Lie(Spin(5))$ has rank $r= 2$, dual Coxeter number $g=3$, Coxeter number $\gamma = 4$, dimension $d= r(\gamma+1) =10$, adjacency matrix (Dynkin) $G = \{\{0, 2\}, \{1, 0\}\}$ (it is not simply laced), 
Cartan matrix $A = 2 * \one - G = \{\{2,-2\},\{-1,2\}\}$, quadratic form matrix $Q =\{\{1,1/2\},\{1/2,1/2\}\} = A^{-1}.\{\{1,0\},\{0,1/2\}\}$, highest root\footnote{All components are written in the base of fundamental weights (Dynkin).} $\theta = \{0, 2\}$,  Weyl vector $\varrho = \{1,1\}$, exponents $(\epsilon_1=1,\epsilon_2=3)$, Casimir degrees $(\epsilon_1+1=2,\epsilon_2+1=4)$, Weyl group order $(\epsilon_1+1)(2,\epsilon_2+1)=8=2^2 \, 2!$, number of positive roots ${\Sigma_{+}} =4$,  $\Delta =1/4$.
The central charge of $B_2$ at level $k$ is $c=\frac{d \, k}{g+k}=\frac{10 k}{k+3}$.
The level $\langle \lambda,  \theta  \rangle$ of an irreducible representation with highest weight $\{\lambda_1, \lambda_2\}$ is $\lambda_1+\lambda_2$ and its (classical or quantum) dimension is given by the classical or quantum version of the Weyl formula and reads
$$ dim[\{\lambda_1, \lambda_2\}] = 
\frac{\left(\lambda _1+1\right)_q \left(\frac{\lambda
   _2}{2}+\frac{1}{2}\right)_q \left(\lambda _1+\frac{\lambda
   _2}{2}+\frac{3}{2}\right)_q \left(\lambda _1+\lambda
   _2+2\right)_q}{\left(\frac{1}{2}\right)_q 1_q
   \left(\frac{3}{2}\right)_q 2_q}
   $$
   where $(x)_q=\frac{q^x - q^{-x}} {q - q^{-1}}$ and $(x)_q=x$ if $q=1$.
   At level $k$, one takes $q = exp(i \pi/\kappa)$ with  $\kappa=k+g=k+3$.
   The  fundamental representations are the vectorial, with highest weight $\{1,0\}$, classical dimension $5$,  quantum dimension $\frac{\left(\frac{5}{2}\right)_q 3_q}{1_q \left(\frac{3}{2}\right)_q}$, and the spinorial (basic), with  highest weight  $\{1,0\}$, classical dimension $4$, quantum dimension $\frac{1_q 3_q}{\left(\frac{1}{2}\right)_q \left(\frac{3}{2}\right)_q}$.
   
   \noindent
Properties of $C_2$ are obtained by exchanging the Dynkin labels $\lambda_1$ and $\lambda_2$ of $B_2$.

\subsection{$B_2$ at level $1$}

Both fundamental representations of $B_2$ already appear at level $1$. The category ${\mathcal A}_1={\mathcal A}_1(B_2)$ has three simple objects ($r_A=3$) corresponding to the trivial, spinorial (basic) and vectorial representations, that we order as $(\{0,0\},\{0,1\},\{1,0\})$. Their quantum dimensions are respectively $(1, \sqrt 2, 2)$ so that $\vert{\mathcal A}_1\vert=4$. The conformal weights modulo $1$ are $(0, 5/16, 1/2)$ and the central charge is $c=5/2$. 

Modular generators :
$$
\begin{array}{ccc}
S/S_{{\0}{\0}}=
 \left(
\begin{array}{ccc}
 1 & \sqrt{2} & 1 \\
 \sqrt{2} & 0 & -\sqrt{2} \\
 1 & -\sqrt{2} & 1
\end{array}
\right) & , & T=\left(
\begin{array}{ccc}
 e^{-\frac{5 i \pi }{24}} & 0 & 0 \\
 0 & e^{\frac{5 i \pi }{12}} & 0 \\
 0 & 0 & e^{\frac{19 i \pi }{24}}
\end{array}
\right)
\end{array}
$$
with $S_{{\0}{\0}} = 1/\sqrt{\vert{\mathcal A}_1\vert} = 1/2$.

\begin{figure}
\centerline{\scalebox{0.8}{\includegraphics[width=6cm]{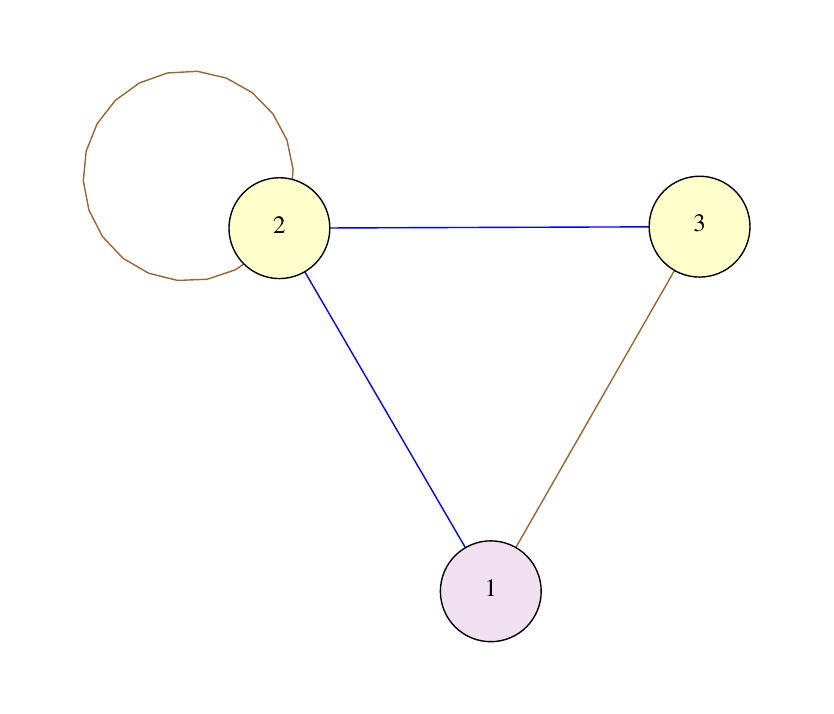}}} 
\caption{The fusion graph of  $\mathcal{A}_1(B_2))$.}
\label{fig:B2level1}
\end{figure}

The graph \ref{fig:B2level1}  describes the fundamental fusion matrices $N_{\{0,1\}}$ and $N_{\{1,0\}}$.
Here and below, labels $1,2,3\ldots$ of vertices refer to the order already defined on the set of simple objects,  multiplication by the generator $01$, labelled $2$, the spinorial,  is encoded by unoriented blue edges, multiplication by $10$, labelled $3$, the vectorial, is encoded  by unoriented brown edges.
 Notice that $N_{\{0,1\}}$ corresponding to the spinorial, is connected, whereas $N_{\{1,0\}}$ corresponding to the vectorial, has two components. Also, the first does not contains self-loops (tadpoles), the second does. These features are generic (\ie for all values of the level $k$). These Cayley graphs describe multiplication by the two generators. There is an intertwiner from $\{0,1\} \otimes n$ to $m$ (resp. from $\{1,0\} \otimes n$) if and only if we have a blue (resp. brown) edge from $n$ to $m$. Edges are not oriented because conjugation is trivial.
The blocks of the associated quantum groupo\"\i d are of size $d_n=(3,4,3)$ and its dimension is $d_{\mathcal B}=34=2^1 \, 17^1$. 

\subsection{$B_2$ at level $2$ and its exceptional quantum subgroup ${\mathcal E}_2$}

\subsubsection*{$B_2$ at level $2$}
The category ${\mathcal A}_2={\mathcal A}_2(B_2)$ has six simple objects ($r_A=6$).
The following table gives  the chosen ordering for the highest weights, the quantum dimensions, the conformal weights modulo $1$, and the size of  blocks for the associated quantum groupo\"\i d.

$$ 
\begin{array}{cccccc}
 \{0,0\} & \{0,1\} & \{1,0\} & \{0,2\} & \{1,1\} & \{2,0\} \\
 1 & \sqrt{5} & 2 & 2 & \sqrt{5} & 1 \\
 0 & \frac{1}{4} & \frac{2}{5} & \frac{3}{5} & \frac{3}{4} & 0 \\
 6 & 12 & 11 & 11 & 12 & 6
\end{array}
$$

\begin{figure}
\centerline{\scalebox{0.8}{\includegraphics[width=6cm]{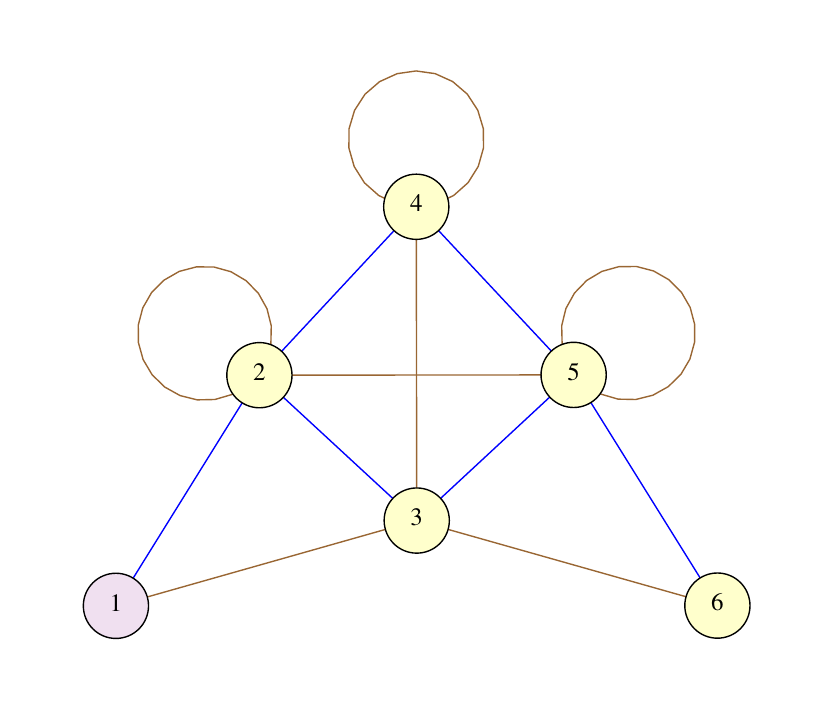}}} 
\caption{The fusion graph of  $\mathcal{A}_2(B_2))$.}
\label{fig:B2level2}
\end{figure}

 The central charge is  $c=4$, the bialgebra dimension is $d_{\mathcal B}=602=2^1\, 7^1 \ 43^1$, and $\vert{\mathcal A}_2\vert=20$.
 The  modular generators, with $S_{{\0}{\0}} = 1/\sqrt{\vert{\mathcal A}_1\vert}= \frac{1}{2 \sqrt{5}}$, are as follows :
\begin{eqnarray*}
S/S_{{\0}{\0}} &=&
\tiny{
\left(
\begin{array}{cccccc}
 1 & \sqrt{5} & 2 & 2 & \sqrt{5} & 1 \\
 \sqrt{5} & \sqrt{5} & 0 & 0 & -\sqrt{5} & -\sqrt{5} \\
 2 & 0 & -1+\sqrt{5} & -1-\sqrt{5} & 0 & 2 \\
 2 & 0 & -1-\sqrt{5} & -1+\sqrt{5} & 0 & 2 \\
 \sqrt{5} & -\sqrt{5} & 0 & 0 & \sqrt{5} & -\sqrt{5} \\
 1 & -\sqrt{5} & 2 & 2 & -\sqrt{5} & 1
\end{array}\right)
}
\\
T&=&diag
\left\{e^{-\frac{i \pi }{3}},e^{\frac{i \pi }{6}},e^{\frac{7 i \pi }{15}},e^{\frac{13 i \pi }{15}},e^{-\frac{5 i \pi
   }{6}},e^{-\frac{i \pi }{3}}\right\}
\end{eqnarray*}
The graph \ref{fig:B2level2} describes the fundamental fusion matrices $N_{\{0,1\}}$ and $N_{\{1,0\}}$.

\subsubsection*{The exceptional quantum subgroup ${\mathcal E}_2={\mathcal D}_2$}

We shall provide more details here than for the cases relative to higher levels, studied later, thank's to the fact that most results for $k=2$ are small enough to be displayed in printed form. Also, we shall discover, at the end of the section, that this case could have been defined directly as an orbifold of $B_2$ at level $2$, so that this exceptional example can also be considered as the smallest member ${\mathcal D}_2$ of a ${\mathcal D}$ family with self-fusion, like the  $D_4$ for $SU(2)$ (an orbifold of $A_5={\mathcal A}_4(SU(2))$,  the ${\mathcal D}_3(SU(3))$ or the $ {\mathcal D}_2(SU(4))$.

\bigskip

We use  the  embedding  $B_2 \subset A_4$ (with Dynkin index $k=2$), which is conformal when the level of $B_2 \simeq Spin(5)$ is $k=2$ and the level of $A_4\simeq SU(5)$ is $1$ (same conformal charges $c$).
The integrable representations of $A_4=Lie(SU(5))$ at level $1$ are the trivial representation and the four fundamental representations of $SU(5)$. A priori the type $I$ partition function obtained by reduction  will have five modular blocks.
The above five simple objects of $A_4$ have all a (quantum) dimension equal to $1$ (using $q=exp(i\pi/(5+1)$ for ${\mathcal A}_1(A_4)$), so that summing their square gives the  ambichiral dimension $\vert J \vert = 5$.
Their conformal weights are ${0,2/5,3/5,3/5,2/5}$. Comparing them with those obtained (modulo $1$) for $B_2$ at level $2$ (\ie using $q=exp(i\pi/(3+2)$) gives a necessary condition for branching rules. One obtains in this way one exceptional $B_2$ partition function (quadratic form\footnote{To alleviate notations we denote  characters by the corresponding classical highest weight.}):
$$Z=(\{0,0\}+\{2,0\})^2  + 2 \{0,2\}^2+2 \{1,0\}^2$$
One verifies that the modular invariant matrix $M$ defined by $Z = \sum_{m,n} m \, M_{mn} \, n$ indeed commutes with $S$ and $T$.
The simple objects appearing in the first modular block of $Z$  define a Frobenius  algebra  ${\mathcal F} =
\{0,0\} \oplus \{2,0\}$, with $\vert{\mathcal F}\vert = \vert{\mathcal A}_2(B_2)\vert / \vert{\mathcal E}_2\vert = qdim(\{0,0\} )+qdim(\{2,0\} )=1+1=2$. Since $\vert{\mathcal A}_2\vert =20$, we find $\vert{\mathcal E}_2\vert =10$. The general relation 
$\vert{\mathcal A}_k(G)\vert  = (\vert{\mathcal E}_k\vert)^2 / \vert J \vert $ leads again to the result  $\vert J \vert =5$.
The number of simple objects for $\vert{\mathcal A}_2$ is $r_A=6$.
From the modular invariant matrix we read the generalized exponents 
$\{\{0,0\}, \{2,0\}; \{0,2\}_2 ;  \{1,0\}_2\}$
and obtain $r_E=6$ (number of simple objects $a$ of the quantum subgroup), 
$r_O=12$ (number of quantum symmetries $x$), $r_W=6$ (number of independent toric matrices $W_x$). Notice in this special case the unusual coincidence $r_E = r_A$, in particular the intertwiner $E_{\underline{0}}$ is here a square matrix.

\begin{figure}[htp]
  \begin{center}
    \subfigure[The left and right generators $\{0,1\}$]{\label{fig:OLR01E2B2}\includegraphics[width=6cm]{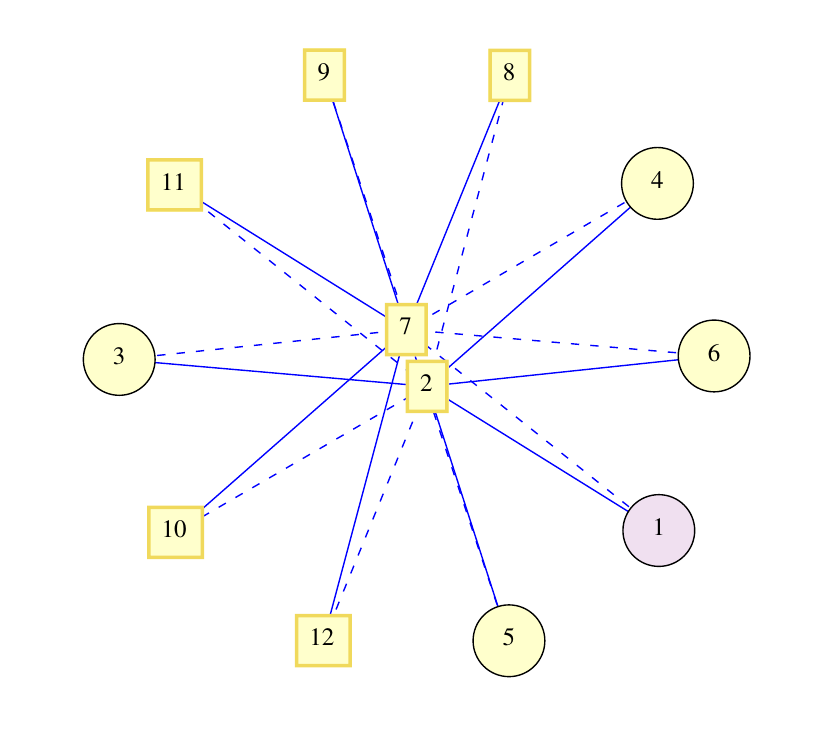}}
    \subfigure[The left and right generators $\{1,0\}$]{\label{fig:OLR10E2B2}\includegraphics[width=6cm]{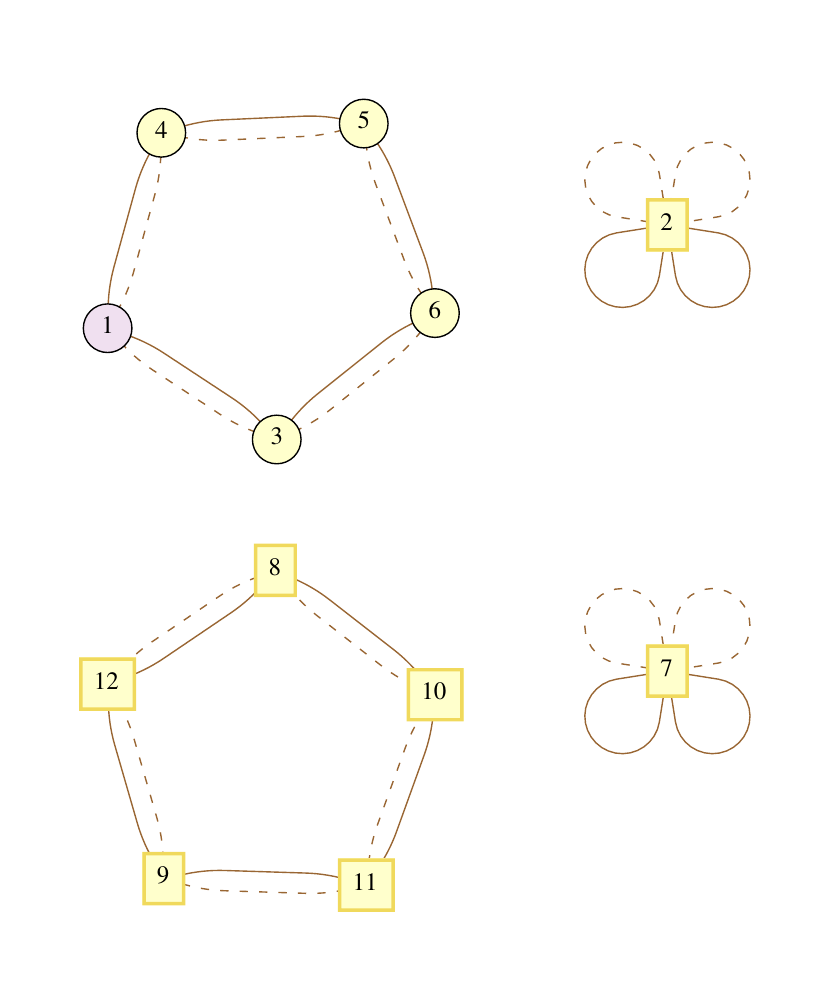}}
  \end{center}
  \caption{Quantum symmetries of  $\mathcal{E}_2(B_2)$}
\end{figure}

Resolution of the modular splitting equation, leading to a list of toric matrices $W_x$ parametrized by simple objects $x$ of  ${\mathcal O(E)}$ represents, as usual, the most involved part of the calculation. 
There are $5$ possible norms associated with matrices ${\mathcal K}_{mn}$. Three toric matrices are discovered while analyzing norm $1$ (among them\footnote{Warning: The indices of $W_x$ are shifted by $1$ and chosen to be compatible with figures  \ref{fig:OLR01E2B2}, \ref{fig:OLR10E2B2}.} $W_1=M$), then two others in norm $2$ (but each of them should appears twice in the list) and a last one in norm $5$ (but it appears five times in the list). So we indeed recover $r_W=6$ and $r_O=12$.  These matrices, given below, define twisted partition functions (only the first is modular invariant) $Z_x = \sum_{m,n} m \, (W_x)_{mn} \, n$ :
{\tiny
$$
W_1=
\left(
\begin{array}{cccccc}
 1 & 0 & 0 & 0 & 0 & 1 \\
 0 & 0 & 0 & 0 & 0 & 0 \\
 0 & 0 & 2 & 0 & 0 & 0 \\
 0 & 0 & 0 & 2 & 0 & 0 \\
 0 & 0 & 0 & 0 & 0 & 0 \\
 1 & 0 & 0 & 0 & 0 & 1
\end{array}
\right),
W_2=
\left(
\begin{array}{cccccc}
 0 & 1 & 0 & 0 & 1 & 0 \\
 0 & 0 & 0 & 0 & 0 & 0 \\
 0 & 2 & 0 & 0 & 2 & 0 \\
 0 & 2 & 0 & 0 & 2 & 0 \\
 0 & 0 & 0 & 0 & 0 & 0 \\
 0 & 1 & 0 & 0 & 1 & 0
\end{array}
\right),
W_7 = W_2^{tr}
$$
$$
W_3=W_4=\left(
\begin{array}{cccccc}
 0 & 0 & 1 & 0 & 0 & 0 \\
 0 & 0 & 0 & 0 & 0 & 0 \\
 1 & 0 & 0 & 1 & 0 & 1 \\
 0 & 0 & 1 & 1 & 0 & 0 \\
 0 & 0 & 0 & 0 & 0 & 0 \\
 0 & 0 & 1 & 0 & 0 & 0
\end{array}
\right),
W_5=W_6=\left(
\begin{array}{cccccc}
 0 & 0 & 0 & 1 & 0 & 0 \\
 0 & 0 & 0 & 0 & 0 & 0 \\
 0 & 0 & 1 & 1 & 0 & 0 \\
 1 & 0 & 1 & 0 & 0 & 1 \\
 0 & 0 & 0 & 0 & 0 & 0 \\
 0 & 0 & 0 & 1 & 0 & 0
\end{array}
\right),
W_8=...=W_{12}\left(
\begin{array}{cccccc}
 0 & 0 & 0 & 0 & 0 & 0 \\
 0 & 1 & 0 & 0 & 1 & 0 \\
 0 & 0 & 0 & 0 & 0 & 0 \\
 0 & 0 & 0 & 0 & 0 & 0 \\
 0 & 1 & 0 & 0 & 1 & 0 \\
 0 & 0 & 0 & 0 & 0 & 0
\end{array}
\right)
$$
}

Generators for the quantum symmetries ${\mathcal O(E)}$ are then obtained from the toric matrices by solving the set of intertwining equations stemming from the general modular splitting equation (see our general discussion).
To ease the reading of the figures \ref{fig:OLR01E2B2}, \ref{fig:OLR10E2B2}, we display separately the graphs associated with $\{0,1\}$  (the  two chiral generators $V_{01,00}$ (left blue: plain lines) and $V_{00,01}$ (right blue: dashed lines)), and those associated with $\{1,0\}$ (the  two chiral generators $V_{10,00}$  (left brown: plain lines) and $V_{00,10}$ (right brown: dashed lines)).

\begin{figure}
\centerline{\scalebox{0.6}{\includegraphics[width=10cm]{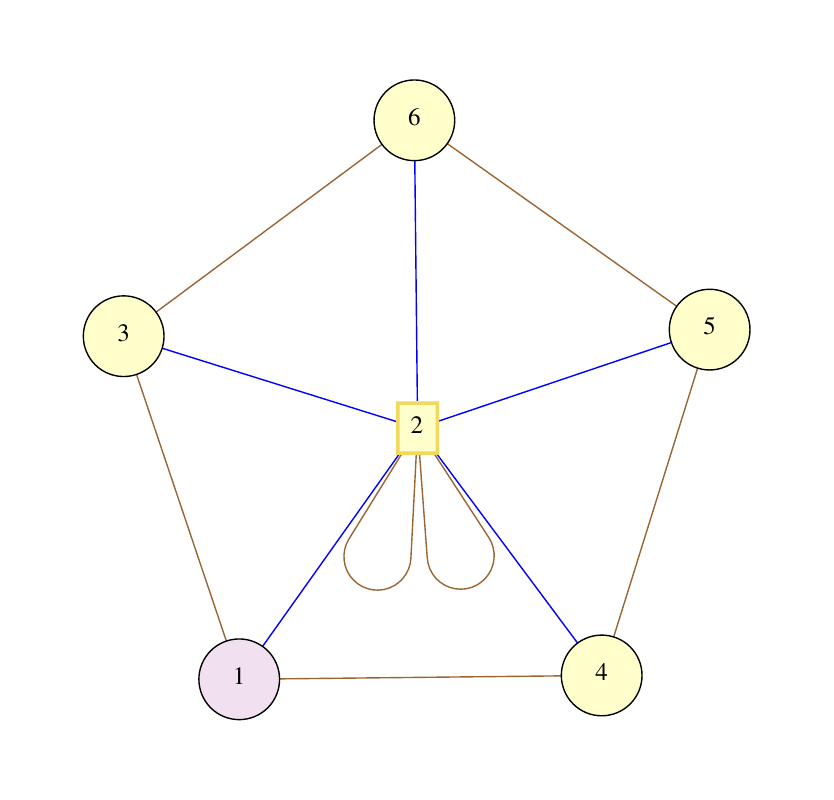}}} 
\caption{The fusion graph of  $\mathcal{E}_2(B_2)$.}
\label{fig:E2B2}
\end{figure}

Fundamental annular matrices $F_{01}$ and $F_{10}$ are then obtained by selecting the component of the identity (for instance from the left graph) in the graph of quantum symmetries.
They are adjacency matrices for the fusion graphs of ${\mathcal E}_2$, given on figure \ref{fig:E2B2}.

The other annular matrices (there is one $F_n$ for every vertex $n$ of the graph $\mathcal{A}_2(B_2)$, so $6$ in our case) can be obtained from the fundamental ones by using fusion polynomials (obtained from the fusion matrices $N_{01}$ and $N_{10}$). Using the ordering adopted for $B_2$ at level $2$, \ie $f_1=F_{0,0}=\one, f_2=F_{01}, f_3=F_{10}, f_4= F_{02}, f_5=F_{1,1}, f_6=F_{20}$, one finds
$
f_4 = -f_1 + f_2.f_2 - f_3;
f_5 = -f_2 + f_2.f_3;
f_6 = -f_3 - f_4 + f_2.f_5
$
which are polynomials in two variables $f_2,f_3$.

From the list of $F_n$'s we obtain the essential matrices $E_a$ and in particular the intertwiner $E_{\underline{0}}$, a rectangular matrix that happens to be square in this special case since $r_A =r_E$.

The numbers $d_n$ of admissible triangles of type $\mathcal{A},\mathcal{E},\mathcal{E}$ (horizontal dimensions) are read from the $F_n$'s,  the blocks of the associated quantum groupo\"\i d are therefore of size $d_n=(6, 10, 12, 12, 10, 6)$ and its dimension (total number of diffusion graphs)  is the sum of their squares: $d_{\mathcal B}=560=2^4 \, 5^1 \, 7^1$. 

\begin{figure}
{\small
$$
\begin{array}{ccc}
F_{01}=
\left(
\begin{array}{cccccc}
 . & 1 & . & . & . & . \\
 1 & . & 1 & 1 & 1 & 1 \\
 . & 1 & . & . & . & . \\
 . & 1 & . & . & . & . \\
 . & 1 & . & . & . & . \\
 . & 1 & . & . & . & .
\end{array}
\right)
&
F_{10}=
\left(
\begin{array}{cccccc}
 . & . & 1 & 1 & . & . \\
 . & 2 & . & . & . & . \\
 1 & . & . & . & . & 1 \\
 1 & . & . & . & 1 & . \\
 . & . & . & 1 & . & 1 \\
 . & . & 1 & . & 1 & .
\end{array}
\right)
&
E_{\underline{0}}=
 \left(
\begin{array}{cccccc}
 1 & . & . & . & . & . \\
 . & 1 & . & . & . & . \\
 . & . & 1 & 1 & . & . \\
 . & . & . & . & 1 & 1 \\
 . & 1 & . & . & . & . \\
 1 & . & . & . & . & .
\end{array}
\right)
\end{array}
$$
}
\caption{Annular matrices and intertwiner for $\mathcal{E}_2(B_2)$.}
\label{fig:AnnularIntertwinerE2B2}
\end{figure}

We know a priori that this quantum module is also a quantum subgroup (self-fusion): one can associate graph matrices $G_a$ with all vertices of the fusion graph $\mathcal{E}_2$, so again $6$ in our case. The vertex denoted $1$ is the unit (\ie $G_1 = \one$). From the fusion graph we read $F_{01}=G_2$, associated with the basic representation, but   $F_{10}=G_3+G_4$. Notice that in this particular example, one of the fundamental annular matrices (namely $F_{10}$)  cannot be identified with any of the $G_a$'s. The matrices $G_a$ are obtained by solving the whole set of equations stemming from the fusion graphs: $F_{01}.G_2 =  G_1+G_3+G_4+G_5+G_6$,  $F_{01}.G_a =  G_a$, for $a\neq 2$, $F_{10}.G_2 = 2 G_2$, $F_{10}.G_4 = G_1+G_5$, \etc
After imposing non negativity of the coefficients and compatibility equations $E_a=E_{\underline{0}}.Ga$,  we find, up to symmetries, only one solution, that we do not present explicitly here (see details in \cite{RochdiThesis}), but it should be noticed that the obtained algebra is non commutative. Diagonalizing the non commutative  $*$-algebra $\mathcal{O(E)}$ gives a direct sum of four blocks isomorphic with $\CC$ and two blocks isomorphic with $M(2,\CC)$, in agreement with the expression of $Z$. 
The dual annular matrices $S_x$ encoding the action of $\mathcal{O(E)}$ on  $\mathcal{E}$ can be expressed in terms of the $G_a$'s, and the vertical dimensions   (admissible triangles of type $\mathcal{O},\mathcal{E},\mathcal{E}$)are  $d_x= \{6, 10, 6, 6, 6, 6, 10, 6, 6, 6, 6, 6\}$;  one recovers $d_B=d_n.d_n=d_x.d_x=560$ (duality) as it should.

The quantum dimensions of the simple objects of ${\mathcal E}_2$ are calculated from the normalized eigenvector of $F_{01}=G_2$,  associated with the basic\footnote{Using $\{1,0\}$ instead gives a largest eigenvalue $2$, compatible with $F_{10}=G_3+G_4$, but the associated eigenspace has dimension $2$.} representation $\{0,1\}$, for the Perron-Frobenius eigenvalue ($\sqrt 5$, as it should). One finds $(1,\sqrt 5,1,1,1,1)$. They can also be obtained from the intertwiner $E_{\underline{0}}$ by adding the dimensions (not their square) of the representations of $B_2$ corresponding to the entries appearing in each column and dividing the result by $\vert {\mathcal F} \vert = \vert  {\mathcal A}_2 / {\mathcal E}_2  \vert = 2$.
By summing squares of the quantum dimensions, one recovers $\vert {\mathcal E}_2 \vert = 10$, a value that we had obtained at the very beginning.  The modular vertices are $\{1,3,4,5,6\}$. By adding the squares of their quantum dimensions we recover the ambichiral dimension $\vert J \vert = 5$ already obtained in two ways at the very beginning,  from the embedding into $SU(4)$ and from the general relation $\vert{\mathcal A}_k(G)\vert  = (\vert{\mathcal E}_k\vert)^2 / \vert J \vert $.

The reader familiar with the structure of the smallest  quantum subgroups  of type ${\mathcal D}$ (with self-fusion) for the families $SU(2)$, $SU(3)$ or $SU(4)$ will have recognized here a direct generalization of what happens in those cases, that can also be obtained from ``small'' conformal embeddings or via an orbifold procedure. Now that we have discovered what the fusion graph of ${\mathcal E}_2(B_2)$ is, we may look at the table given for $B_2$, level $2$, and at its fusion graph (figure \ref{fig:B2level2}); it is then clear that we can indeed obtain the graph(s) of ${\mathcal E}_2$ as an orbifold of the later :  fold the graph(s)  \ref{fig:B2level2} along the vertical line 3-4, identify $1$ and $6$,  $2$ and $5$ and blow up (duplicate) the two fixed points $3$ and $4$.

\subsection{$B_2$ at level $3$ and its exceptional quantum subgroup ${\mathcal E}_3$}

\subsubsection*{$B_2$ at level $3$}
The category ${\mathcal A}_3={\mathcal A}_3(B_2)$ has ten simple objects ($r_A=10$).
The following table gives  the chosen ordering for the highest weights, the quantum dimensions, the conformal weights modulo $1$, and the size of  blocks for the associated quantum groupo\"\i d.
$$ \begin{array}{cccccccccc}
 \{0,0\} & \{0,1\} & \{1,0\} & \{0,2\} & \{1,1\} & \{2,0\} & \{0,3\} & \{1,2\} & \{2,1\} &
   \{3,0\} \\
 1 & 1+\sqrt{3} & 1+\sqrt{3} & 2+\sqrt{3} & 3+\sqrt{3} & 1+\sqrt{3} & 1+\sqrt{3} &
   2+\sqrt{3} & 1+\sqrt{3} & 1 \\
 0 & \frac{5}{24} & \frac{1}{3} & \frac{1}{2} & \frac{5}{8} & \frac{5}{6} & \frac{7}{8} & 0
   & \frac{5}{24} & \frac{1}{2} \\
 10 & 24 & 24 & 32 & 40 & 24 & 24 & 32 & 24 & 10
\end{array} $$

 The central charge is  $c=5$, the bialgebra dimension is $d_{\mathcal B}=6728=2^3 \, 29^2$, and $\vert{\mathcal A}_3\vert=24 \left(2+\sqrt{3}\right)$.

\begin{figure}
\centerline{\scalebox{0.8}{\includegraphics[width=10cm]{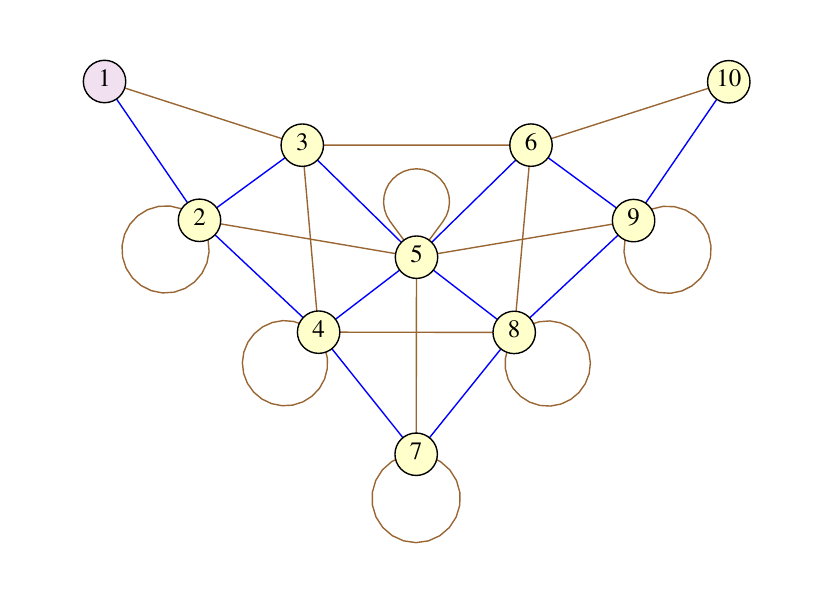}}} 
\caption{The fusion graph of  $\mathcal{A}_3(B_2))$.}
\label{fig:B2level3}
\end{figure}

The modular generators, with $S_{{\0}{\0}} = 1/\sqrt{\vert{\mathcal A}_1\vert}=\frac{1}{12} \left(3-\sqrt{3}\right)$, are as follows (we give $S/S_{{\0}{\0}}$ and $T$):

$$
\begin{array}{c}
\tiny{
\left(
\begin{array}{cccccccccc}
 1 & 1+\sqrt{3} & 1+\sqrt{3} & 2+\sqrt{3} & 3+\sqrt{3} & 1+\sqrt{3} & 1+\sqrt{3} &
   2+\sqrt{3} & 1+\sqrt{3} & 1 \\
 1+\sqrt{3} & 3+\sqrt{3} & 1+\sqrt{3} & 1+\sqrt{3} & 0 & -1-\sqrt{3} & 0 & -1-\sqrt{3} &
   -3-\sqrt{3} & -1-\sqrt{3} \\
 1+\sqrt{3} & 1+\sqrt{3} & 1+\sqrt{3} & -1-\sqrt{3} & 0 & 1+\sqrt{3} & -2
   \left(1+\sqrt{3}\right) & -1-\sqrt{3} & 1+\sqrt{3} & 1+\sqrt{3} \\
 2+\sqrt{3} & 1+\sqrt{3} & -1-\sqrt{3} & 1 & -3-\sqrt{3} & -1-\sqrt{3} & 1+\sqrt{3} & 1 &
   1+\sqrt{3} & 2+\sqrt{3} \\
 3+\sqrt{3} & 0 & 0 & -3-\sqrt{3} & 0 & 0 & 0 & 3+\sqrt{3} & 0 & -3-\sqrt{3} \\
 1+\sqrt{3} & -1-\sqrt{3} & 1+\sqrt{3} & -1-\sqrt{3} & 0 & 1+\sqrt{3} & 2
   \left(1+\sqrt{3}\right) & -1-\sqrt{3} & -1-\sqrt{3} & 1+\sqrt{3} \\
 1+\sqrt{3} & 0 & -2 \left(1+\sqrt{3}\right) & 1+\sqrt{3} & 0 & 2 \left(1+\sqrt{3}\right) &
   0 & -1-\sqrt{3} & 0 & -1-\sqrt{3} \\
 2+\sqrt{3} & -1-\sqrt{3} & -1-\sqrt{3} & 1 & 3+\sqrt{3} & -1-\sqrt{3} & -1-\sqrt{3} & 1 &
   -1-\sqrt{3} & 2+\sqrt{3} \\
 1+\sqrt{3} & -3-\sqrt{3} & 1+\sqrt{3} & 1+\sqrt{3} & 0 & -1-\sqrt{3} & 0 & -1-\sqrt{3} &
   3+\sqrt{3} & -1-\sqrt{3} \\
 1 & -1-\sqrt{3} & 1+\sqrt{3} & 2+\sqrt{3} & -3-\sqrt{3} & 1+\sqrt{3} & -1-\sqrt{3} &
   2+\sqrt{3} & -1-\sqrt{3} & 1
\end{array}
\right)
}
\\
\\
diag
\left\{e^{-\frac{5 i \pi }{12}},1,e^{\frac{i \pi }{4}},e^{\frac{7 i \pi }{12}},e^{\frac{5 i \pi }{6}},e^{-\frac{3 i \pi
   }{4}},e^{-\frac{2 i \pi }{3}},e^{-\frac{5 i \pi }{12}},1,e^{\frac{7 i \pi }{12}}\right\}
   \end{array}
   $$
The graph \ref{fig:B2level3} describes the fundamental fusion matrices $N_{\{0,1\}}$ and $N_{\{1,0\}}$.

\subsubsection*{The exceptional quantum subgroup ${\mathcal E}_3$}

We use  the  embedding  $B_2 \subset D_5$ (with Dynkin index $k=3$), which is conformal when the level of $B_2 \simeq Spin(5)$ is $k=3$ and the level of $D_5\simeq Spin(10)$ is $1$ (same conformal charges $c$). This is the adjoint embedding : $dim(Spin(5))=10$ and the level $k$ coincides with  the dual Coxeter number of  $B_2$.

The integrable representations of $D_5=Lie(Spin(10))$ at level $1$ are the trivial representation, the two spinorial and the vectorial. A priori the type $I$ partition function obtained by reduction will have four modular blocks.
The above four simple objects of $D_5$ have a (quantum) dimension equal to $1$ (using $q=exp(i\pi/(8+1)$ for ${\mathcal A}_1(D_5)$), so that summing their square gives the ambichiral dimension $\vert J \vert = 4$.
Their conformal weights are ${0, 5/8, 5/8, 1/2}$. Comparing them with those obtained (modulo $1$) for $B_2$ at level $3$ (\ie using $q=exp(i\pi/(3+3)$) gives a necessary condition for branching rules. One obtains in this way one exceptional $B_2$ partition function (quadratic form):
$$Z=(\{0,0\}+\{1,2\})^2+2 \{1,1\}^2+(\{0,2\}+\{3,0\})^2 \quad $$
One verifies that the modular invariant matrix $M$ defined by $Z = \sum_{m,n} m \, M_{mn} \, n$ indeed commutes with $S$ and $T$.

The simple objects appearing in the first modular block of $Z$  define a Frobenius  algebra  ${\mathcal F} =
\{0,0\} \oplus \{1,2\}$, with $\vert{\mathcal F}\vert = \vert{\mathcal A}_3(B_2)\vert / \vert{\mathcal E}_3\vert = qdim(\{0,0\} )+qdim(\{1,2\} )=3+\sqrt 3$. Since $\vert{\mathcal A}_3\vert =24(2+\sqrt 3)$, we find $\vert{\mathcal E}_3\vert =4 \left(3+\sqrt{3}\right)$. 
The general relation $\vert{\mathcal A}_k(G)\vert  = (\vert{\mathcal E}_k\vert)^2 / \vert J \vert $ leads again to the result  $\vert J \vert =4$. The number of simple objects for ${\mathcal A}_3(B_2)$ is $r_A=10$.
From the modular invariant matrix
we read the generalized exponents 
$\{\{0,0\}, \{1,2\} ;   \{1,1\}_2 ;  \{0,2\}, \{3,0\}\}$
and  obtain $r_E=6$ (number of simple objects $a$ of the quantum subgroup), 
$r_O=12$ (number of quantum symmetries $x$), $r_W=9$ (number of independent toric matrices $W_x$). 

Resolution of the modular splitting equation, leads to the list of toric matrices $W_x$. 
There are $7$ possible norms associated with matrices ${\mathcal K}_{mn}$. Six toric matrices are discovered while analyzing norm $1$, then two others in norm $2$ (but each of them appears twice) and a last one in norm $3$ (but it appears twice). So we recover $r_W=6$ and $r_O=6+2\times 2+ 2 \times 1=12$.  The nine independent $W_x$'s are given below (the first is the modular invariant).

\smallskip

\noindent
{\tiny
$$
\left(
\begin{array}{cccccccccc}
 1 & . & . & . & . & . & . & 1 & . & . \\
 . & . & . & . & . & . & . & . & . & . \\
 . & . & . & . & . & . & . & . & . & . \\
 . & . & . & 1 & . & . & . & . & . & 1 \\
 . & . & . & . & 2 & . & . & . & . & . \\
 . & . & . & . & . & . & . & . & . & . \\
 . & . & . & . & . & . & . & . & . & . \\
 1 & . & . & . & . & . & . & 1 & . & . \\
 . & . & . & . & . & . & . & . & . & . \\
 . & . & . & 1 & . & . & . & . & . & 1
\end{array}
\right),\left(
\begin{array}{cccccccccc}
 . & 1 & . & . & 1 & . & 1 & . & 1 & . \\
 . & . & . & . & . & . & . & . & . & . \\
 . & . & . & . & . & . & . & . & . & . \\
 . & 1 & . & . & 1 & . & 1 & . & 1 & . \\
 . & . & 2 & 2 & . & 2 & . & 2 & . & . \\
 . & . & . & . & . & . & . & . & . & . \\
 . & . & . & . & . & . & . & . & . & . \\
 . & 1 & . & . & 1 & . & 1 & . & 1 & . \\
 . & . & . & . & . & . & . & . & . & . \\
 . & 1 & . & . & 1 & . & 1 & . & 1 & .
\end{array}
\right),\left(
\begin{array}{cccccccccc}
 . & . & 1 & 1 & . & 1 & . & 1 & . & . \\
 . & . & . & . & . & . & . & . & . & . \\
 . & . & . & . & . & . & . & . & . & . \\
 . & . & 1 & 1 & . & 1 & . & 1 & . & . \\
 . & 2 & . & . & 2 & . & 2 & . & 2 & . \\
 . & . & . & . & . & . & . & . & . & . \\
 . & . & . & . & . & . & . & . & . & . \\
 . & . & 1 & 1 & . & 1 & . & 1 & . & . \\
 . & . & . & . & . & . & . & . & . & . \\
 . & . & 1 & 1 & . & 1 & . & 1 & . & .
\end{array}
\right)
$$},
{\tiny
$$
\left(
\begin{array}{cccccccccc}
 . & . & . & 1 & . & . & . & . & . & 1 \\
 . & . & . & . & . & . & . & . & . & . \\
 . & . & . & . & . & . & . & . & . & . \\
 1 & . & . & . & . & . & . & 1 & . & . \\
 . & . & . & . & 2 & . & . & . & . & . \\
 . & . & . & . & . & . & . & . & . & . \\
 . & . & . & . & . & . & . & . & . & . \\
 . & . & . & 1 & . & . & . & . & . & 1 \\
 . & . & . & . & . & . & . & . & . & . \\
 1 & . & . & . & . & . & . & 1 & . & .
\end{array}
\right),
\left(
\begin{array}{cccccccccc}
 . & . & . & . & 1 & . & . & . & . & . \\
 . & . & . & . & . & . & . & . & . & . \\
 . & . & . & . & . & . & . & . & . & . \\
 . & . & . & . & 1 & . & . & . & . & . \\
 1 & . & . & 1 & . & . & . & 1 & . & 1 \\
 . & . & . & . & . & . & . & . & . & . \\
 . & . & . & . & . & . & . & . & . & . \\
 . & . & . & . & 1 & . & . & . & . & . \\
 . & . & . & . & . & . & . & . & . & . \\
 . & . & . & . & 1 & . & . & . & . & .
\end{array}
\right),\left(
\begin{array}{cccccccccc}
 . & . & . & . & . & . & . & . & . & . \\
 1 & . & . & 1 & . & . & . & 1 & . & 1 \\
 . & . & . & . & 2 & . & . & . & . & . \\
 . & . & . & . & 2 & . & . & . & . & . \\
 1 & . & . & 1 & . & . & . & 1 & . & 1 \\
 . & . & . & . & 2 & . & . & . & . & . \\
 1 & . & . & 1 & . & . & . & 1 & . & 1 \\
 . & . & . & . & 2 & . & . & . & . & . \\
 1 & . & . & 1 & . & . & . & 1 & . & 1 \\
 . & . & . & . & . & . & . & . & . & .
\end{array}
\right)$$,
}
{\tiny
$$\left(
\begin{array}{cccccccccc}
 . & . & . & . & . & . & . & . & . & . \\
 . & 1 & . & . & 1 & . & 1 & . & 1 & . \\
 . & . & 1 & 1 & . & 1 & . & 1 & . & . \\
 . & . & 1 & 1 & . & 1 & . & 1 & . & . \\
 . & 1 & . & . & 1 & . & 1 & . & 1 & . \\
 . & . & 1 & 1 & . & 1 & . & 1 & . & . \\
 . & 1 & . & . & 1 & . & 1 & . & 1 & . \\
 . & . & 1 & 1 & . & 1 & . & 1 & . & . \\
 . & 1 & . & . & 1 & . & 1 & . & 1 & . \\
 . & . & . & . & . & . & . & . & . & .
\end{array}
\right),\left(
\begin{array}{cccccccccc}
 . & . & . & . & . & . & . & . & . & . \\
 . & . & 1 & 1 & . & 1 & . & 1 & . & . \\
 . & 1 & . & . & 1 & . & 1 & . & 1 & . \\
 . & 1 & . & . & 1 & . & 1 & . & 1 & . \\
 . & . & 1 & 1 & . & 1 & . & 1 & . & . \\
 . & 1 & . & . & 1 & . & 1 & . & 1 & . \\
 . & . & 1 & 1 & . & 1 & . & 1 & . & . \\
 . & 1 & . & . & 1 & . & 1 & . & 1 & . \\
 . & . & 1 & 1 & . & 1 & . & 1 & . & . \\
 . & . & . & . & . & . & . & . & . & .
\end{array}
\right),\left(
\begin{array}{cccccccccc}
 . & . & . & . & . & . & . & . & . & . \\
 . & . & . & . & 2 & . & . & . & . & . \\
 1 & . & . & 1 & . & . & . & 1 & . & 1 \\
 1 & . & . & 1 & . & . & . & 1 & . & 1 \\
 . & . & . & . & 2 & . & . & . & . & . \\
 1 & . & . & 1 & . & . & . & 1 & . & 1 \\
 . & . & . & . & 2 & . & . & . & . & . \\
 1 & . & . & 1 & . & . & . & 1 & . & 1 \\
 . & . & . & . & 2 & . & . & . & . & . \\
 . & . & . & . & . & . & . & . & . & .
\end{array}
\right)
$$
}

\begin{figure}[htp]
  \begin{center}
\scalebox{0.5}{
  \scalebox{2.0}{ \subfigure[The left and right generators $\{0,1\}$]{\label{fig:OLR01E3B2}   \includegraphics[width=7cm]{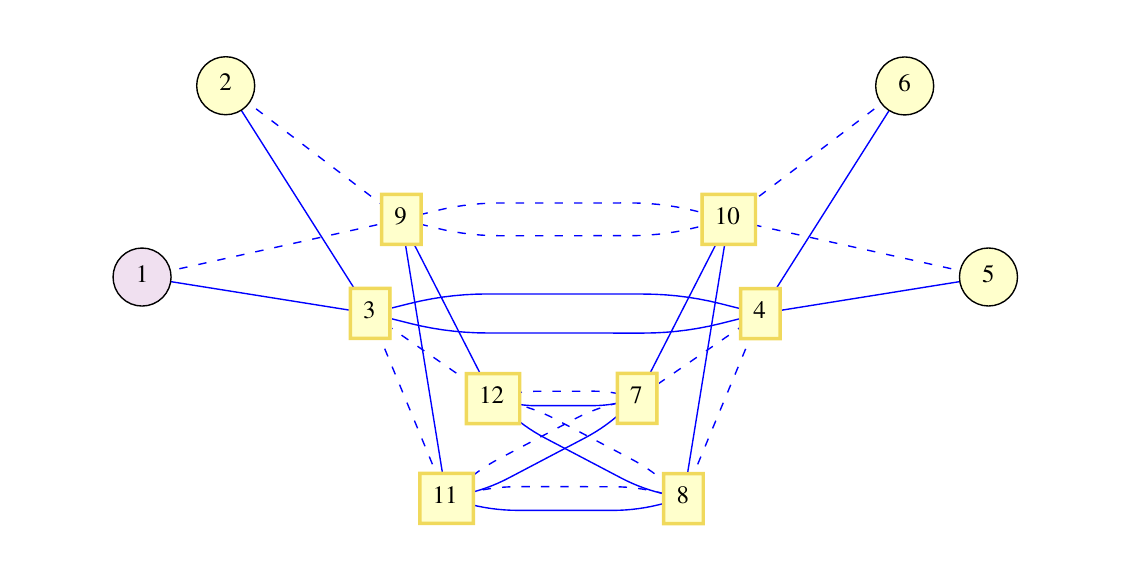}   }}
    \scalebox{2.0}{ \subfigure[The left and right generators $\{1,0\}$]{\label{fig:OLR10E3B2} \includegraphics[width=6cm]{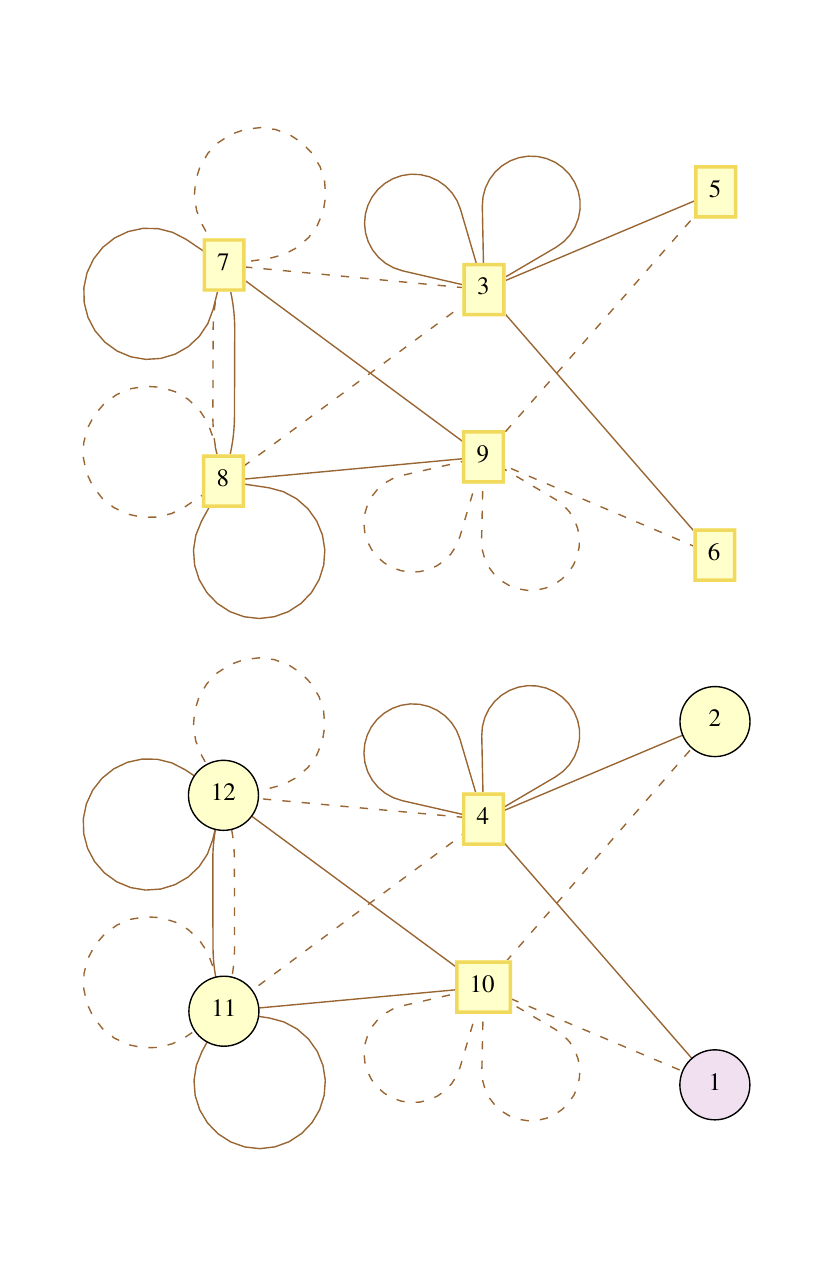}}}
    }
  \end{center}
  \caption{Quantum symmetries of  $\mathcal{E}_3(B_2)$}
  \label{fig:OcneanugraphsE3B2}
\end{figure}

Generators for the quantum symmetries ${\mathcal O(E)}$ are given on  figures \ref{fig:OLR01E3B2}, \ref{fig:OLR10E3B2}. We use the same notations as those defined in the  $\mathcal{E}_2(B_2)$ section.
Blue lines (representation $\{0,1\}$): the unit $1$ moves to the vertex $3$ (resp. $9$) when multiplied by the left (resp. right)  chiral generator.
Brown lines (representation $\{1,0\}$): the unit $1$ moves to the vertex $4$ (resp. $10$) when multiplied by the left (resp. right)  chiral generator.

It is clear from the graph of quantum symmetries that it allows us to define both a quantum subgroup (vertices from $1$ to $6$) and a quantum module (vertices from $7$ to $12$) whose structure is totally different. Both are exceptional modules, but only the first enjoys self-fusion.

\begin{figure}[htp]
  \begin{center}
\scalebox{0.8}{
    \subfigure[$\mathcal{E}_3(B_2)$]{\label{fig:E3B2}\scalebox{1.3}{ \includegraphics[width=6cm]{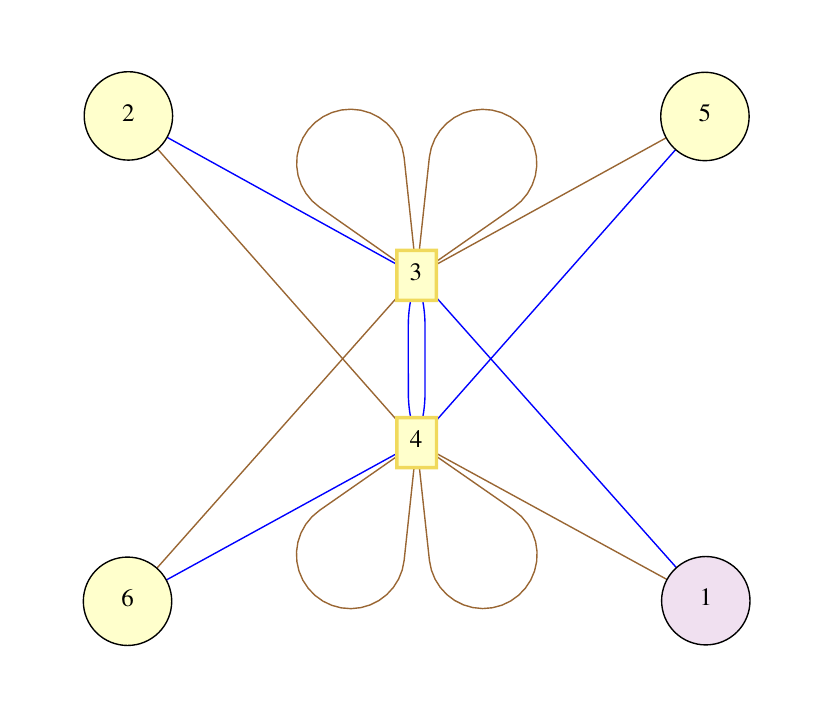}}}
    \subfigure[${\mathcal{E}_3}^M(B_2)$]{\label{fig:E3mB2}\includegraphics[width=6cm]{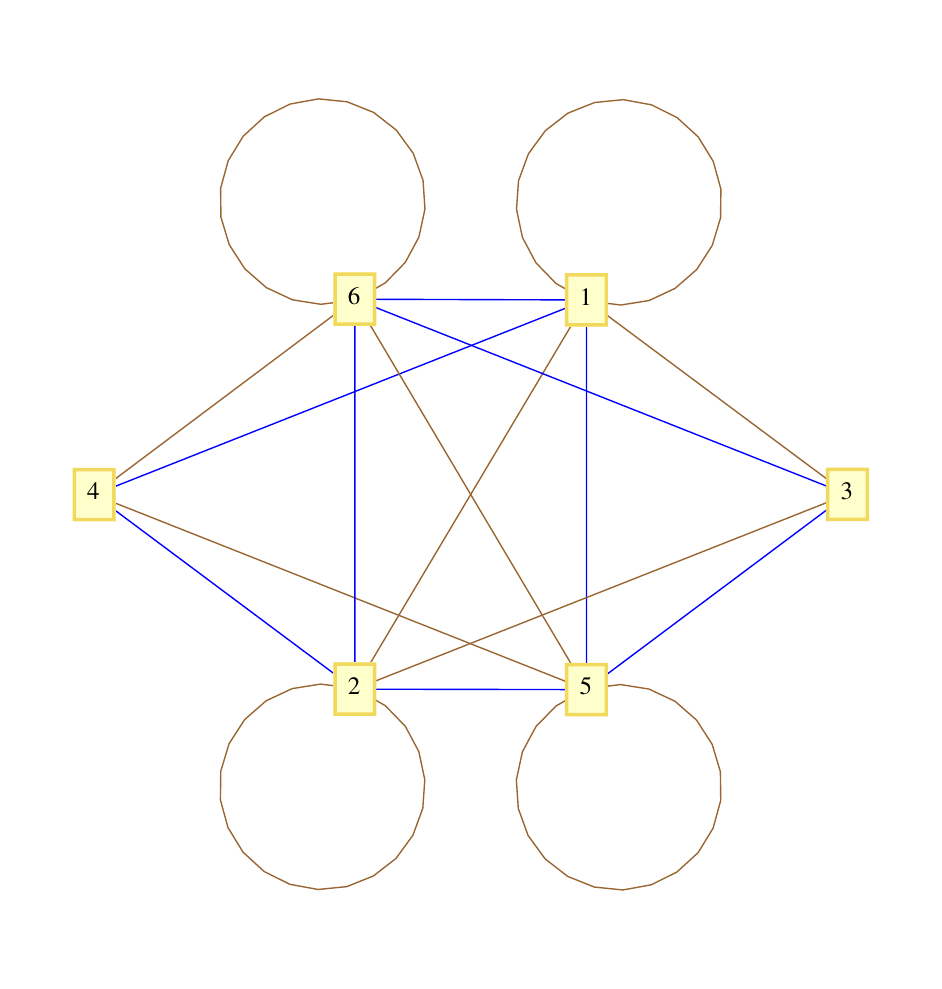}}
    }
  \end{center}
  \caption{Fusion graphs for the quantum subgroup $\mathcal{E}_3(B_2)$ and its module ${\mathcal{E}_3}^M(B_2)$ }
  \label{fig:fusiongraphsE3B2}
\end{figure}

Fundamental annular matrices $F_{01}$ and $F_{10}$ for the quantum subgroup are obtained by selecting the component of the identity (first six vertices, for instance from the left graph).
They are adjacency matrices for the fusion graphs of ${\mathcal E}_3$, given on figure \ref{fig:E3B2}.

In the same way we obtain annular matrices ${F_{01}}^M$ and ${F_{10}}^M$ for the quantum module ${{\mathcal E}_3}^M$; its fusion graph is given on figure \ref{fig:E3mB2}. Warning:  vertex labels $1-6$ are back-shifted on  \ref{fig:E3mB2}, they correspond to vertices $7-12$ of fig.   \ref{fig:OcneanugraphsE3B2}.

The other annular matrices $F_n$ and ${F_n}^M$, for the quantum subgroup and its module are obtained from the fundamental ones by using fusion polynomials (themselves obtained from the fusion matrices $N_{01}$ and $N_{10}$). Using the ordering adopted for $B_2$ at level $3$, one finds: $f_1=1,f_2,f_3,f_2^2-f_3-1,f_2 f_3-f_2,-f_2^2+f_3^2+f_3,f_2^3-2 f_3 f_2-f_2,f_3
   f_2^2-f_2^2-f_3^2-f_3+1,-f_2^3+f_3^2 f_2+f_2,f_3^3+2 f_3^2-2 f_2^2 f_3+f_2^2-1$.
We shall not exhibit these two collections of $10$ matrices  $6\times 6$, but only the two intertwiners $E_{\underline{0}}$ and $E_{{0}}^M$, which are rectangular $10 \times 6$ and describe induction rules (see Fig. \ref{fig:IntertwinersE3B2andE3mB2}).
The numbers $d_n$ (horizontal dimensions) read from the $F_n$'s are $d_n=(6, 12, 12, 18, 24, 12, 12, 18, 12, 6)$ and the dimension of the associated quantum groupo\"\i d  is : $d_{\mathcal B}=d_n.d_n=2016=2^5 \, 3^2 \, 7^1$. 
For the module, we get  $d^M_n=(6, 16, 16, 22, 28, 16, 16, 22, 16, 6)$ and $d^M_{\mathcal B}=3104=2^5 \, 97^1$. Graph matrices $G_a$ describing self-fusion and other details will appear elsewhere \cite{RochdiThesis}.

\begin{figure}
{\small
$$
\begin{array}{cc}
E_{\underline{0}}=
\left(
\begin{array}{cccccc}
 1 & 0 & 0 & 0 & 0 & 0 \\
 0 & 0 & 1 & 0 & 0 & 0 \\
 0 & 0 & 0 & 1 & 0 & 0 \\
 0 & 1 & 0 & 1 & 0 & 0 \\
 0 & 0 & 1 & 0 & 1 & 1 \\
 0 & 0 & 0 & 1 & 0 & 0 \\
 0 & 0 & 1 & 0 & 0 & 0 \\
 1 & 0 & 0 & 1 & 0 & 0 \\
 0 & 0 & 1 & 0 & 0 & 0 \\
 0 & 1 & 0 & 0 & 0 & 0
\end{array}
\right)
&
E^M_{0}=
\left(
\begin{array}{cccccc}
 1 & 0 & 0 & 0 & 0 & 0 \\
 0 & 0 & 0 & 1 & 1 & 1 \\
 1 & 1 & 1 & 0 & 0 & 0 \\
 1 & 2 & 1 & 0 & 0 & 0 \\
 0 & 0 & 0 & 1 & 2 & 2 \\
 1 & 1 & 1 & 0 & 0 & 0 \\
 0 & 0 & 0 & 1 & 1 & 1 \\
 2 & 1 & 1 & 0 & 0 & 0 \\
 0 & 0 & 0 & 1 & 1 & 1 \\
 0 & 1 & 0 & 0 & 0 & 0
\end{array}
\right)
\end{array}
$$
}
\caption{Intertwiners for $\mathcal{E}_3(B_2)$ and ${\mathcal{E}_3}^M(B_2)$ .}
\label{fig:IntertwinersE3B2andE3mB2}
\end{figure}

The quantum dimensions of the simple objects of ${\mathcal E}_3$ are calculated from the normalized eigenvector of $F_{01}$
 associated with the basic representation $\{0,1\}$, for the Perron-Frobenius eigenvalue ($1+\sqrt 3$, as it should)  or from the intertwiner $E_{\underline{0}}$ by adding the dimensions (not their square) of the representations of $B_2$ corresponding to the entries appearing in each column and dividing the result by $\vert {\mathcal F} \vert = \vert  {\mathcal A}_3 / {\mathcal E}_3  \vert =  3 + \sqrt 3$.
 One finds $1, 1, 1 + \sqrt 3, 1 + \sqrt 3, 1, 1$. By summing squares of these quantum dimensions, one recovers the global dimension already obtained at the very beginning
 $\vert {\mathcal E}_3 \vert =4(3+\sqrt 3)$.  The modular vertices are $\{1,2,5,6\}$. By adding the squares of their quantum dimensions we recover the ambichiral dimension $\vert J \vert = 4$ already obtained in two ways at the very beginning,  from the embedding into $Spin(10)$ and from the general relation $\vert{\mathcal A}_k(G)\vert  = (\vert{\mathcal E}_k\vert)^2 / \vert J \vert $.
A similar analysis can be done for the exceptional module ${\mathcal E}_3^M$.

\subsection{$B_2$ at level $7$ and its exceptional quantum subgroup ${\mathcal E}_7$}
\subsubsection*{$B_2$ at level $7$}

\begin{figure}
\centerline{\scalebox{1.0}{\includegraphics[width=10cm]{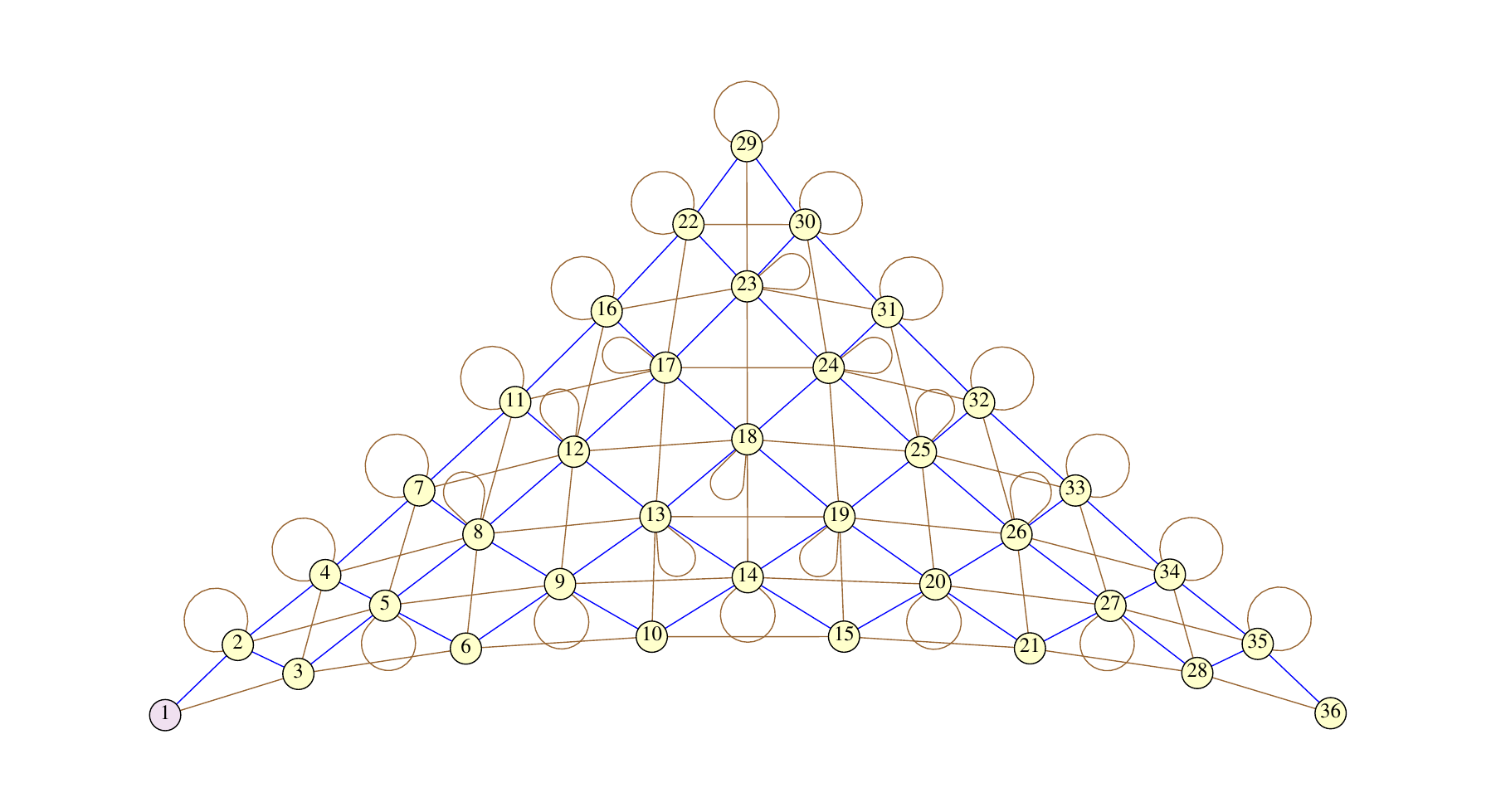}}} 
\caption{The fusion graph of  $\mathcal{A}_7(B_2))$.}
\label{fig:B2level7}
\end{figure}

The category ${\mathcal A}_7={\mathcal A}_7(B_2)$ has $36$ simple objects ($r_A=36$) and it is unreasonable\footnote{\label{webpub} Some tables are given 
in  \cite{Coquereaux:FusionGraphs} } to print lists giving quantum dimensions and conformal weights for all of them, or to give explicitly the matrices $S$ and $T$.
For the fundamental objects one finds
$qdim(\{0,1\}) = \frac{1}{2} \left(1+\sqrt{5}+\sqrt{2 \left(5+\sqrt{5}\right)}\right)$, $qdim(\{1,0\}) =1+\sqrt{5+2 \sqrt{5}} $, 
the central charge is  $c=7$, the bialgebra dimension $d_{\mathcal B}= 7004440 =2^3 \, 5^1 \, 41^1 \, 4271^1$, and  $S_{{\0}{\0}} =  1/\sqrt{\vert{\mathcal A}_k\vert}$.  with $$\vert{\mathcal A}_7\vert= 80 \left(30+13 \sqrt{5}+2 \sqrt{5 \left(85+38 \sqrt{5}\right)}\right)$$
The graph \ref{fig:B2level7} describes the fundamental fusion matrices $N_{\{0,1\}}$ and $N_{\{1,0\}}$. It clearly generalizes those obtained in the previous cases. Notice that  the picture provides an intuitive notion for the notion of level (the ordering chosen for vertices should be clear).

\subsubsection*{The exceptional quantum subgroup ${\mathcal E}_7$}

To keep the size of this paper reasonable, we shall not give many explicit results in this section but only those that somehow summarize them (same thing in the section devoted to a quantum subgroup of $B_2$ at level $12$). 
The comments about techniques are also rather sketchy since they have been already given in previous sections.

\smallskip 

We use  the  embedding  $B_2 \subset D_7$ (with Dynkin index $k=7$), which is conformal when the level of $B_2 \simeq Spin(5)$ is $k=7$ and the level of $D_7\simeq Spin(14)$ is $1$ (same conformal charges $c$).
The integrable representations of $D_5=Lie(Spin(10))$ at level $1$ are the trivial representation, the two spinorial and the vectorial. A priori the type $I$ partition function obtained by reduction will have four modular blocks.
These simple objects of $D_7$ have a (quantum) dimension equal to $1$ (using $q=exp(i\pi/(12+1)$ for ${\mathcal A}_1(D_7)$), so that the ambichiral dimension will be $\vert J \vert = 4$.
Their conformal weights are ${0, 7/8, 7/8, 1/2}$. Comparing them with those obtained (modulo $1$) for $B_2$ at level $7$ (\ie using $q=exp(i\pi/(3+7)$) gives a necessary condition for branching rules. One obtains in this way one exceptional $B_2$ partition function (quadratic form):
$$Z=(\{0,0\}+\{1,6\}+\{2,2\}+\{5,0\})^2+2 (\{1,3\}+\{3,3\})^2+(\{0,6\}+\{2,0\}+\{3,2\}+\{7,0\})^2$$
One verifies that the modular invariant matrix $M$ defined by $Z = \sum_{m,n} m \, M_{mn} \, n$ indeed commutes with $S$ and $T$.

The simple objects appearing in the first modular block of $Z$  define a Frobenius  algebra  ${\mathcal F} =\{0,0\} \oplus \{1,6\}  \oplus \{2,2\}  \oplus \{5,0\}$, with $\vert{\mathcal F}\vert = \vert{\mathcal A}_7(B_2)\vert / \vert{\mathcal E}_7\vert = qdim(\{0,0\} )+qdim(\{1,6\}+qdim(\{2,2\}+qdim(\{5,0\} )= 5 \left(3+\sqrt{5}\right)+\sqrt{250+110 \sqrt{5}}$. Since $\vert{\mathcal A}_7\vert =80 \left(30+13 \sqrt{5}+2 \sqrt{5 \left(85+38 \sqrt{5}\right)}\right)$, we find $\vert{\mathcal E}_7\vert =4 \left(5 \left(3+\sqrt{5}\right)+\sqrt{250+110 \sqrt{5}}\right)$. 
The general relation $\vert{\mathcal A}_k(G)\vert  = (\vert{\mathcal E}_k\vert)^2 / \vert J \vert $ leads again to  $\vert J \vert =4$. The number of simple objects for ${\mathcal A}_7(B_2)$ is $r_A=36$.
From the modular invariant matrix we read the generalized exponents and obtain $r_E=12$ (number of simple objects $a$ of the quantum subgroup), 
$r_O=48$ (number of quantum symmetries $x$), $r_W=36$ (number of independent toric matrices $W_x$). 

Resolution of the modular splitting equation, leads to the list of toric matrices $W_x$. 
There are $44$ possible norms associated with matrices ${\mathcal K}_{mn}$. $16$ independent toric matrices are discovered while analyzing norm $1$, then $12$ others in norm $2$ (among them, four appear twice), then $6$ in norm $3$ (all they appear twice),  and $1$ in norm $4$ (and it appears twice). No information is obtained from analyzing norms $5,6,7$;  the last toric matrix of the list appears in norm $8$ (and it appears twice). At this point we indeed obtain the expected $r_W=36$ independent members of the list and reach the total $r_O=16+(8+2\times 4)+2\times 6+ 2\times 1 + 2 \times 1=48$. The first is, as expected, the modular invariant. These $48$ matrices of size $36\times 36$ are needed to determine the quantum symmetries (next paragraph), but, of course, we shall not give them here. 

Generators for the quantum symmetries ${\mathcal O({\mathcal E})}$ are then obtained. The full graph incorporating both left and right generators, for each of the fundamental objects $\{0,1\}$ and $\{1,0\}$, is nice... but unfortunately a bit clumsy if not properly scaled. It is enough to say that it contains two copies ($2 \times 12$ vertices) of a quantum subgroup, and two copies ($2 \times 12$ vertices) of a quantum module. Both are therefore discovered in this way and will be described below. Fundamental annular matrices for the quantum subgroup and for its module are then obtained as usual, as diagonal sub-blocks of the quantum symmetry generators $V_{\{0,1\},\{0,0\}}$ and  $V_{\{1,0\},\{0,0\}}$.

\begin{figure}[htp]
  \begin{center}
    \subfigure[$\mathcal{E}_7(B_2)$]{\label{fig:E7B2}\scalebox{1.3}{ \includegraphics[width=6cm]{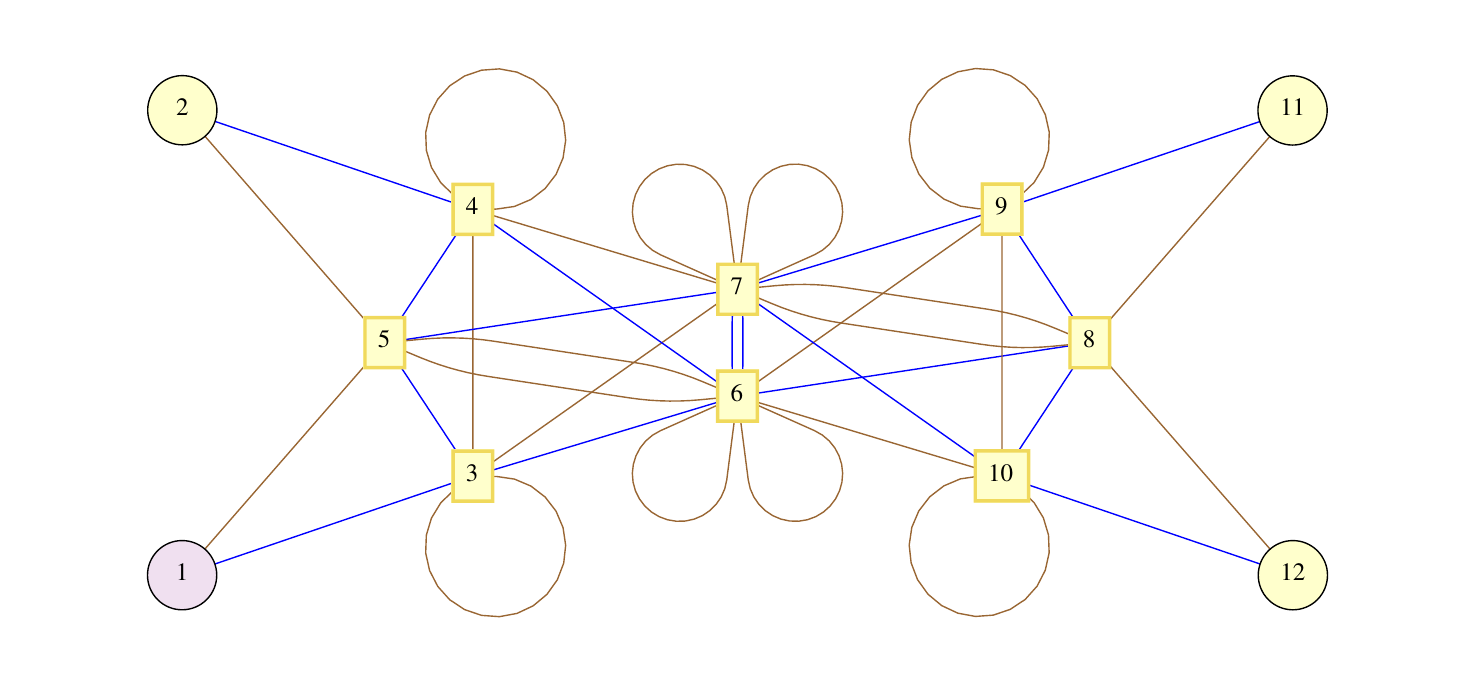}}}
    \subfigure[${\mathcal{E}_7}^M(B_2)$]{\label{fig:E7mB2} \scalebox{1.2}{\includegraphics[width=6cm]{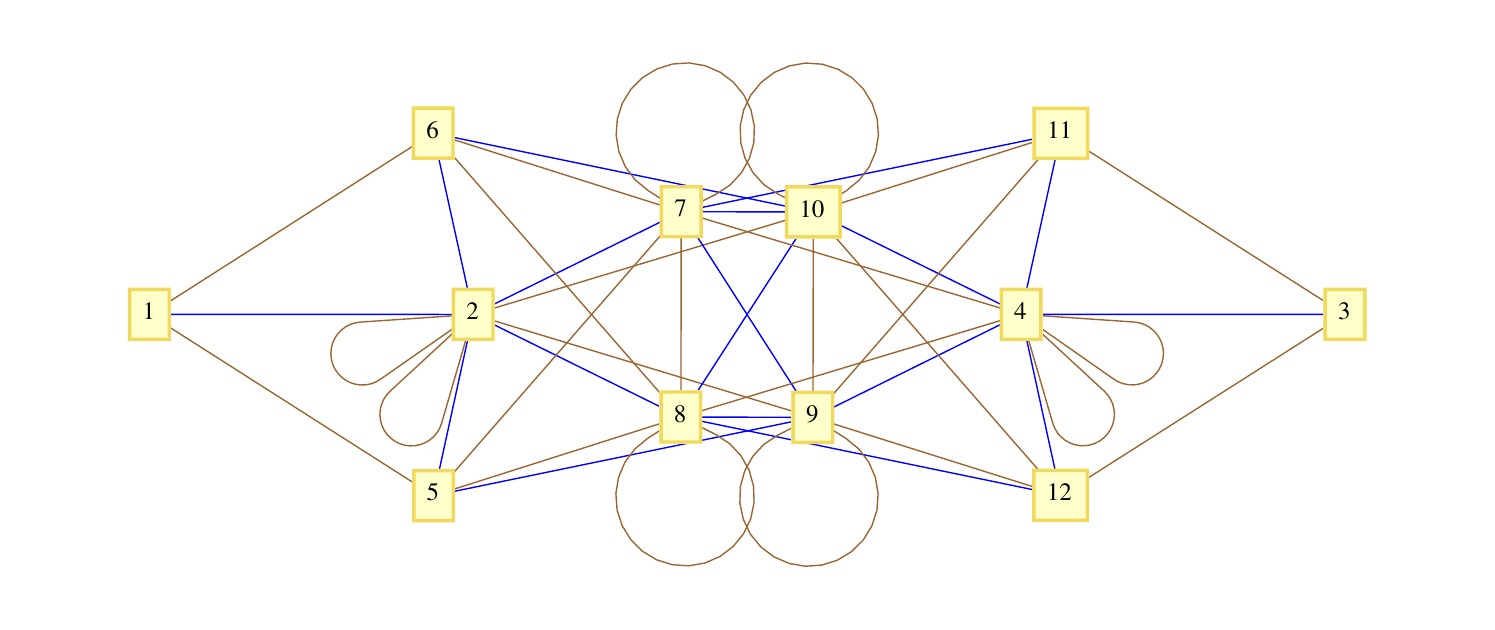}}}
  \end{center}
  \caption{Fusion graphs for the quantum subgroup $\mathcal{E}_7(B_2)$ and its module ${\mathcal{E}_7}^M(B_2)$ }
  \label{fig:fusiongraphsE7B2}
\end{figure}

The fusion graphs of the quantum subgroup ${\mathcal E}_7$, are given on figure \ref{fig:E7B2} and those describing its module  ${{\mathcal E}_7}^M$ appear on figure \ref{fig:E7mB2}. Vertex labels of the later, taken from the graph  ${\mathcal O({\mathcal E})}$, are here back-shifted by $2\times 12$.
As usual the blue color corresponds to $\{0,1\}$ and the brown to \{1,0\}.
${\mathcal E}_7$ enjoys self-fusion (its unit is in position ``1'') and the four modular points of ${\mathcal E}_7$, corresponding to the four ambichiral points obtained as intersection of chiral parts of the graph of quantum symmetries, are located in positions $\{1,2,11,12\}$.

The other annular matrices (both for the quantum subgroup and its module) are obtained from the fundamental ones by using fusion polynomials. We only list below the horizontal dimensions $d_n $ 

$12, 32, 36, 60, 92, 72, 96, 160, 160, 96, 124, 212, 232, 184, 96, 128, 220, 256, $

$232, 160, 72, 108,184, 220, 212, 160, 92, 36, 64, 108, 128, 124, 96, 60, 32, 12,$

\noindent
and the dimension of the quantum groupo\"\i d $d_{\mathcal B}=d_n.d_n=697216=2^7\, \, 13^1 \, 419^1$  associated with $\mathcal{E}_7(B_2)$.

 The intertwiners $E_{\underline{0}}$ (for $\mathcal{E}_7$ and ${\mathcal{E}_7}^M$)  are rectangular matrices of dimension $36 \times 12$, a  bit large to be displayed, but we can condense the information by listing only its non-zero entries (induction\footnote{Entries in positions $1,2,11,12$ coincide with the modular blocks of $Z$} tables). We list $E_{\underline{0}}$ associated with $\mathcal{E}_7$  : 
 
{\tiny

$\{0, 0\} , \{1, 6\} , \{2, 2\} , \{5, 0\}$

 $\{0, 6\} , \{2, 0\} , \{3, 2\} , \{7, 0\}$
 
$ \{0, 1\} , \{0, 7\} , \{1, 3\} , \{1, 5\} , \{2, 1\} , 
   \{2, 3\} , \{2, 5\} , \{3, 1\} , \{4, 1\} , \{5, 1\}$

 $ \{0, 5\} , \{0, 7\} , \{1, 1\} , \{1, 5\} , \{2, 1\} , 
   \{2, 3\} , \{3, 1\} , \{3, 3\} , \{4, 1\} , \{6, 1\}$
   
  $\{0, 6\} , \{1, 0\} , \{1, 2\} , \{1, 4\} , \{1, 6\} , 
   \{2, 2\} , \{2, 4\} , \{3, 0\} , \{3, 2\} , \{4, 0\} , 
   \{4, 2\} , \{6, 0\}$ 
   
  $\{0, 2\} , \{0, 4\} , \{0, 6\} , \{1, 2\} , 2 \{1, 4\} , 
   \{1, 6\} , \{2, 0\} , 2 \{2, 2\} , 2 \{2, 4\} , \{3, 0\} , 
   2 \{3, 2\} , \{3, 4\} , \{4, 0\} , \{4, 2\} , \{5, 0\} , 
   \{5, 2\}$ 
   
  $\{0, 3\} , \{0, 5\} , \{1, 1\} , 2 \{1, 3\} , 2 \{1, 5\} , 
   \{2, 1\} , 2 \{2, 3\} , \{2, 5\} , 2 \{3, 1\} , 2 \{3, 3\} , 
   \{4, 1\} , \{4, 3\} , \{5, 1\}$ 
   
  $\{0, 3\} , \{0, 5\} , \{1, 3\} , \{2, 1\} , 2 \{2, 3\} , 
   \{2, 5\} , \{3, 3\} , \{4, 1\} , \{4, 3\}$ 
   
  $\{0, 4\} , \{1, 2\} , \{1, 4\} , \{2, 2\} , \{2, 4\} , 
   \{3, 2\} , \{3, 4\} , \{4, 2\}$ 
   
  $\{0, 4\} , \{1, 2\} , \{1, 4\} , \{2, 2\} , \{2, 4\} , 
   \{3, 2\} , \{3, 4\} , \{4, 2\}$ 
   
  $\{1, 3\} , \{3, 3\}$ 
  
  $\{1, 3\} , \{3, 3\}$
  
}

The quantum dimensions of the simple objects of ${\mathcal E}_7$ can be calculated, as before,  in several ways.
When ordered as above (with labels $a$ from 1 to 12), one finds a symmetry $qdim[a+1]=qdim[12-a]$ so that it is enough to list the first six dimensions:

{\tiny

$
\left\{1,1,\frac{1}{2} \left(1+\sqrt{5}+\sqrt{2 \left(5+\sqrt{5}\right)}\right),\frac{1}{2} \left(1+\sqrt{5}+\sqrt{2
   \left(5+\sqrt{5}\right)}\right),1+\sqrt{5+2 \sqrt{5}},2+\sqrt{5}+\sqrt{5+2 \sqrt{5}}\right\}
   $
   
   }

 \noindent The global dimension $\vert {\mathcal E}_7 \vert$ already obtained directly at the beginning can be recovered by summing their squares.
 As a check, by taking only squares of dimensions relative to modular vertices, we can  recalculate in a third way the ambichiral dimension $\vert J \vert = 4$.
 
\subsection{$B_2$ at level $12$ and its exceptional quantum subgroup ${\mathcal E}_{12}$}
\subsubsection*{$B_2$ at level $12$}

\begin{figure}
\centerline{\scalebox{1.0}{\includegraphics[width=12cm]{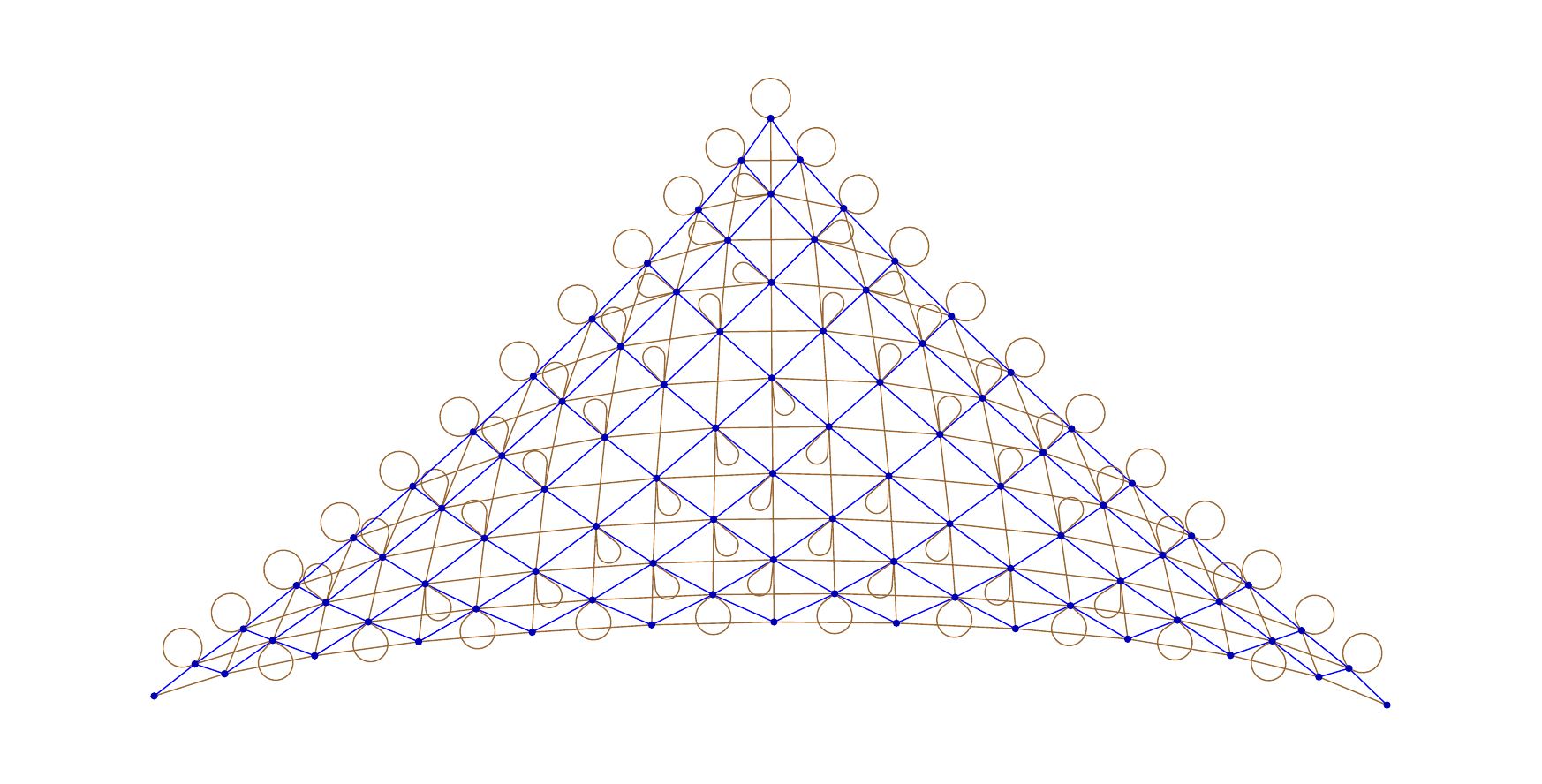}}} 
\caption{The fusion graph of  $\mathcal{A}_{12}(B_2))$.}
\label{fig:B2level12}
\end{figure}

The category ${\mathcal A}_{12}={\mathcal A}_{12}(B_2)$ has $91$ simple objects ($r_A=91$) and, like before,  it is unreasonable to print lists giving quantum dimensions, conformal weights or the blocks for the associated quantum groupo\"\i d (see nevertheless footnote \ref{webpub}). This is even more so for a matrix like $S$, which is of dimension $r_A \times r_A$ and which is not sparse!. The calculation can nevertheless be done, not only for $S$ and $T$, but for all fusion matrices. The fundamental objects have dimensions
$qdim(\{0,1\}) = \frac{1}{2} \left(\sqrt{5}+\sqrt{15+6 \sqrt{5}}\right)$, $dim(\{1,0\}) =\frac{1}{4} \left(5+3 \sqrt{5}+\sqrt{6 \left(5+\sqrt{5}\right)}\right)$, 
the central charge is $c=8$, the bialgebra dimension $d_{\mathcal B}= 1918226582 = 2^1 \, 10067^1 \, 95273^1$, and $$\vert{\mathcal A}_{12}\vert= 450 \left(262+117 \sqrt{5}+\sqrt{136995+61266 \sqrt{5}}\right)$$
 The graph \ref{fig:B2level12} describes the fundamental fusion matrices $N_{\{0,1\}}$ and $N_{\{1,0\}}$.
We did not print the labels for vertices since the pattern should by now be obvious.

\subsubsection*{The exceptional quantum subgroup ${\mathcal E}_{12}$}

We use  the  embedding  $B_2 \subset E_8$, with Dynkin index $k=12$, which is conformal when the level of $B_2 \simeq Spin(5)$ is $k=12$ and the level of $E_8$ is $1$ (same conformal charges $c$).
The only integrable representation of $E_8$ at level $1$ is trivial. A priori the type $I$ partition function obtained by reduction will have only one block.
The trivial representation of $E_8$ has a (quantum) dimension equal to $1$ (note that $q=exp(i\pi/(30+1)$ for ${\mathcal A}_1(E_8)$), so that the ambichiral dimension is $\vert J \vert = 1$, and has conformal weight  ${0}$. Comparing this value with those obtained (modulo $1$) for $B_2$ at level $12$ (\ie using $q=exp(i\pi/(3+12)$) gives a necessary condition for branching rules. One obtains the exceptional partition function :
$$Z= (\{0,0\}+\{0,6\}+\{1,8\}+\{3,2\}+\{3,8\}+2 \{4,4\}+\{6,6\}+\{7,2\}+\{12,0\})^2$$
Notice that one of the $9$ terms has multiplicity $2$.
One verifies that $M$ (a matrix of size $91\times 91$) defined by $Z = \sum_{m,n} m \, M_{mn} \, n$ indeed commutes with $S$ and $T$.

The direct sum of the simple objects appearing in the first (and unique) modular block of $Z$ define a Frobenius  algebra  ${\mathcal F}$, and by summing their dimensions, one finds  $\vert{\mathcal F}\vert = \vert{\mathcal A}_{12}(B_2)\vert / \vert{\mathcal E}_{12}\vert = 15 \left(12+5 \sqrt{5}+\sqrt{3 \left(85+38 \sqrt{5}\right)}\right)$.
Since $\vert{\mathcal A}_{12}\vert$ was obtained in the previous subsection, we calculate immediately  $\vert{\mathcal E}_{12}\vert$ and find that in this very particular case 
 $\vert{\mathcal F}\vert=\vert{\mathcal E}\vert=\sqrt{\vert{\mathcal A}\vert}$, something which was expected from the general relation
 $\vert{\mathcal A}_k(G)\vert  = \vert{\mathcal E}_k\vert^2 / \vert J \vert $ since here  $\vert J \vert =1$. 
The number of simple objects for ${\mathcal A}_{12}(B_2)$ is $r_A=91$.
From the modular invariant matrix we read the generalized exponents and obtain $r_E=12$ (number of simple objects $a$ of the quantum subgroup), 
$r_O=144$ (number of quantum symmetries $x$), and $r_W=81$ (number of independent toric matrices $W_x$). 

\begin{figure}[htp]
  \begin{center}
    \scalebox{1.4}{ \includegraphics[width=8cm]{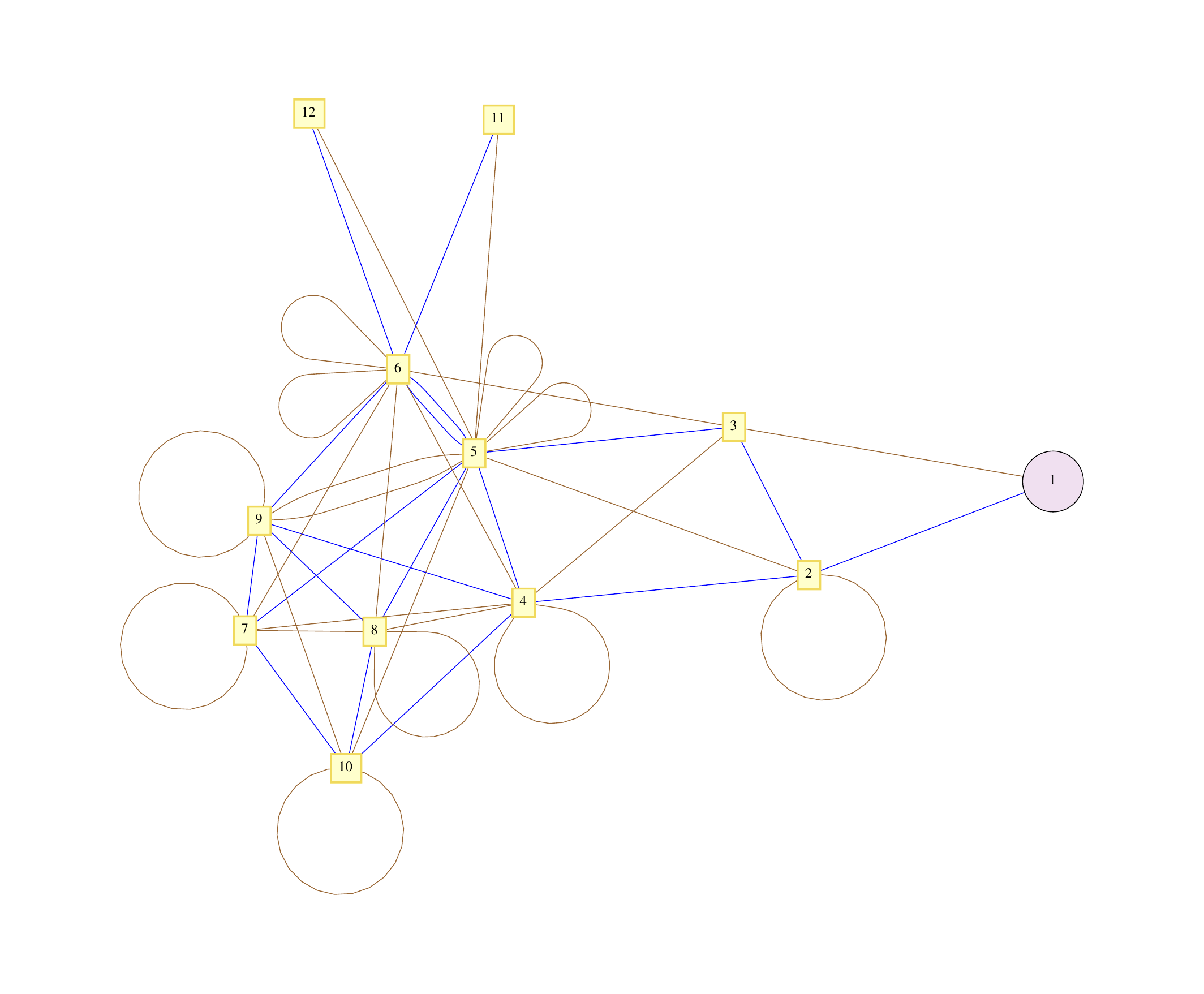}}
  \end{center}
  \caption{Fusion graph for the quantum subgroup $\mathcal{E}_{12}(B_2)$ }
  \label{fig:fusiongraphE12B2}
\end{figure}

Resolution of the modular splitting equation, leads to the list of toric matrices $W_x$. 
There are $304$ possible norms associated with matrices ${\mathcal K}_{mn}$; many  independent $W_x$ can be determined after analyzing the first few norms $1\ldots 5$ but one has to go up to norm $15$ to conclude, and the discussion is quite involved. After some efforts,  one obtains  $81$ independent matrices $W_x$ of size $91\times 91$, the difficulty, then, is actually to determine the multiplicities of each of them within the whole family of $144$ elements. This can be done.
Generators for the quantum symmetries ${\mathcal O(E)}$ are then obtained as usual by solving intertwining equations. We checked that the obtained solution for the full algebra of quantum symmetries could be obtained by using directly the ``reduced modular splitting technique'' described in the general section. We also determined the 
fusion graph of ${\mathcal E}_{12}$, not only by restriction of the algebra of quantum symmetries (fundamental annular matrices, which are adjacency matrices for the fusion graphs of ${\mathcal E}_{12}$, are obtained as diagonal sub-blocks of the quantum symmetry generators)
but also by using the simpler chiral modular splitting technique. The full graph of quantum symmetries cannot be displayed, because of its size, but its general features are a priori clear, since we already know that the ambichiral algebra has only one element (the origin) :  the full graph contains twelve copies ($12 \times 12$ vertices) of the quantum subgroup ${\mathcal E}_{12}$ which is described on figure \ref{fig:fusiongraphE12B2} with the same conventions as usual ( the blue color corresponds to $\{0,1\}$ and the brown to \{1,0\}).
It enjoys self-fusion, of course, its unit is in position ``1'', which is also the unique modular points corresponding to the unique ambichiral vertex ( intersection of chiral parts of the graph of quantum symmetries).

The other annular matrices (both for the quantum subgroup and its module) are obtained from the fundamental ones by using fusion polynomials, a large table involving powers up to $12$ in $f_2$ and $f_3$.
 Using these matrices, we have calculated (but we shall not list) the horizontal $91$ horizontal dimensions $d_n$.  The dimension of the associated quantum groupo\"\i d $d_{\mathcal B}=d_n.d_n=  51068=2^2 \, 17^1 \, 751^1$. 
The intertwiner $E_{\underline{0}}$  for $\mathcal{E}_{12}$ is a rectangular matrix of dimension $91 \times 12$, a  bit large to be displayed;  unfortunately  we cannot even condense the information listing its non-zero entries (induction)\footnote{Entries in column $1$ are known since they coincide with the (unique) modular blocks of $Z$.}. 
 The quantum dimensions of the simple objects of ${\mathcal E}_{12}$ can now be calculated, as before,  in several ways.
 Their explicit expression as particular roots of fourth degree polynomials is not particularly enlightening, but the sum of their squares simplifies to give the global dimension $\vert\mathcal{E}_{12}\vert = \frac{15}{2} \left(2+\sqrt{5}+\sqrt{15+6 \sqrt{5}}\right)^2$ an expression that looks formally different from the expression directly obtained at the very beginning, from conformal embedding considerations. The two expressions are nevertheless equal since $2 \sqrt{15+6 \sqrt{5}}+\sqrt{75+30 \sqrt{5}}=\sqrt{3 \left(85+38
   \sqrt{5}\right)}$.

\section{Quantum subgroups of $G_2$}

\subsection{General properties of  $G_2$}

$G_2 $ has rank $r= 2$, dual Coxeter number $g=4$, Coxeter number $\gamma =6 $, dimension $d= r(\gamma+1) =14$, adjacency matrix (Dynkin) $G = \{\{0, 3\}, \{1, 0\}\}$ (it is not simply laced), 
Cartan matrix $A = 2 * \one - G = \{\{2,-3\},\{-1,2\}\}$, quadratic form matrix $Q =\{\{2,1\},\{1,2/3\}\} = A^{-1}.\{\{1,0\},\{0,1/3\}\}$, highest root $\theta = \{1, 0\}$,  Weyl vector $\varrho = \{1,1\}$, exponents $(\epsilon_1=1,\epsilon_2=5)$, Casimir degrees $(\epsilon_1+1=2,\epsilon_2+1=6)$, Weyl group order $(\epsilon_1+1)(2,\epsilon_2+1)=12$ (a dihedral group),  number of positive roots ${\Sigma_{+}} =6$,  $\Delta =1/3$.
The central charge of $B_2$ at level $k$ is $c=\frac{d \, k}{g+k}=\frac{14 k}{k+4}$.
The level $\langle \lambda,  \theta  \rangle$ of an irreducible representation with highest weight $\{\lambda_1, \lambda_2\}$ is $2 \lambda_1+\lambda_2$ and its (classical or quantum) dimension is given by the classical or quantum version of the Weyl formula. It is
$$
\frac{\left(\lambda _1+1\right)_q \left(\frac{\lambda
   _2}{3}+\frac{1}{3}\right)_q \left(\lambda _1+\frac{\lambda
   _2}{3}+\frac{4}{3}\right)_q \left(\lambda _1+\frac{2 \lambda
   _2}{3}+\frac{5}{3}\right)_q \left(\lambda _1+\lambda _2+2\right)_q
   \left(2 \lambda _1+\lambda _2+3\right)_q}{\left(\frac{1}{3}\right)_q
   1_q \left(\frac{4}{3}\right)_q \left(\frac{5}{3}\right)_q 2_q 3_q}
   $$
   where, as usual, $(x)_q=\frac{q^x - q^{-x}} {q - q^{-1}}$ and $(x)_q=x$ if $q=1$.
   At level $k$ one takes $q = exp(i \pi/\kappa)$ with an altitude $\kappa=k+g=k+4$.
   The  fundamental representations\footnote{Warning: some people would write $\{\lambda_2, \lambda_1\}$ instead of $\{\lambda_1, \lambda_2\}$ and/or would shift indices by $\{1,1\}$}  are of highest weight $\{1,0\}$, classical dimension $14$,  quantum dimension $\frac{\left(\frac{7}{3}\right)_q \left(\frac{8}{3}\right)_q 5_q}{1_q
   \left(\frac{4}{3}\right)_q \left(\frac{5}{3}\right)_q}$, and of  highest weight  $\{0,1\}$, classical dimension $7$, quantum dimension $\frac{\left(\frac{2}{3}\right)_q \left(\frac{7}{3}\right)_q
   4_q}{\left(\frac{1}{3}\right)_q \left(\frac{4}{3}\right)_q 2_q}$.

\subsection{$G_2$ at level $1$}

The category ${\mathcal A}_1={\mathcal A}_1(G_2)$ has two simple objects ($r_A=2$) corresponding to the trivial $(\{0,0\}$ and the basic representation of highest weight $\{0,1\}$, of classical dimension $7$, and quantum dimension $\phi = \frac{1}{2} \left(1+\sqrt{5}\right)$ (the golden number). One finds
 $\vert{\mathcal A}_1\vert=1+\phi^2=\frac{1}{2} \left(5+\sqrt{5}\right)$. 
 The conformal weights modulo $1$ are $(0, 2/5)$ and the central charge is $c=14/5$. 
  The two fundamental representations of $G_2$ do not appear at level $1$.
 The modular generators, with $S_{{\0}{\0}} = 1/\sqrt{\vert{\mathcal A}_1\vert} = \sqrt{\frac{1}{10} \left(5-\sqrt{5}\right)}$, are
$$
\begin{array}{ccc}
S/S_{{\0}{\0}}=
\left(
\begin{array}{cc}
 1 & \phi  \\
 \phi  & -1
\end{array}
\right)
& , & T=
\left(
\begin{array}{cc}
 e^{-\frac{7 i \pi }{30}} & 0 \\
 0 & e^{\frac{17 i \pi }{30}}
\end{array}
\right)
\end{array}
$$

\begin{figure}
\centerline{\scalebox{1.2}{\includegraphics[width=6cm]{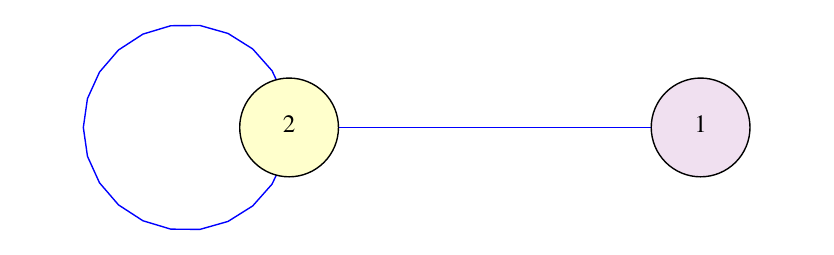}}} 
\caption{The fusion graph of  $\mathcal{A}_1(G_2)$.}
\label{fig:G2level1}
\end{figure}

The graph \ref{fig:G2level1}  describes only the basic fusion matrix $N_{\{0,1\}}$ since the other fundamental representation $\{1,0\}$ appears only at level 2. The blocks of the associated quantum groupo\"\i d are therefore of size $d_n=(2, 3)$ so that the dimension of the later is $d_{\mathcal B}=13$. 
Notice that such a tadpole graph would be forbidden as a graph describing a quantum module associated with $SU(2)$  at level two (same altitude), but we are here discussing  ${\mathcal A}_1(G_2)$.  Here and below, labels $1,2,3,4,\ldots$ of vertices refer to the order\footnote{The level dependent ordering is not trivial for $G_2$} that was chosen on the set of simple objects $\{0,0\},\{0,1\},\{0,2\},\{1,0\},\ldots$. Multiplication by the generator $\{0,1\}$, labelled $2$,  is encoded by unoriented blue edges. Like for $B_2$, edges are not oriented because conjugation is trivial. Multiplication by the fundamental $\{1,0\}$ labelled $4$ (warning: not $3$)  will be encoded  by unoriented brown edges when studying cases  $k>1$ (it does not exist when $k=1$).  
The fusion graphs of ${\mathcal A}_k(G_2)$ obtained for the first few levels are rather simple, and one may notice that, when $k$ increases, the fusion graph associated with the basic representation $\{0,1\}$ grows in a predictable way:  it is a truncated version of the corresponding classical ${\mathcal A}_\infty(G_2)$ fusion graph. This is not so for the other fundamental representation $\{1,0\}$. 
We shall see that in our cases it is  anyway enough to know $N_{\{0,1\}}$ only since $N_{\{1,0\}}$  can be obtained from the other by a fusion polynomial.

\subsection{$G_2$ at level $2$}

The category ${\mathcal A}_2={\mathcal A}_2(G_2)$ has four simple objects ($r_A=4$).
The following table gives  the chosen ordering for the highest weights, the quantum dimensions, the conformal weights modulo $1$, and the size of  blocks for the associated quantum groupo\"\i d.

$$
\begin{array}{cccc}
 \{0,0\} & \{0,1\} & \{0,2\} & \{1,0\} \\
 z_1 & z_2 & z_3 & z_4 \\
 0 & \frac{1}{3} & \frac{7}{9} & \frac{2}{3} \\
 4 & 10 & 9 & 7
\end{array}
$$
where $z_i$ are the largest (real) roots of the polynomials
${1, 1 - 3 z^2 + z^3 ,3 - 3 z^2 + z^3, -1 - 3 z + z^3}$.
 The central charge is  $c=14/3$, the bialgebra dimension is $d_{\mathcal B}=246=2^1\, 3^1 \ 41^1$, and $\vert{\mathcal A}_2\vert=  \vert{\mathcal A}_2(G_2)\vert$ is the largest (real) root of the polynomial $-243 + 162 z - 27 z^2 + z^3$.
 
\begin{figure}
\centerline{\scalebox{1.2}{\includegraphics[width=6cm]{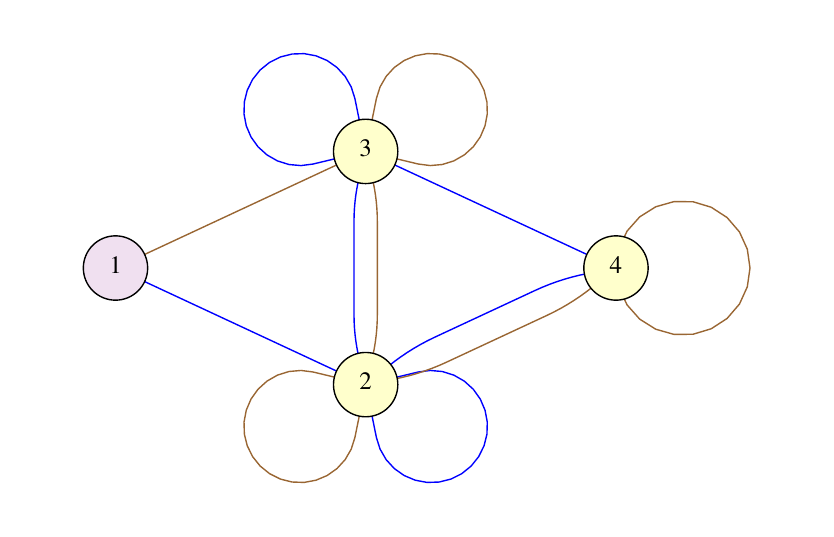}}} 
\caption{The fusion graph of  $\mathcal{A}_2(G_2)$.}
\label{fig:G2level2}
\end{figure}

 The  modular generators are as follows :
\begin{eqnarray*}
S &=&
\tiny{
\left(
\begin{array}{cccc}
 \frac{2}{3} \sin \left(\frac{\pi }{9}\right) & \frac{2}{3} \cos
   \left(\frac{\pi }{18}\right) & \frac{1}{\sqrt{3}} & \frac{2}{3} \sin
   \left(\frac{2 \pi }{9}\right) \\
 \frac{2}{3} \cos \left(\frac{\pi }{18}\right) & \frac{2}{3} \sin
   \left(\frac{2 \pi }{9}\right) & -\frac{1}{\sqrt{3}} & -\frac{2}{3}
   \sin \left(\frac{\pi }{9}\right) \\
 \frac{1}{\sqrt{3}} & -\frac{1}{\sqrt{3}} & 0 & \frac{1}{\sqrt{3}} \\
 \frac{2}{3} \sin \left(\frac{2 \pi }{9}\right) & -\frac{2}{3} \sin
   \left(\frac{\pi }{9}\right) & \frac{1}{\sqrt{3}} & -\frac{2}{3} \cos
   \left(\frac{\pi }{18}\right)
\end{array}
\right)
}
\\
T&=&diag 
\left\{e^{-\frac{7 i \pi }{18}},e^{\frac{5 i \pi }{18}},e^{-\frac{5 i \pi }{6}},e^{\frac{17 i \pi }{18}}\right\}
\end{eqnarray*}
Notice that $S_{{\0}{\0}} = 1/\sqrt{\vert{\mathcal A}_2\vert} = \frac{2}{3} \sin \left(\frac{\pi }{9}\right)$.
The graph \ref{fig:G2level2} describes the fundamental fusion matrices $N_{\{0,1\}}$ and $N_{\{1,0\}}$.

\subsection{$G_2$ at level $3$ and its exceptional quantum subgroup ${\mathcal E}_3$}

\subsubsection*{$G_2$ at level $3$}

The category ${\mathcal A}_3={\mathcal A}_3(G_2)$ has six simple objects ($r_A=6$).
The following table gives  the chosen ordering for the highest weights, the quantum dimensions, the conformal weights modulo $1$, and the size of  blocks for the associated quantum groupo\"\i d.

$$
\begin{array}{cccccc}
 \{0,0\} & \{0,1\} & \{0,2\} & \{1,0\} & \{0,3\} & \{1,1\} \\
 1 & \frac{1}{2} \left(3+\sqrt{21}\right) & \frac{1}{2}
   \left(7+\sqrt{21}\right) & \frac{1}{2} \left(3+\sqrt{21}\right) &
   \frac{1}{2} \left(3+\sqrt{21}\right) & \frac{1}{2}
   \left(5+\sqrt{21}\right) \\
 0 & \frac{2}{7} & \frac{2}{3} & \frac{4}{7} & \frac{1}{7} & 0 \\
 6 & 20 & 30 & 20 & 20 & 25
\end{array}
$$

 The central charge is  $c=6$, the bialgebra dimension is $d_{\mathcal B}=2761=11^1 \, 251^1$, and $\vert{\mathcal A}_3\vert=  \vert{\mathcal A}_3(G_2)\vert = \frac{21}{2} \left(5+\sqrt{21}\right)$. 

 The  modular generators are as follows :
\begin{eqnarray*}
S &=&
\left(
\begin{array}{cccccc}
 \sqrt{\frac{1}{42} \left(5-\sqrt{21}\right)} & \frac{1}{\sqrt{7}} & \frac{1}{\sqrt{3}} & \frac{1}{\sqrt{7}} &
   \frac{1}{\sqrt{7}} & \sqrt{\frac{1}{42} \left(5+\sqrt{21}\right)} \\
 \frac{1}{\sqrt{7}} & \frac{{u_3}}{\sqrt{7}} & 0 & \frac{{u_2}}{\sqrt{7}} & \frac{{u_1}}{\sqrt{7}} &
   -\frac{1}{\sqrt{7}} \\
 \frac{1}{\sqrt{3}} & 0 & -\frac{1}{\sqrt{3}} & 0 & 0 & \frac{1}{\sqrt{3}} \\
 \frac{1}{\sqrt{7}} & \frac{{u_2}}{\sqrt{7}} & 0 & \frac{{u_1}}{\sqrt{7}} & \frac{{u_3}}{\sqrt{7}} &
   -\frac{1}{\sqrt{7}} \\
 \frac{1}{\sqrt{7}} & \frac{{u_1}}{\sqrt{7}} & 0 & \frac{{u_3}}{\sqrt{7}} & \frac{{u_2}}{\sqrt{7}} &
   -\frac{1}{\sqrt{7}} \\
 \sqrt{\frac{1}{42} \left(5+\sqrt{21}\right)} & -\frac{1}{\sqrt{7}} & \frac{1}{\sqrt{3}} & -\frac{1}{\sqrt{7}} &
   -\frac{1}{\sqrt{7}} & \sqrt{\frac{1}{42} \left(5-\sqrt{21}\right)}
\end{array}
\right)
\\
T &=& diag 
\left\{-i,e^{\frac{i \pi }{14}},e^{\frac{5 i \pi }{6}},e^{\frac{9 i \pi }{14}},e^{-\frac{3 i \pi }{14}},-i\right\}
\end{eqnarray*}
where $\{u_1,u_2,u_3\}$ are the three real roots of the polynomial $u^3-u^2-2 u+1$, ordered from the smallest to the largest. Notice that $S_{{\0}{\0}} = 1/\sqrt{\vert{\mathcal A}_3\vert} =\sqrt{\frac{1}{42} \left(5-\sqrt{21}\right)}$. 

\begin{figure}
\centerline{\scalebox{1.1}{\includegraphics[width=6cm]{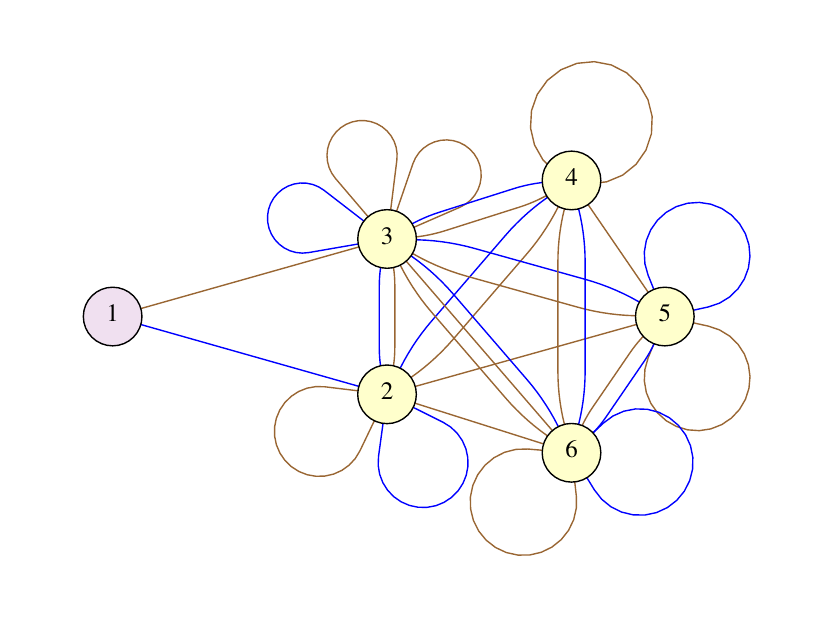}}} 
\caption{The fusion graph of  $\mathcal{A}_3(G_2))$.}
\label{fig:G2level3}
\end{figure}

The graph \ref{fig:G2level3} describes the fundamental fusion matrices $N_{\{0,1\}}$ and $N_{\{1,0\}}$.

\subsubsection*{The exceptional quantum subgroup ${\mathcal E}_3$}

We use  the  embedding  $G_2 \subset E_6$ (with Dynkin index $k=3$), which is conformal when the level of $G_2$ is $k=3$ and the level of $E_6$ is $1$ (same conformal charges $c$).

The integrable representations of $E_6$ at level $1$ are the trivial representation and the two basic representations $\{1,0,0,0,0;0\}$ and $\{0,0,0,0,1;0\}$ of classical dimensions $27$ (notice that not all fundamental representations appear).
A priori the type $I$ partition function obtained by reduction will have three modular blocks (two of them should be equal, by classical symmetry).
The above three simple objects of $E_6$ have a (quantum) dimension equal to $1$ (using $q=exp(i\pi/(12+1)$ for ${\mathcal A}_1(E_6)$), so that summing their square gives the ambichiral dimension $\vert J \vert = 3$.
Their conformal weights are $\{0, 2/3, 2/3\}$. Comparing them with those obtained (modulo $1$) for $G_2$ at level $3$ (\ie using $q=exp(i\pi/(3+3)$) gives a necessary condition for branching rules. One obtains in this way one exceptional $G_2$ partition function (quadratic form):
$$Z= (\{0,0\}+\{1,1\})^2+ 2 \{0,2\}^2 ,  \qquad 
M=
\left(
\begin{array}{cccccc}
 1 & . & . & . & . & 1 \\
 . & . & . & . & . & . \\
 . & . & 2 & . & . & . \\
 . & . & . & . & . & . \\
 . & . & . & . & . & . \\
 1 & . & . & . & . & 1
\end{array}
\right)
$$
The modular invariant matrix $M$ defined by $Z = \sum_{m,n} m \, M_{mn} \, n$ indeed commutes with $S$ and $T$.

The simple objects appearing in the first modular block of $Z$  define a Frobenius  algebra  ${\mathcal F} =
\{0,0\}  \oplus  \{1,1\}$, with $\vert{\mathcal F}\vert = \vert{\mathcal A}_3(G_2)\vert / \vert{\mathcal E}_3\vert = qdim(\{0,0\} )+qdim(\{1,1\} )= 1+\frac{1}{2} \left(5+\sqrt{21}\right)$. Since $\vert{\mathcal A}_3\vert = \frac{21}{2} \left(5+\sqrt{21}\right)$, we find $\vert{\mathcal E}_3\vert =\frac{3}{2} \left(7+\sqrt{21}\right)$. 
The general relation $\vert{\mathcal A}_k(G)\vert  = (\vert{\mathcal E}_k\vert)^2 / \vert J \vert $ leads again to $\vert J \vert =3$. The number of simple objects for ${\mathcal A}_3(G_2)$ is $r_A=6$.
From the modular invariant matrix  we read the generalized exponents and obtain $r_E=4$ (number of simple objects $a$ of the quantum subgroup), 
$r_O=8$ (number of quantum symmetries $x$), $r_W=5$ (number of independent toric matrices $W_x$). 

\begin{figure}[htp]
  \begin{center}
    \subfigure[The left and right generators $\{0,1\}$]{\label{fig:OLR01E3G2}\includegraphics[width=7.0cm]{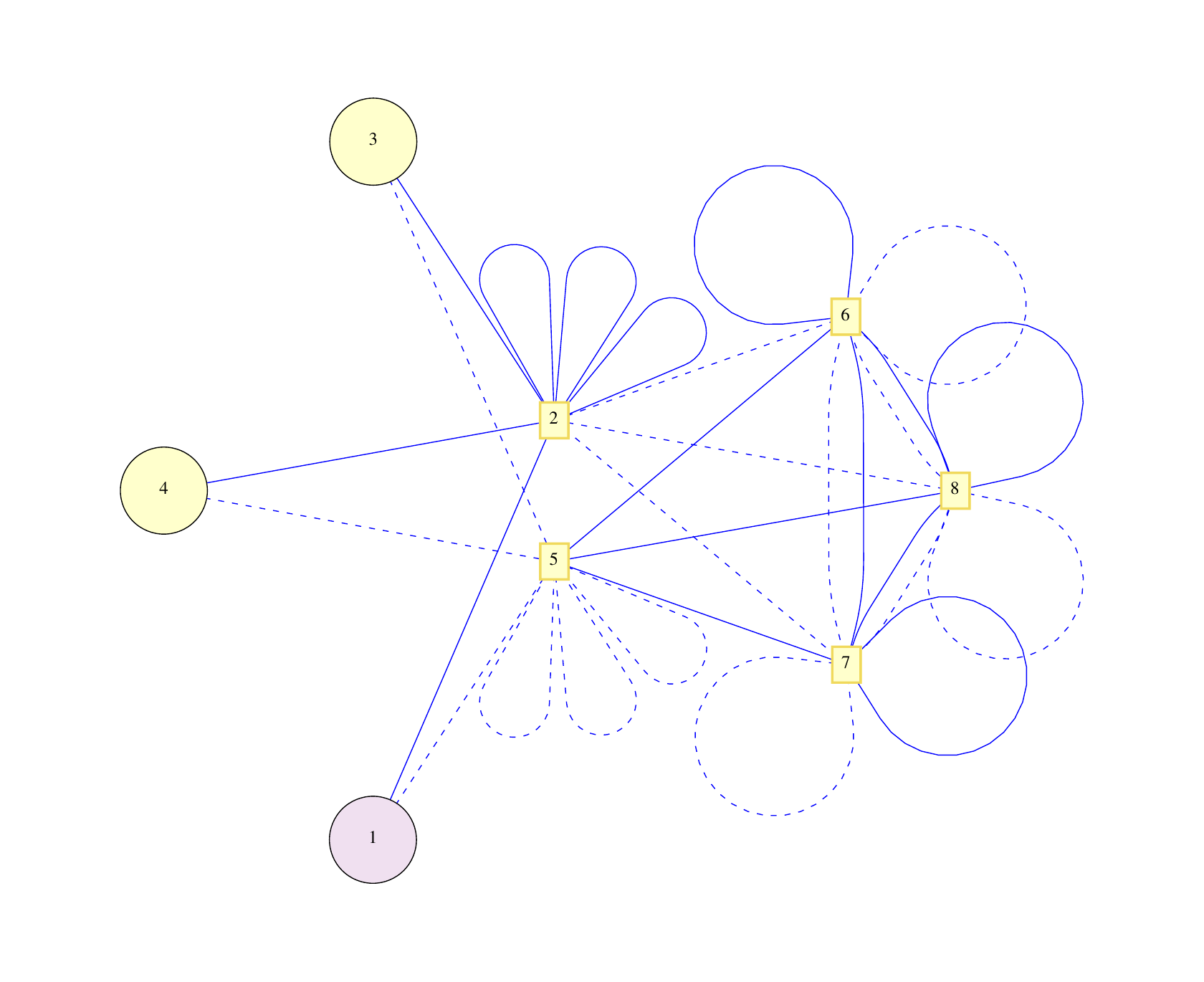}}
    \subfigure[The left and right generators $\{1,0\}$]{\label{fig:OLR10E3G2}\includegraphics[width=7.0cm]{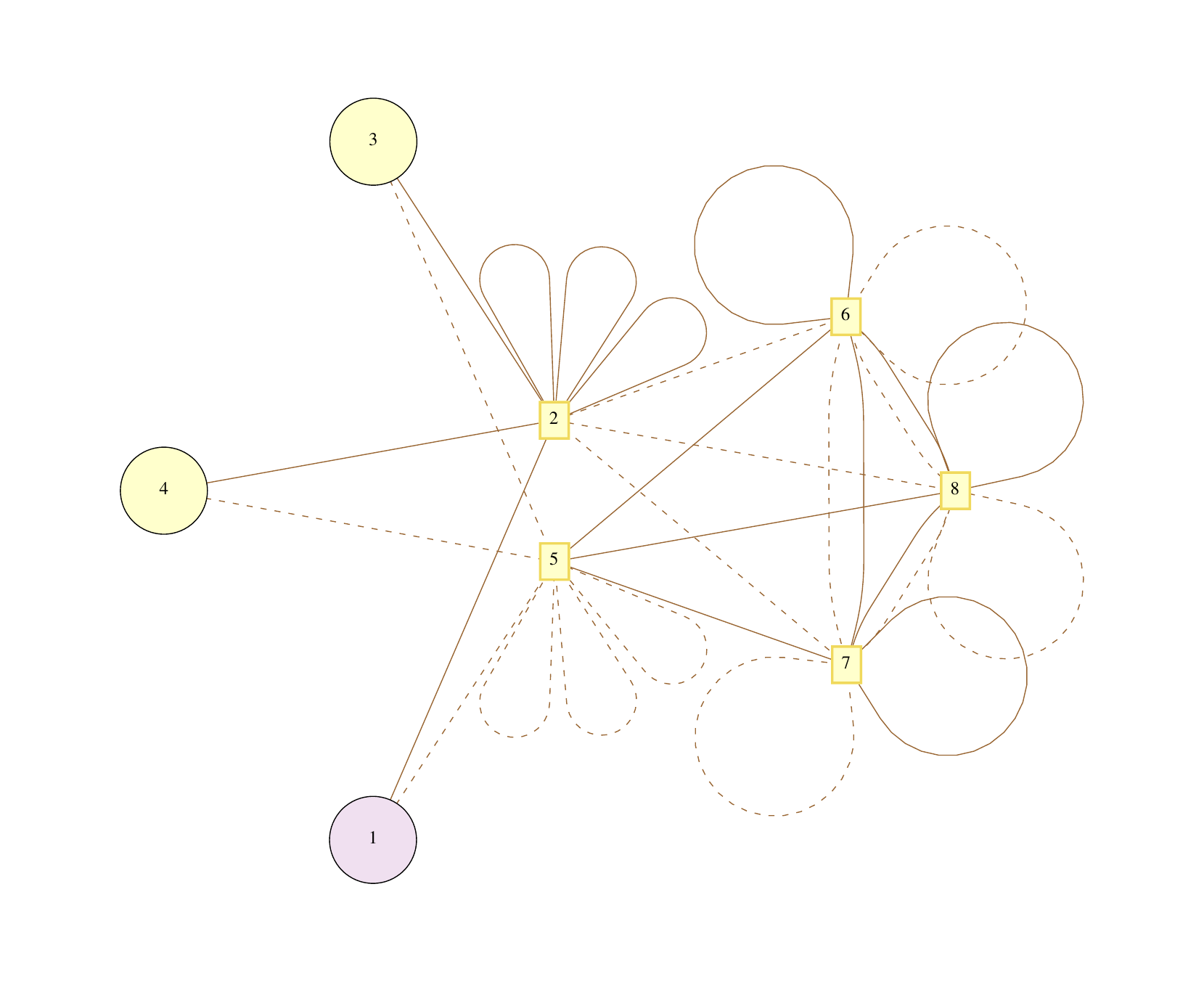}}
  \end{center}
  \caption{Quantum symmetries of  $\mathcal{E}_3(G_2)$}
\end{figure}

Resolution of the modular splitting equation, leads to the list of toric matrices $W_x$. 
One has to analyze $8$ possible norms associated with matrices ${\mathcal K}_{mn}$.  Three toric matrices appear in norm $1$ and two other independent ones in norm $3$, so we reach the total $r_W=5$, but the fourth one has to appear twice, and the last one, three times, in the final list, so we reach the total $r_O=8$.  The first is the modular invariant and was already given, the four others $W_x$'s  are as follows:

{\tiny
$$
\left(
\begin{array}{cccccc}
 . & 1 & 1 & 1 & 1 & 1 \\
 . & . & . & . & . & . \\
 . & 2 & 2 & 2 & 2 & 2 \\
 . & . & . & . & . & . \\
 . & . & . & . & . & . \\
 . & 1 & 1 & 1 & 1 & 1
\end{array}
\right),\left(
\begin{array}{cccccc}
 . & . & . & . & . & . \\
 1 & . & 2 & . & . & 1 \\
 1 & . & 2 & . & . & 1 \\
 1 & . & 2 & . & . & 1 \\
 1 & . & 2 & . & . & 1 \\
 1 & . & 2 & . & . & 1
\end{array}
\right),\left(
\begin{array}{cccccc}
 . & . & 1 & . & . & . \\
 . & . & . & . & . & . \\
 1 & . & 1 & . & . & 1 \\
 . & . & . & . & . & . \\
 . & . & . & . & . & . \\
 . & . & 1 & . & . & .
\end{array}
\right),\left(
\begin{array}{cccccc}
 . & . & . & . & . & . \\
 . & 1 & 1 & 1 & 1 & 1 \\
 . & 1 & 1 & 1 & 1 & 1 \\
 . & 1 & 1 & 1 & 1 & 1 \\
 . & 1 & 1 & 1 & 1 & 1 \\
 . & 1 & 1 & 1 & 1 & 1
\end{array}
\right)
$$
}
Generators for quantum symmetries ${\mathcal O(E)}$ are  obtained from the toric matrices by solving intertwining equations.
Again, we display separately the graphs associated with $\{0,1\}$  (the  two chiral generators $V_{01,00}$ (left blue: plain lines) and $V_{00,01}$ (right blue: dashed lines)), and those associated with $\{1,0\}$ (the  two chiral generators $V_{10,00}$  (left brown: plain lines) and $V_{00,10}$ (right brown: dashed lines)), on figures \ref{fig:OLR01E3G2}, \ref{fig:OLR10E3G2}.

\begin{figure}[htp]
  \begin{center}
    \subfigure[$\mathcal{E}_3(G_2)$]{\label{fig:E3G2}\scalebox{1.4}{ \includegraphics[width=6cm]{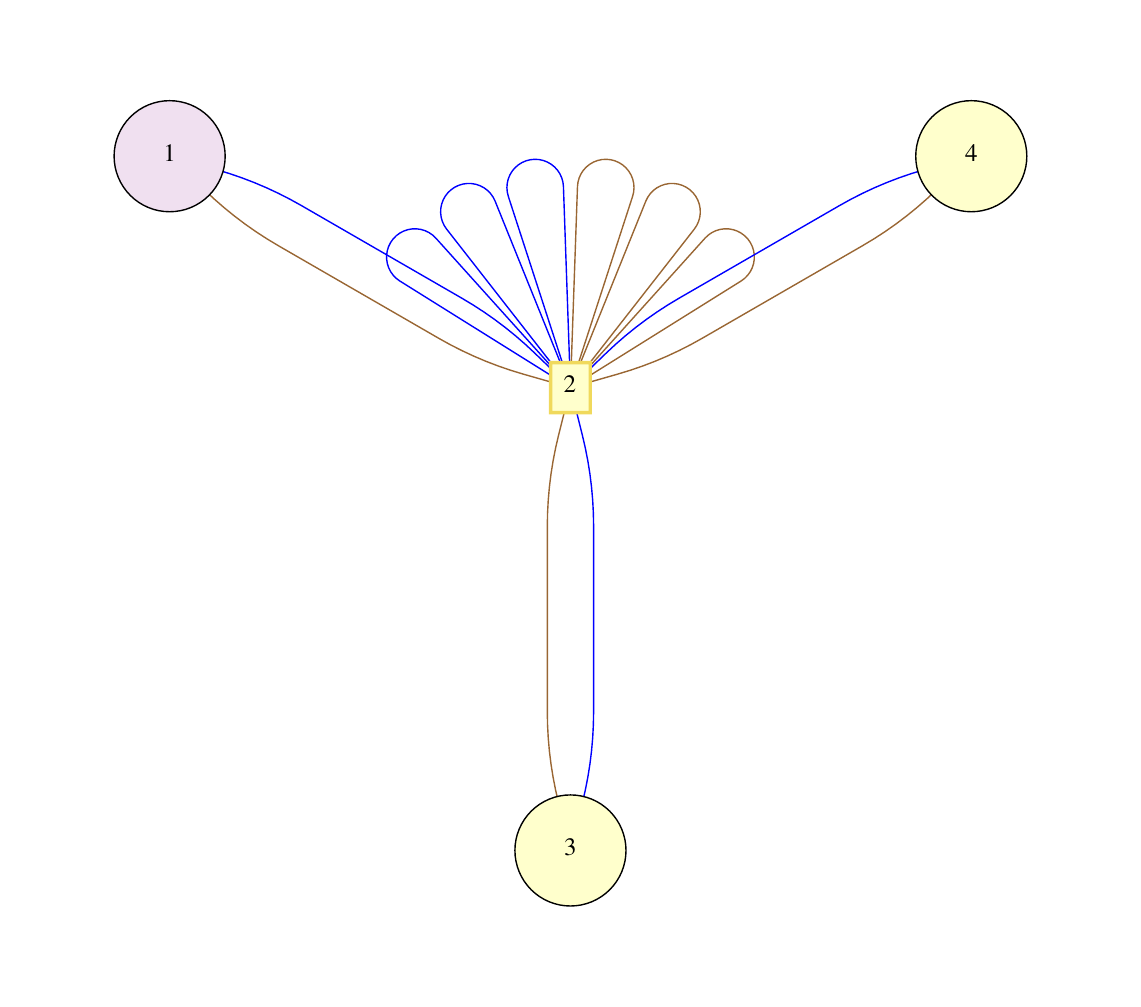}}}
    \subfigure[${\mathcal{E}_3}^M(G_2)$]{\label{fig:E3mG2}\includegraphics[width=6cm]{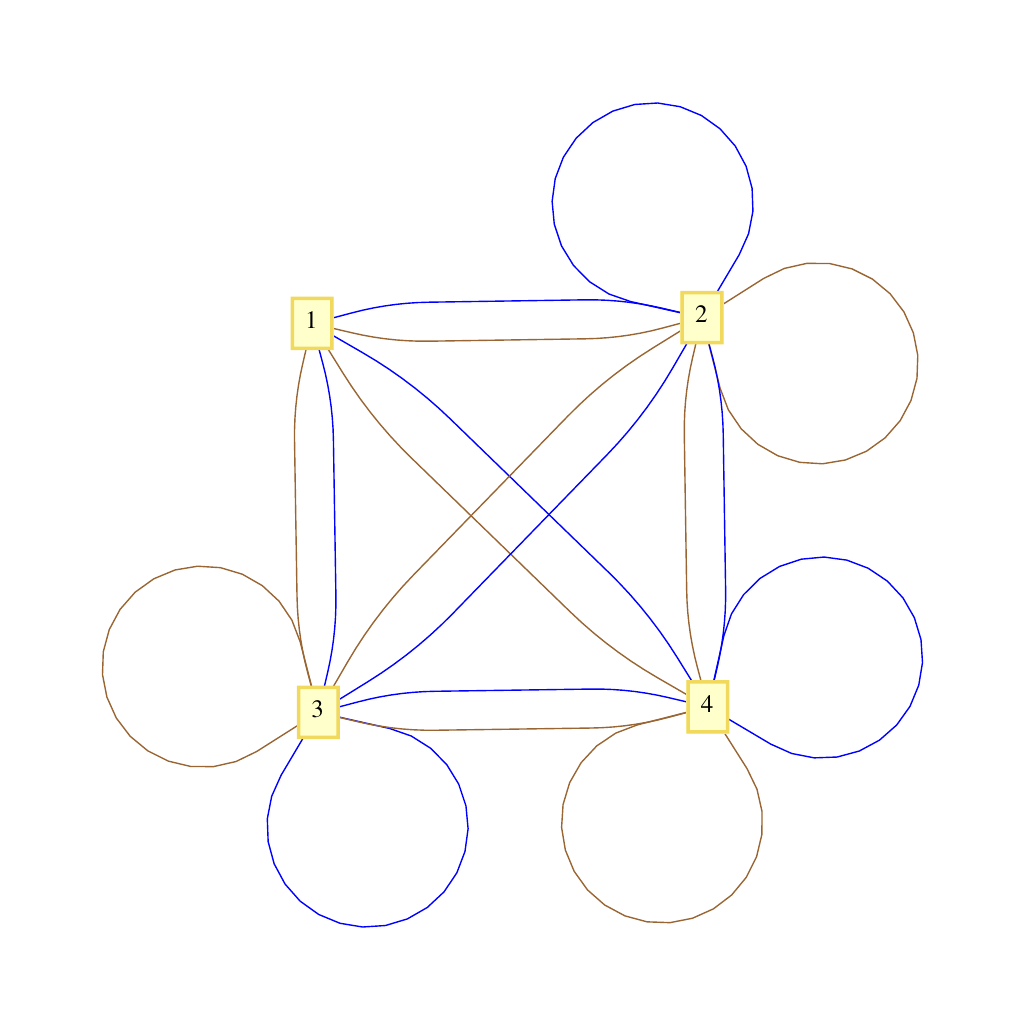}}
  \end{center}
  \caption{Fusion graphs for the quantum subgroup $\mathcal{E}_3(G_2)$ and its module ${\mathcal{E}_3}^M(G_2)$ }
  \label{fig:fusiongraphsE3G2}
\end{figure}

From the generators (graphs) of quantum symmetries, we discover the generators (annular matrices $F_{\{0,1\}}$ and $F_{\{1,0\}}$) for the quantum subgroup $\mathcal{E}_3(G_2)$ itself, and discover as well the existence of a module denoted ${{\mathcal E}_3}^M$. Their graphs are given on figures \ref{fig:E3G2}  and \ref{fig:E3mG2}. As usual the blue color corresponds to $\{0,1\}$ and the brown to \{1,0\}.
${\mathcal E}_3$ enjoys self-fusion (its unit is in position ``1'') and its three modular points  corresponding to the three ambichiral points obtained as intersection of chiral parts of the graph of quantum symmetries, are located in positions $\{1,3,4\}$. The graphs for generators $\{0,1\}$ and $\{1,0\}$ turn out to be identical, and this is certainly not a generic feature for other values of the level $k$. We shall explain this coincidence below.
 So, $F_{\{0,1\}}$ is equal to $F_{\{1,0\}}$ for this particular value of $k=3$.  Actually $F_{\{0,1\}}=F_{\{1,0\}}=F_{\{0,3\}}$, with
$$
F_{\{0,1\}}=
\left(
\begin{array}{cccc}
 0 & 1 & 0 & 0 \\
 1 & 3 & 1 & 1 \\
 0 & 1 & 0 & 0 \\
 0 & 1 & 0 & 0
\end{array}
\right)
\quad
F_{\{0,2\}}=
\left(
\begin{array}{cccc}
 0 & 1 & 1 & 1 \\
 1 & 5 & 1 & 1 \\
 1 & 1 & 0 & 1 \\
 1 & 1 & 1 & 0
\end{array}
\right)
\quad
F_{\{1,1\}}=\left(
\begin{array}{cccc}
 1 & 1 & 0 & 0 \\
 1 & 4 & 1 & 1 \\
 0 & 1 & 1 & 0 \\
 0 & 1 & 0 & 1
\end{array}
\right)
$$
The other annular matrices are obtained, for instance, from the fusion polynomials :

$
f_1=1, \quad  f_2= \text{x}, \quad
f_3=\frac{1}{7} \left(4 \text{x}^5-17 \text{x}^4-2
   \text{x}^3+37 \text{x}^2+\text{x}-7\right)$
   
   $
f_4=     \frac{1}{7} \left(-4
   \text{x}^5+17 \text{x}^4+2 \text{x}^3-30 \text{x}^2-8
   \text{x}\right)$
   
 $
  f_5=   \frac{1}{7} \left(2 \text{x}^5-5 \text{x}^4-15
   \text{x}^3+15 \text{x}^2+25 \text{x}\right),
   \;
f_6=   \frac{1}{7} \left(-3
   \text{x}^5+11 \text{x}^4+12 \text{x}^3-33 \text{x}^2-20
   \text{x}+7\right)
$

\noindent
We can obtain annular matrices for the module in a similar fashion.
Notice that in the present case one can express all fusion matrices, including the fundamental $N_{\{1,0\}}$, in terms of the basic one $N_{\{0,1\}}$. In particular the former is given in terms of the later by the polynomial $f_4$. One can see that $N_{\{1,0\}}\neq N_{\{0,1\}}$ and 
$N_{\{1,0\}}\neq {N_{\{0,1\}}}^{tr}$. However, the minimal polynomial of $F_{\{0,1\}}$ is $x^3-3 x^2 -3x$ and the reminder of the polynomial quotient of $f_4$ by the later is just $x$. For this reason $F_{\{1,0\}} = F_{\{0,1\}}$. 
The first essential matrix for the quantum subgroup (intertwiners) and its module are given below:

$$
E_{\underline{0}}=
\left(
\begin{array}{cccc}
 1 & 0 & 0 & 0 \\
 0 & 1 & 0 & 0 \\
 0 & 1 & 1 & 1 \\
 0 & 1 & 0 & 0 \\
 0 & 1 & 0 & 0 \\
 1 & 1 & 0 & 0
\end{array}
\right)
 \qquad
{E_{\underline{0}}}^M=
\left(
\begin{array}{cccc}
 1 & 0 & 0 & 0 \\
 0 & 1 & 1 & 1 \\
 2 & 1 & 1 & 1 \\
 0 & 1 & 1 & 1 \\
 0 & 1 & 1 & 1 \\
 1 & 1 & 1 & 1
\end{array}
\right)$$
The horizontal dimensions of ${\mathcal E}_3$, are $d_n= \{4, 9, 17, 9, 9, 13\}$ so that $d_\mathcal{B}=717=3^1 \, 239^1$.
For the module ${{\mathcal E}_3}^M$ we find $d_n=\{{4, 15, 23, 15, 15, 19}\}$ and $d_\mathcal{B}=1581=3^1 \, 17^1 \, 31^1$.
The q-dimensions of simple objects of ${\mathcal E}_3$,  calculated from the normalized Perron-Frobenius vector are $\left\{1,\frac{1}{2} \left(3+\sqrt{21}\right),1,1\right\}$, so that by summing squares we recover the global dimension $\vert {\mathcal E}_3 \vert = 3/2 (7 + Sqrt[21])$ already obtained.
For the module, we find q-dimensions, $\left\{1,\frac{1}{6} \left(3+\sqrt{21}\right),\frac{1}{6} \left(3+\sqrt{21}\right),\frac{1}{6} \left(3+\sqrt{21}\right)\right\}$ and $\vert{{\mathcal E}_3}^M\vert = \frac{1}{2} \left(7+\sqrt{21}\right)$.

\subsection{$G_2$ at level $4$ and its exceptional quantum subgroup ${\mathcal E}_4$}

\subsubsection*{$G_2$ at level $4$}

The category ${\mathcal A}_4={\mathcal A}_4(G_2)$ has nine simple objects ($r_A=9$).
The following table gives  the chosen ordering for the highest weights, the quantum dimensions, the conformal weights modulo $1$, and the size of  blocks for the associated quantum groupo\"\i d.

$$
\begin{array}{ccccccccc}
\{0 , 0\} & \{0 , 1\} & \{0 , 2\} & \{1 , 0\} & \{0 , 3\} & \{1 , 1\} & \{0 , 4\} & \{1 , 2\} & \{2 , 0\} \\
1 & 2+\sqrt{6} & 2 \left(2+\sqrt{6}\right) & 3+\sqrt{6} & 5+2 \sqrt{6} & 2
   \left(3+\sqrt{6}\right) & 3+\sqrt{6} & 2
   \left(2+\sqrt{6}\right) & 2+\sqrt{6}\\
   0 & \frac{1}{4} & \frac{7}{12} & \frac{1}{2} & 0 & \frac{7}{8} & \frac{1}{2} & \frac{1}{3} & \frac{1}{4} \\
9 & 34 & 66 & 41& 73 & 80& 41& 66 & 34
\end{array}
$$

\begin{figure}[htp]
  \begin{center}
    \subfigure[$\mathcal{A}_4(G_2)$, rep. $\{0,1\}$]{\label{fig:G2rep01level4}\scalebox{1.4}{ \includegraphics[width=6cm]{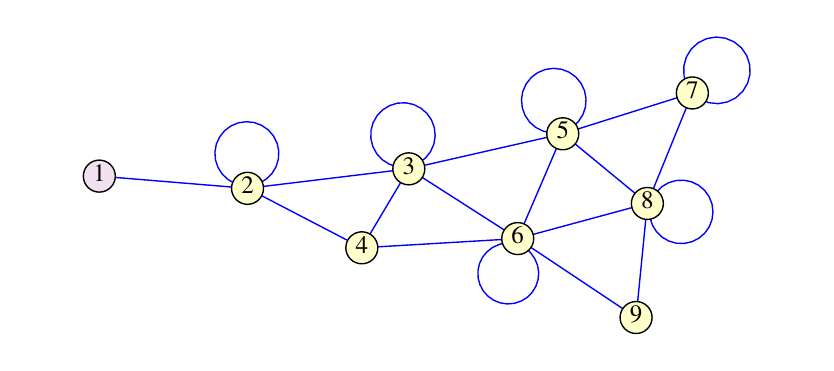}}}
    \subfigure[$\mathcal{A}_4(G_2)$, rep. $\{1,0\}$]{\label{fig:G2rep10level4}\includegraphics[width=6cm]{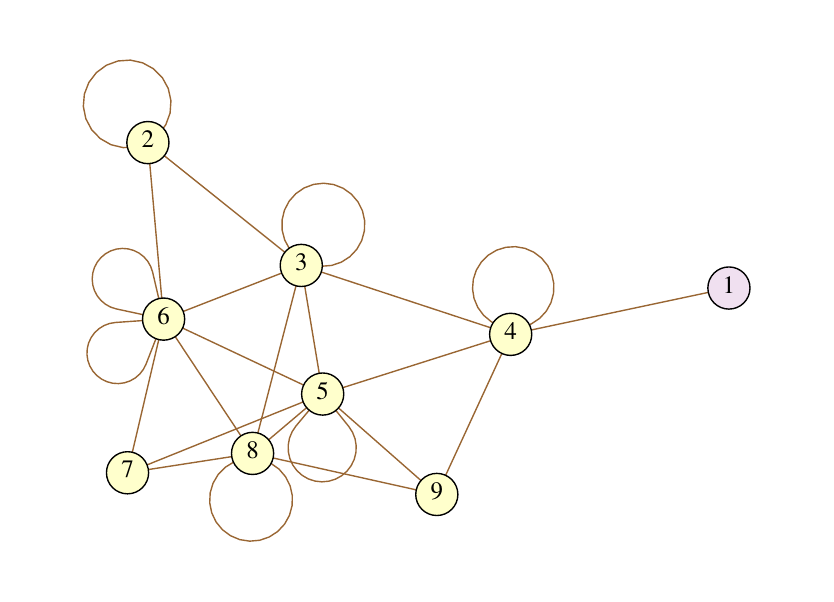}}
  \end{center}
  \caption{Fusion graphs for  $\mathcal{A}_4(G_2)$}
  \label{fig:G2level4}
\end{figure}

 The central charge is  $c=7$, the bialgebra dimension is $d_{\mathcal B}=26196=2^2 3^1 37^1 59^1$, and $\vert{\mathcal A}_4\vert=  \vert{\mathcal A}_4(G_2)\vert =48 \left(5+2 \sqrt{6}\right)$.  The  modular generators are as follows :
\begin{eqnarray*}
S &=&
\left(
\begin{array}{ccccccccc}
 \frac{1}{12}   \left(3-\sqrt{6}\right) & \frac{1}{2 \sqrt{6}} & \frac{1}{\sqrt{6}} &
   \frac{1}{4} & \frac{1}{12} \left(3+\sqrt{6}\right) & \frac{1}{2} &
   \frac{1}{4} & \frac{1}{\sqrt{6}} & \frac{1}{2 \sqrt{6}} \\
 \frac{1}{2 \sqrt{6}} & \frac{1+\sqrt{3}}{2 \sqrt{6}} &
   \frac{1}{\sqrt{6}} & \frac{1}{2 \sqrt{2}} & -\frac{1}{2 \sqrt{6}} & 0
   & -\frac{1}{2 \sqrt{2}} & -\frac{1}{\sqrt{6}} & \frac{1- \sqrt{3}}{2  \sqrt{6}} \\
 \frac{1}{\sqrt{6}} & \frac{1}{\sqrt{6}} & -\frac{1}{\sqrt{6}} & 0 &
   -\frac{1}{\sqrt{6}} & 0 & 0 & \frac{1}{\sqrt{6}} & \frac{1}{\sqrt{6}}
   \\
 \frac{1}{4} & \frac{1}{2 \sqrt{2}} & 0 & \frac{1}{4}
   \left(1-\sqrt{2}\right) & \frac{1}{4} & -\frac{1}{2} & \frac{1}{4}
   \left(1+\sqrt{2}\right) & 0 & -\frac{1}{2 \sqrt{2}} \\
 \frac{1}{12} \left(3+\sqrt{6}\right) & -\frac{1}{2 \sqrt{6}} &
   -\frac{1}{\sqrt{6}} & \frac{1}{4} & \frac{1}{12}
   \left(3-\sqrt{6}\right) & \frac{1}{2} & \frac{1}{4} &
   -\frac{1}{\sqrt{6}} & -\frac{1}{2 \sqrt{6}} \\
 \frac{1}{2} & 0 & 0 & -\frac{1}{2} & \frac{1}{2} & 0 & -\frac{1}{2} & 0
   & 0 \\
 \frac{1}{4} & -\frac{1}{2 \sqrt{2}} & 0 & \frac{1}{4}
   \left(1+\sqrt{2}\right) & \frac{1}{4} & -\frac{1}{2} & \frac{1}{4}
   \left(1-\sqrt{2}\right) & 0 & \frac{1}{2 \sqrt{2}} \\
 \frac{1}{\sqrt{6}} & -\frac{1}{\sqrt{6}} & \frac{1}{\sqrt{6}} & 0 &
   -\frac{1}{\sqrt{6}} & 0 & 0 & \frac{1}{\sqrt{6}} & -\frac{1}{\sqrt{6}}
   \\
 \frac{1}{2 \sqrt{6}} & \frac{1- \sqrt{3}}{2  \sqrt{6}}&
   \frac{1}{\sqrt{6}} & -\frac{1}{2 \sqrt{2}} & -\frac{1}{2 \sqrt{6}} & 0
   & \frac{1}{2 \sqrt{2}} & -\frac{1}{\sqrt{6}} & \frac{1+\sqrt{3}}{2
   \sqrt{6}}
\end{array}
\right)
\\
T &=& diag 
\left\{e^{-\frac{7 i \pi }{12}},e^{-\frac{i \pi }{12}},e^{\frac{7 i \pi
   }{12}},e^{\frac{5 i \pi }{12}},e^{-\frac{7 i \pi }{12}},e^{-\frac{5 i
   \pi }{6}},e^{\frac{5 i \pi }{12}},e^{\frac{i \pi }{12}},e^{-\frac{i
   \pi }{12}}\right\}
\end{eqnarray*}
 In particular one checks that $S_{{\0}{\0}} = 1/\sqrt{\vert{\mathcal A}_4\vert} =  \frac{1}{12} \left(3-\sqrt{6}\right)$. \\
  \noindent
The graphs \ref{fig:G2level4} describe (separately) the fundamental fusion matrices $N_{\{0,1\}}$ and $N_{\{1,0\}}$.

\subsubsection*{The exceptional quantum subgroup ${\mathcal E}_4$}

We use  the  embedding  $G_2 \subset D_7=Spin(14)$ (with Dynkin index $k=4$), which is conformal when the level of $G_2$ is $k=4$ and the level of $D_7$ is $1$ (same conformal charges $c$). This is the adjoint embedding : $dim(G_2)=14$ and the level $k$ coincides with  the dual Coxeter number $G_2$.

The integrable representations of $D_7$ at level $1$ are the trivial, the vectorial and the two spinorial  representations.
A priori the type $I$ partition function obtained by reduction will have four modular blocks (two of them being  equal).
The above four simple objects of $D_7$ have a (quantum) dimension equal to $1$ (using $q=exp(i\pi/(12+1)$ for ${\mathcal A}_1(D_7)$), so that summing their square gives the ambichiral dimension $\vert J \vert = 4$.
Their conformal weights are $\{0, 7/8, 7/8,1/2\}$. Comparing them with those obtained (modulo $1$) for $G_2$ at level $4$ (\ie using $q=exp(i\pi/(4+4)$) gives a necessary condition for branching rules. One obtains in this way one exceptional $G_2$ partition function (quadratic form):
$$Z=(\{0,0\}+\{0,3\})^2+(\{0,4\}+\{1,0\})^2+2 \{1,1\}^2, 
\quad
M=
\left(
\begin{array}{ccccccccc}
 1 & . & . & . & 1 & . & . & . & . \\
 . & . & . & . & . & . & . & . & . \\
 . & . & . & . & . & . & . & . & . \\
 . & . & . & 1 & . & . & 1 & . & . \\
 1 & . & . & . & 1 & . & . & . & . \\
 . & . & . & . & . & 2 & . & . & . \\
 . & . & . & 1 & . & . & 1 & . & . \\
 . & . & . & . & . & . & . & . & . \\
 . & . & . & . & . & . & . & . & .
\end{array}
\right)
$$
One verifies that the modular invariant matrix $M$ defined by $Z = \sum_{m,n} m \, M_{mn} \, n$ indeed commutes with $S$ and $T$.

The simple objects appearing in the first modular block of $Z$  define a Frobenius  algebra  ${\mathcal F} =
\{0,0\}  \oplus  \{0,3\}$, with $\vert{\mathcal F}\vert = \vert{\mathcal A}_3(B_2)\vert / \vert{\mathcal E}_3\vert = qdim(\{0,0\} )+qdim(\{0,3\} )= 6+2 \sqrt{6}$. Since $\vert{\mathcal A}_4\vert =  48 \left(5+2 \sqrt{6}\right)$, we find $\vert{\mathcal E}_4\vert =8 \left(3+\sqrt{6}\right) $. 
The general relation $\vert{\mathcal A}_k(G)\vert  = (\vert{\mathcal E}_k\vert)^2 / \vert J \vert $ leads again to the result  $\vert J \vert =4$. The number of simple objects for ${\mathcal A}_4(G_2)$ is $r_A=9$.
From the modular invariant matrix we read the generalized exponents and obtain $r_E=6$ (number of simple objects $a$ of the quantum subgroup), 
$r_O=12$ (number of quantum symmetries $x$), $r_W=9$ (number of independent toric matrices $W_x$). 

Resolution of the modular splitting equation, leads to the list of toric matrices $W_x$. 
One has to analyze $14$ possible norms associated with matrices ${\mathcal K}_{mn}$.  
Five toric matrices are found in norm $1$, three others in norm $2$ (two of them appear twice), and the last, which also appear twice, in norm $3$ is  half the corresponding ${\mathcal K}_{mn}$, which has even coefficients. So we recover $r_W=5+3+1$ and $r_O=5+(1+2\times 2)+ 2\times2$, \ie $12$ toric matrices ($9$ independent ones) of dimension $9 \times 9$ that we cannot give here (the first is the modular invariant, as usual).

\begin{figure}[htp]
  \begin{center}
 \includegraphics[width=10cm]{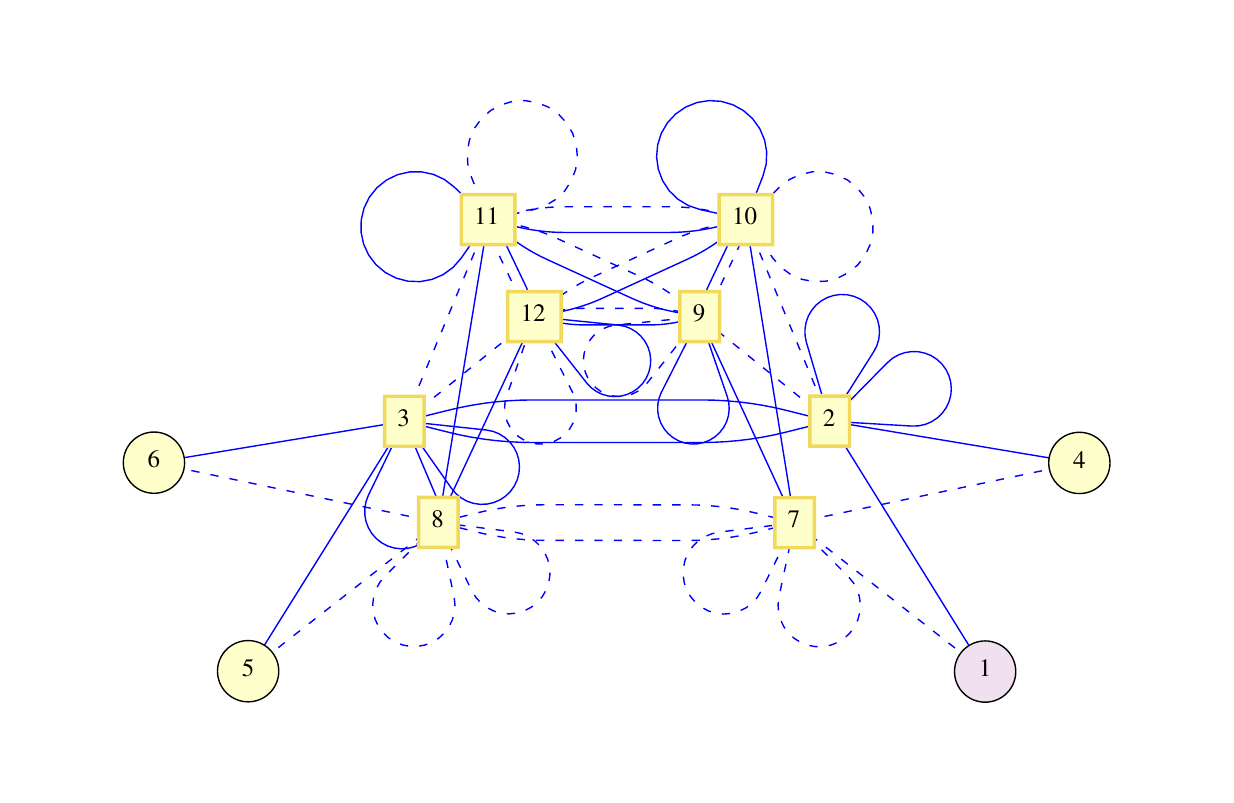}
    \end{center}
  \caption{Quantum symmetries (left and right) of  $\mathcal{E}_4(G_2)$ for the basic representation $\{0,1\}$}
  \label{fig:OLR01E4G2}
\end{figure}

Generators for quantum symmetries are obtained from the toric matrices by solving intertwining equations. 
We give only, on figure \ref{fig:OLR01E4G2}, the graph of quantum symmetries for the left and right generators $V_{\{0,1\},\{0,0\}}$ and $V_{\{0,0\},\{0,1\}}$ corresponding to the basic representation $\{0,1\}$. The graph corresponding to the other fundamental representation $\{1,0\}$ is more cluttered, and is anyway not needed because matrices (fusion $N_{\{1,0\}}$, annular $F_{\{1,0\}}$ or double fusion $V_{\{1,0\}}$) can be expressed polynomially in terms of their  $\{1,0\}$ analogs with a fusion polynomial\footnote{Fusion polynomials have usually rather complicated rational coefficients, but matrix elements of $N_n,F_n,V_{m,n}$ are non negative integers, as it should be.}:  $\{1,0\}=f(\{0,1\})$ with 
$$f(x)=\frac{1}{40} \left(-39 \text{x}^8+255 \text{x}^7-179 \text{x}^6-1059
   \text{x}^5+810 \text{x}^4+1352 \text{x}^3-552 \text{x}^2-548
   \text{x}-40\right)$$

\begin{figure}[htp]
  \begin{center}
    \subfigure[$\mathcal{E}_4(G_2)$]{\label{fig:E4G2}\scalebox{1.4}{ \includegraphics[width=5.5cm]{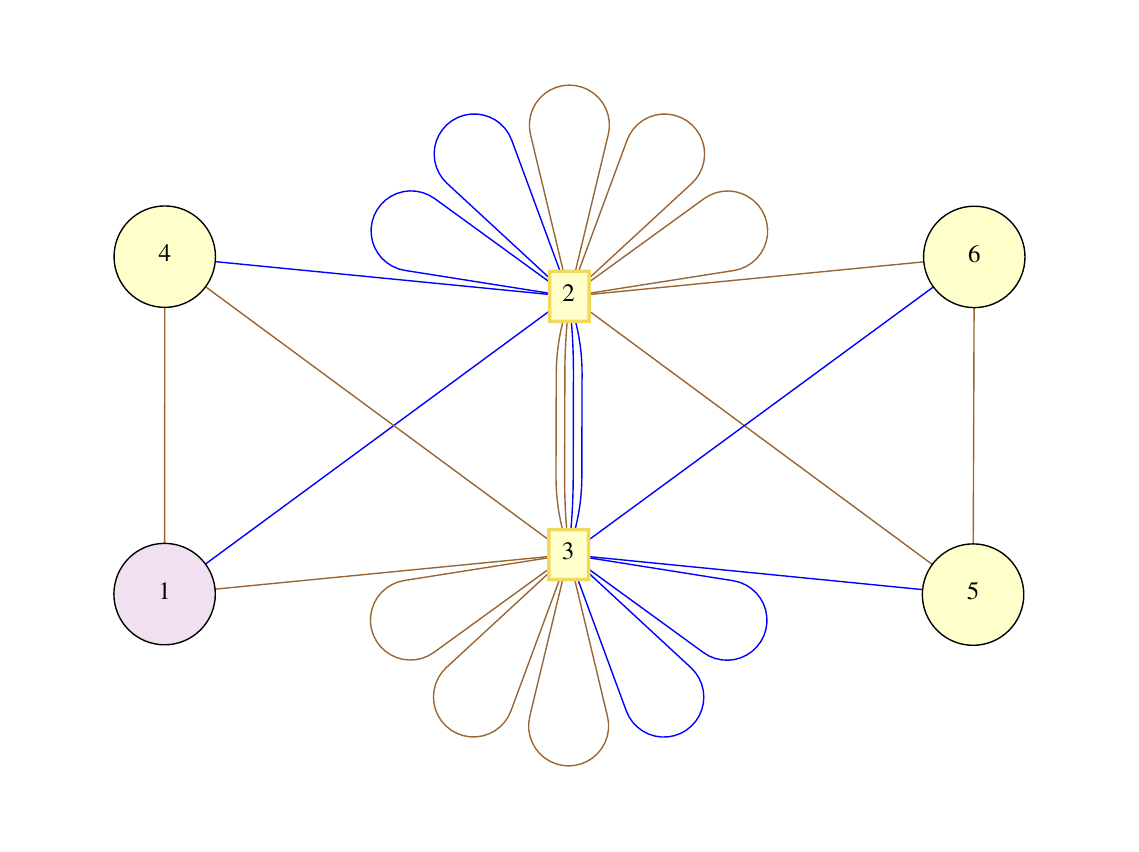}}}
    \subfigure[${\mathcal{E}_4}^M(G_2)$]{\label{fig:E4mG2}\includegraphics[width=6.cm]{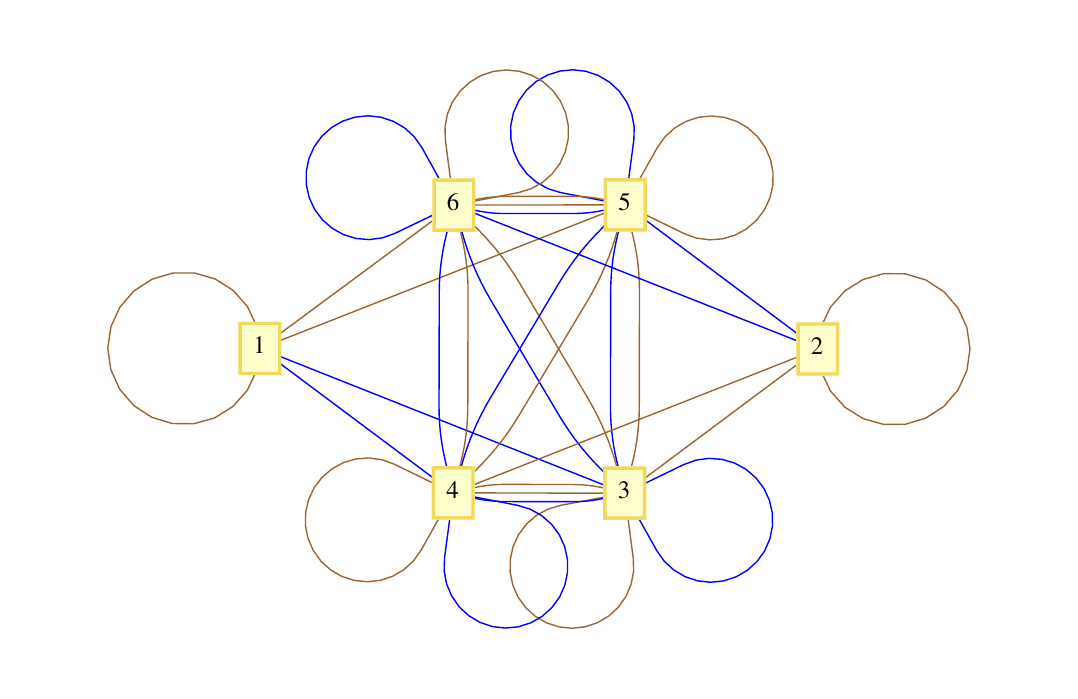}}
  \end{center}
  \caption{Fusion graphs for the quantum subgroup $\mathcal{E}_4(G_2)$ and its module ${\mathcal{E}_4}^M(G_2)$ }
  \label{fig:fusiongraphsE4G2}
\end{figure}

Intersection of the left and right connected components of the unit (marked ``1'') define the ambichiral vertices $\{1,4,5,6\}$.
The left (or right) graph has two non isomorphic connected components, so that we discover both a quantum subgroup $\mathcal{E}_4$ and a module ${\mathcal{E}_4}^M$. Adjacency matrices of their graphs (given on figure \ref{fig:fusiongraphsE4G2}) are the annular matrices $F_{\{0,1\}}$ and ${F_{\{0,1\}}}^M$ extracted from the double fusion generator $V_{\{0,1\},\{0,0\}}$. We display both actions of  $\{0,1\}$ (blue lines) and  $\{1,0\}$ (brown lines), although the later could be obtained from the former by using the previously given fusion polynomial.
The nine annular matrices, both for the quantum subgroup and its module, are calculated with fusion polynomials.
We only give $F_{\{0,1\}}$ and ${F_{\{0,1\}}}^M$, as well as the first essential matrix,  for both.

$$
F_{\{0,1\}} =
\left(
\begin{array}{cccccc}
 0 & 1 & 0 & 0 & 0 & 0 \\
 1 & 2 & 2 & 1 & 0 & 0 \\
 0 & 2 & 2 & 0 & 1 & 1 \\
 0 & 1 & 0 & 0 & 0 & 0 \\
 0 & 0 & 1 & 0 & 0 & 0 \\
 0 & 0 & 1 & 0 & 0 & 0
\end{array}
\right)
\qquad
{F_{\{0,1\}}}^M=
\left(
\begin{array}{cccccc}
 0 & 0 & 1 & 1 & 0 & 0 \\
 0 & 0 & 0 & 0 & 1 & 1 \\
 1 & 0 & 1 & 1 & 1 & 1 \\
 1 & 0 & 1 & 1 & 1 & 1 \\
 0 & 1 & 1 & 1 & 1 & 1 \\
 0 & 1 & 1 & 1 & 1 & 1
\end{array}
\right)
$$

$$
E_{\underline{0}}=
\left(
\begin{array}{cccccc}
 1 & 0 & 0 & 0 & 0 & 0 \\
 0 & 1 & 0 & 0 & 0 & 0 \\
 0 & 1 & 1 & 0 & 0 & 0 \\
 0 & 0 & 1 & 1 & 0 & 0 \\
 1 & 1 & 1 & 0 & 0 & 0 \\
 0 & 1 & 1 & 0 & 1 & 1 \\
 0 & 1 & 0 & 1 & 0 & 0 \\
 0 & 1 & 1 & 0 & 0 & 0 \\
 0 & 0 & 1 & 0 & 0 & 0
\end{array}
\right)
 \qquad
{E_{\underline{0}}}^M=
\left(
\begin{array}{cccccc}
 1 & 0 & 0 & 0 & 0 & 0 \\
 0 & 0 & 1 & 1 & 0 & 0 \\
 0 & 0 & 1 & 1 & 1 & 1 \\
 1 & 0 & 0 & 0 & 1 & 1 \\
 1 & 0 & 1 & 1 & 1 & 1 \\
 0 & 2 & 1 & 1 & 1 & 1 \\
 1 & 0 & 1 & 1 & 0 & 0 \\
 0 & 0 & 1 & 1 & 1 & 1 \\
 0 & 0 & 0 & 0 & 1 & 1
\end{array}
\right)
$$
The horizontal dimensions of ${\mathcal E}_4$, are $d_n= \{6, 16, 32, 22, 38, 44, 22, 32, 16\}$ so that $d_\mathcal{B}=6944=2^5 \, 7^1\, 31^1$.
For the module ${{\mathcal E}_4}^M$ we find $d_n=\{{6, 24, 48, 30, 54, 60, 30, 48, 24}\}$ and $d_\mathcal{B}=14112 = 2^5\, 3^2 \, 7^2$.
The q-dimensions of simple objects of ${\mathcal E}_4$,  calculated from the normalized Perron-Frobenius vector are $\left\{1,2+\sqrt{6},2+\sqrt{6},1,1,1\right\}$, so that by summing squares we recover the global dimension $\vert {\mathcal E}_4 \vert = 8 \left(3+\sqrt{6}\right)$ already obtained.
For the module, we find q-dimensions, $\left\{1,1,1+\sqrt{\frac{3}{2}},1+\sqrt{\frac{3}{2}},1+\sqrt{\frac{3}{2}}
   ,1+\sqrt{\frac{3}{2}}\right\}$ and $\vert{{\mathcal E}_4}^M\vert = 4 \left(3+\sqrt{6}\right)$.
   Going back to the graph of ${\mathcal E}_4$, we can check the known value of $\vert J \vert=4$ by summing the squares of quantum dimensions relative to the modular vertices $\{1,4,5,6\}$ corresponding to the ambichiral vertices of the graph of quantum symmetries.

\subsection*{Acknowledgements}
This research was supported in part by the CNRS-CNRST, under convention M. L. N 1055/09. The authors R. R. and E. H. T. are also supported by U. R. A. C. 07.


\end{document}